\documentclass[10pt, leqno]{article}
\usepackage{amsmath}
\usepackage{amsfonts}
\usepackage{amssymb,esint}
\usepackage{graphicx}
\usepackage{color, colordvi}

\newtheorem{thm}{Theorem}[section]
\newtheorem{cor}[thm]{Corollary}
\newtheorem{lem}[thm]{Lemma}
\newtheorem{prop}[thm]{Proposition}
\newtheorem{defn}[thm]{Definition}

\numberwithin{equation}{section}

\newcommand{\dx}{\,{\rm d}x}
\newcommand{\dy}{\,{\rm d}y}

\newcommand{\rd}{{\rm d}}

\def\LL{\mathrm{L}} 
\def\supp{\mathrm{supp}} 
\newcommand{\RR}{\mathbb{R}}

\def\ee{\mathrm{e}} 
\def\dist{\mathrm{dist}} 
\def\diam{\mathrm{diam}} 

\def\qed{\unskip\kern 6pt \penalty 500
\raise -2pt\hbox{\vrule \vbox to8pt{\hrule width 6pt
\vfill\hrule}\vrule}\par}
\definecolor{darkblue}{rgb}{0.05, .05, .65}
\definecolor{darkgreen}{rgb}{0.05, .55, .05}
\definecolor{darkred}{rgb}{0.8,0,0}
\begin{document}
\title{\bf  Quantitative Local Bounds for Subcritical\\ Semilinear Elliptic Equations}
\author{Matteo Bonforte$^{\,a}$
Gabriele Grillo$^{\,b}$
~and~ Juan Luis Vazquez$^{\,c}$}

\maketitle

\begin{abstract}
The purpose of this paper is to prove local  upper and lower bounds for  weak solutions of semilinear elliptic equations of the form $-\Delta u= c u^p$, with $0<p<p_s=(d+2)/(d-2)$, defined on bounded domains of $\RR^d$, $d\ge 3$, without reference to the boundary behaviour. We give an explicit expression for all the involved constants. As a consequence, we obtain local Harnack inequalities with explicit constant, as well as gradient bounds.
\end{abstract}

\vskip 2cm

\noindent {\bf Keywords.} Local bounds, semilinear elliptic equations, regularity, Harnack inequality.\\[3mm]
\noindent {\bf Mathematics Subject Classification}. {\sc 35B45, 35B65,
35K55, 35K65}.

\hfill

\noindent (a) Departamento de Matem\'{a}ticas, Universidad
Aut\'{o}noma de Madrid,  28049 Madrid, Spain.
E-mail address:{\tt~matteo.bonforte@uam.es}. \\
Web-page:{\tt~http://www.uam.es/matteo.bonforte}

\noindent (b) Dipartimento di Matematica, Politecnico di Milano, Piazza Leonardo da Vinci 32, 20183 Milano, Italy.
E-mail address:{\tt ~gabriele.grillo@polimi.it}

\noindent (c) Departamento de Matem\'{a}ticas, Universidad
Aut\'{o}noma de Madrid,  28049 Madrid, Spain.
E-mail address:{\tt~juanluis.vazquez@uam.es}.\\ Web-page:{\tt~http://www.uam.es/juanluis.vazquez}


\newpage
\small
\tableofcontents
\normalsize





\newpage

\section{Introduction}

In this paper we obtain local  upper and lower estimates for the weak solutions of semilinear elliptic equations of the form
\begin{equation}\label{eq.f}
-\Delta u=f(u)
\end{equation}
posed in a bounded domain  $\Omega\subset \RR^d$. The choice of right-hand side we have in mind is $f(u)=\lambda\,u^p$ with $\lambda, p>0$. The range of exponents of interest will be $1<p<p_s:=(d+2)/(d-2)$ if $d\ge 3$, or $p>1$ if $d=1,2$. This problem is one of the most popular problems in nonlinear elliptic theory and  enjoys a large bibliography \cite{AP, BT, BN, dFLN, Evans, GNN, GS, GT, Giusti, LU, NN, MR0109940, Serrin, Serrin1, Tr1, Tr2} for $0\le p < p_s$ and \cite{BrNi,CGS} for $p=p_s$. We refrain from  attempting to give a complete bibliography for this nowadays classical problem.

We focus our attention on obtaining local estimates for solutions that are defined inside the domain without reference to their boundary behaviour. This is the notion of solution we use.

\begin{defn}\label{WeakSolloc} A local weak solution to equation $-\Delta u=f(u)$ in $\Omega$ is defined as a function  $u\in W^{1,2}_{\rm loc}(\Omega)$ with $f(u)\in L^1_{loc}(\Omega) $ which satisfies
\begin{equation}\label{local.weak.ell.p.1}
\int_K\left[\nabla u \cdot\nabla\varphi-f(u)\varphi\right]\dx=0
\end{equation}
for any subdomain with compact closure $K\subset \Omega$ and all bounded $\varphi \in C^1_0(K)$.
\end{defn}

Our aim is to contribute quantitative estimates in the form of upper bounds for solutions of any sign, lower bounds for positive solutions, and also local Harnack inequalities  and gradient bounds. By quantitative estimates we mean keeping track of all the constants during the proofs. As far as we know, there does not exist in literature a systematic set of quantitative estimates of local upper and lower bounds, and neither of the Harnack constant, in the form we explicitly provide here.  We recall that the quantitative control of the constants of such inequalities may have an important role in the applications; it is needed for instance in the results of \cite{BGV-Domains} on the asymptotic properties of solutions of the fast diffusion equation in bounded domains.

\medskip

\noindent {\bf Contents and main results}. We start with a section devoted to basic energy estimates.  We then consider in Section \ref{lub} the upper estimates for nonnegative solutions of the equation $-\Delta u=\lambda u^p$. The exponent range is $0\le p<p_s$, a main restriction of the theory, as it is already well known. See also \cite{Cabre} for $\LL^\infty$-bounds of different type for Equation \eqref{eq.f} with more general nonlinearities.  Our first main result, Theorem \ref{thm.local.upper}, can be considered as a smoothing effect with very precise constants; it is much simpler for $p\le 1$, but we  also obtain the more complicated and novel estimates for $1<p<p_s$. Next we obtain local upper estimates for $-\Delta u=b(x)u$ with unbounded coefficients in Theorem \ref{Thm.local.upper.b.DG} and we apply them to the case $b(x)=u^{p-1}$ in Theorem \ref{Thm.local.upper.b.p}.

In Section \ref{lb} we prove quantitative lower estimates, Theorems \ref{thm.local.lower}, \ref{thm.local.lower.pc}. We prove Harnack inequalities in Theorems \ref{Thm.Harnack.ps}, \ref{Harnack.p<1} and \ref{Harnack.p<pc}. All of these results appear to be well known from a qualitative point of view.
Let us mention that, as far as we know, the Harnack inequality for solutions to \eqref{local.weak.ell.p.1} when $p>1$ is not stated explicitly in the literature. The fact that the ``constant" involved has to depend on $u$ when $p_c\le p<p_s$ is confirmed by the results of \cite{BV-ADV}, \cite{DGV} applied to separation of variable solutions of parabolic problems, see also the very recent monograph \cite{DGVbook}. This is also related to the fact that, in the range $p_c\le p<p_s$, there exist (very weak) singular solutions. Notice also that in such a range the notion of weak and very weak solution is really different, cf. \cite{dPMP, MP1, P1,P2,P3,P4}.

In Section \ref{sect.absolute} we derive quantitative absolute upper (for $1<p<p_c$) and lower bounds (for $0\le p<1$) which are new as far as we know, cf. Theorem \ref{Thm.Absolute.0.pc}. The last section is devoted to quantitative gradient estimates, cf. Theorem \ref{thm.local.upper.grad}, and absolute upper bounds for the gradient when $1<p<p_c$, cf. Theorem \ref{thm.local.abs.grad}.

As a consequence of the above theory, we conclude that functions in the so-called De Giorgi class (satisfying Sobolev and local reverse Sobolev inequalities, at least at the level of truncates) are indeed locally bounded functions.

Much of the known theory takes into account  boundary conditions of different types: Dirichlet, Neumann, Robin, or other. Our results apply to all those cases. We will study the precise estimates for the Dirichlet problem in an upcoming paper \cite{BGV-Dirichlet}.

\section{Preliminaries. Local energy estimate}\label{sec.lee}

We shall pursue in the sequel the well-known idea that local weak solutions satisfy reverse Sobolev or Poincar\'e inequalities. Such local reverse inequalities are the key to prove local upper and lower estimates of next sections, and indeed imply that such functions are H\"older continuous. We comment that a similar line of reasoning could be adapted to deal with function belonging to suitably defined De Giorgi classes.
\begin{lem}[Energy Estimates]\label{Lem.energy}
Let $\Omega\subset\RR^d$ be a bounded domain, and let $p\,\ge0$ and $\lambda>0$. Let $u$ be a local nonnegative  weak solution in $\Omega$ to $-\Delta u = \lambda u^p$. Then the following energy equality holds true for any $\delta>0$\,, $\alpha\neq -1$ and any  positive  test function $\varphi\in C^2(\Omega)$ and compactly supported in $\Omega$:
\begin{equation}\label{local.energy.identity}
4\alpha\,\int_\Omega \big|\nabla \left((u+\delta)^{\frac{\alpha+1}{2}}\right)\big|^2\varphi\dx=
\lambda(\alpha+1)^2\int_\Omega u^p(u+\delta)^{\alpha}\varphi\dx
    +(\alpha+1)\int_\Omega (u+\delta)^{\alpha+1}\Delta\varphi \dx.
\end{equation}
Moreover, for any $\delta\ge 0$ we have the Caccioppoli estimates
\begin{equation}\label{Cacciopoli}\begin{split}
\lambda\int_\Omega \frac{u^p}{u+\delta}\varphi\dx + \int_\Omega \big|\nabla \log(u+\delta)\big|^2 \varphi\,\dx
    \le \int_\Omega \frac{|\nabla\varphi|^2}{\varphi}\dx.
\end{split}
\end{equation}
Local subsolutions $\underline{u}$ of $-\Delta \underline{u} \le \lambda \underline{u}^p$ satisfy,  for $\alpha\not=-1$ and $\delta>0$:
\begin{equation}\label{local.energy.identity.subsols}
4\alpha\int_\Omega \big|\nabla \left((\underline{u}+\delta)^{\frac{\alpha+1}{2}}\right)\big|^2\varphi\dx
\le\lambda (\alpha+1)^2\int_\Omega \underline{u}^p(\underline{u}+\delta)^{\alpha}\varphi\dx
    +(\alpha+1)\int_\Omega (\underline{u}+\delta)^{\alpha+1}\Delta\varphi \dx\,
\end{equation}
while local supersolution $-\Delta \overline{u} \ge \lambda \overline{u}^p$ satisfy, for any $\alpha\not=-1$ and $\delta>0$:
\begin{equation}\label{local.energy.identity.supersols}
\frac{4\alpha}{(\alpha+1)^2}\int_\Omega \big|\nabla\left(( \overline{u}+\delta)^{\frac{\alpha+1}{2}}\right)\big|^2\varphi\dx
\ge\lambda \int_\Omega \overline{u}^p(\overline{u}+\delta)^{\alpha}\varphi\dx
    +\frac{1}{\alpha+1}\int_\Omega (\overline{u}+\delta)^{\alpha+1}\Delta\varphi \dx\,,
\end{equation}
and the Caccioppoli estimates also work.
\end{lem}

\noindent\textbf{Remark. }Notice that when $\alpha>-1$, we can let $\delta=0$ in the energy identity \eqref{local.energy.identity} to get
\begin{equation}\label{local.energy.identity.2}
4\alpha\,\int_\Omega \big|\nabla \left(u^{\frac{\alpha+1}{2}}\right)\big|^2\varphi\dx=
\lambda(\alpha+1)^2\int_\Omega u^{p+\alpha}\varphi\dx
    +(\alpha+1)\int_\Omega u^{\alpha+1}\Delta\varphi \dx.
\end{equation}
The same remark applies to subsolutions:
\begin{equation}\label{local.energy.identity.subsols.2}
4\alpha\int_\Omega \big|\nabla \left(\underline{u}^{\frac{\alpha+1}{2}}\right)\big|^2\varphi\dx
\le\lambda (\alpha+1)^2\int_\Omega \underline{u}^{p+\alpha}\varphi\dx
    +(\alpha+1)\int_\Omega \underline{u}^{\alpha+1}\Delta\varphi \dx\,
\end{equation}

\noindent {\sl Proof.~}Let $\varphi\in C^2(\Omega)\cap C^1_0(\overline{\Omega})$ and $\delta\ge 0$. Multiply $-\Delta u$ by $(u+\delta)^\alpha\varphi$, with $\alpha\neq -1$ and integrate by parts to get
\begin{equation}\label{int.parts.alpha.=}\begin{split}
-\int_\Omega \varphi (u+\delta)^\alpha\Delta u\dx
    &=\int_\Omega \nabla\varphi \cdot \big(\nabla u\big) (u+\delta)^\alpha\dx
        + \alpha\int_\Omega \varphi (u+\delta)^{\alpha-1}\big|\nabla u\big|^2\dx\\
    &=-\frac{1}{\alpha+1}\int_\Omega (u+\delta)^{\alpha+1}\Delta\varphi \dx
        + \frac{4\alpha}{(\alpha+1)^2}\int_\Omega \big|\nabla (u+\delta)^{\frac{\alpha+1}{2}}\big|^2\varphi\dx.\\
\end{split}
\end{equation}
For local weak solutions of $-\Delta u=\lambda u^p$, the above equality immediately gives the energy identity \eqref{local.energy.identity} for $\alpha\neq -1$. Similar considerations hold, in the stated range of $\alpha$, for sub and supersolutions. To derive the Cacciopoli estimate we use the test function $\varphi/(u+\delta)$ to get
\[\begin{split}
0 &\le \lambda\int_\Omega \frac{u^p}{u+\delta}\varphi\dx
        = -\int_\Omega \frac{\varphi}{u+\delta}\Delta u\dx
        = -\int_\Omega \frac{\varphi}{(u+\delta)^2}\big|\nabla u\big|^2\dx
        + \int_\Omega \frac{\nabla\varphi\cdot\nabla u}{u+\delta}\frac{\sqrt{\varphi}}{\sqrt{\varphi}}\dx\\
    &\le -\int_\Omega \varphi\big|\nabla \log(u+\delta)\big|^2\dx
        + \frac{1}{2}\int_\Omega \frac{|\nabla\varphi|^2}{\varphi}\dx
        +\frac{1}{2} \int_\Omega \big|\nabla \log(u+\delta)\big|^2\varphi\dx \\
    &\le-\frac{1}{2}\int_\Omega \varphi\big|\nabla \log(u+\delta)\big|^2\dx
        +\frac{1}{2}\int_\Omega \frac{|\nabla\varphi|^2}{\varphi}\dx\,,
\end{split}
\]
where we have used the inequality  $a\cdot b\le (|a|^2+|b|^2)/2$.\qed

\noindent We shall also need the following particular computation.

\begin{lem}\label{lem.test.funct} Fix two balls $B_{R_1}\subset B_{R_0}\subset\subset\Omega$. Then there exists a test function $\varphi\in C_0^1(B_{R_0})$,  with $\nabla \varphi\equiv0$ on $\partial\Omega$, which is radially symmetric and piecewise $C^2$ as a function of $r$, satisfies $\supp(\varphi)=B_{R_0}$ and $\varphi=1$ on $B_{R_1}$, and moreover satisfies the bounds
\begin{equation}\label{test.estimates}
\|\nabla\varphi\|_\infty\le \frac{4}{R_0-R_1}\qquad\mbox{and}\qquad \|\Delta\varphi\|_\infty\le \frac{4d}{(R_0-R_1)^2}.
\end{equation}
\end{lem}
\noindent {\sl Proof.~}Consider the radial test function defined on $B_{R_0}$
\begin{equation}\label{test.funct}
\varphi(|x|)=\left\{
\begin{array}{lll}
1\,&\qquad\mbox{if }0\le |x|\le R_1\\[3mm]
1-\frac{2(|x|-R_1)^2}{(R_0-R_1)^2}\,&\qquad\mbox{if }R_1<|x|\le \frac{R_0+R_1}{2}\\[3mm]
\frac{2(R_0-|x|)^2}{(R_0-R_1)^2}\,&\qquad\mbox{if }\frac{R_0+R_1}{2}< |x|\le R_0\\[3mm]
0\,&\qquad\mbox{if }|x|>R_0\\[3mm]
\end{array}
\right.
\end{equation}
for any $0<R_1<R_0$. We have
\begin{equation*}
\nabla \varphi(|x|)=\left\{
\begin{array}{lll}
0\,&\qquad\mbox{if }0\le |x|\le R_1 \mbox{ or if }|x|>R_0\\[3mm]
-\frac{4(|x|-R_1)}{(R_0-R_1)^2}\frac{x}{|x|}\,&\qquad\mbox{if }R_1<|x|\le \frac{R_0+R_1}{2}\\[3mm]
-\frac{4(R_0-|x|)}{(R_0-R_1)^2}\frac{x}{|x|}\,&\qquad\mbox{if }\frac{R_0+R_1}{2}<|x|\le R_0\\[3mm]
\end{array}
\right.
\end{equation*}
and, recalling that $\Delta\varphi(|x|)=\varphi''(|x|)+(d-1)\varphi'(|x|)/|x|$,
\begin{equation*}
\Delta \varphi(|x|)=\left\{
\begin{array}{lll}
0\,&\qquad\mbox{if }0\le |x|\le R_1 \mbox{ or if }|x|>R_0\\[3mm]
-\frac{4}{(R_0-R_1)^2}-\frac{d-1}{|x|}\frac{4(|x|-R_1)}{(R_0-R_1)^2}\,&\qquad\mbox{if }R_1<|x|\le \frac{R_0+R_1}{2}\\[3mm]
-\frac{4}{(R_0-R_1)^2}-\frac{d-1}{|x|}\frac{4(R_0-|x|)}{(R_0-R_1)^2}\,&\qquad\mbox{if }\frac{R_0+R_1}{2}<|x|\le R_0\\[3mm]
\end{array}
\right.
\end{equation*}
As a consequence we easily obtain the bounds \eqref{test.estimates}.\qed

\begin{cor}[Quantitative Caccioppoli Estimates]\label{Lem.energy2} Let $\delta\ge 0$. Let $\Omega\subset\RR^d$ be a bounded domain, and let $p\,\ge0$ and $\lambda>0$. Let $u$ be a local positive  weak solution in $\Omega$ to $-\Delta u = \lambda u^p$. For any $B_{R}\subset B_{R_0}\subset\subset\Omega$ we have
\begin{equation}\label{Cacciopoli.Quantitativa}\begin{split}
\lambda\int_{B_R} \frac{u^p}{u+\delta}\dx+\int_{B_R} \big|\nabla \log(u+\delta)\big|^2\dx
    \le \frac{8\omega_d R_0^d}{(R_0-R)^2}
\end{split}
\end{equation}
where $\omega_d$ denotes the volume of the unit ball in $\mathbb R^d$\,.
\end{cor}
\noindent {\sl Proof.~}We use \eqref{Cacciopoli}, using the test function $\varphi$ of Lemma \ref{lem.test.funct} with $R$ replacing $R_1$:
\begin{equation*}\begin{split}
\lambda\int_{B_R} \frac{u^p}{u+\delta}\dx + \int_{B_R} \big|\nabla \log(u+\delta)\big|^2\dx
    &\le \lambda\int_\Omega \frac{u^p}{u+\delta}\varphi\dx + \int_\Omega \varphi\big|\nabla \log(u+\delta)\big|^2\dx\\
    &\le \int_\Omega \frac{|\nabla\varphi|^2}{\varphi}\dx
    \le \frac{8|\supp(\varphi)|}{(R_0-R)^2}
    =\frac{8\omega_d R_0^d}{(R_0-R)^2}\mbox{.\qed}
\end{split}
\end{equation*}

\noindent Note that the case $\delta>0$ follows immediately from the case $\delta=0$ since $u\ge 0$.

\noindent\textbf{Remark. } Letting $\delta=0$ in the Caccioppoli estimates \eqref{Cacciopoli.Quantitativa} shows that
\begin{equation}\label{Abs.Upper.p-1}
\lambda\int_{B_R} u^{p-1}\dx\le\frac{8\omega_d R_0^d}{(R_0-R)^2}\,
\end{equation}
When $p>1$ this yields a \textit{local absolute upper bound} for the local $\LL^{p-1}$-norm, a fact that will allow to conclude an absolute local $\LL^\infty$-bound in the range $1<p<p_c:=d/(d-2)$, as we shall see in Section \ref{sect.absolute}. This absolute upper bound represents a novelty both because it is quantitative and because it is local: to our knowledge this is the first absolute local bound for elliptic equations. When $p=1$ such absolute bound is easily seen to be impossible, while in the case $0<p<1$ we get an absolute lower bound for the local $\LL^{p-1}$-integral, which is new, at least as far as we know. It will be used below.

\subsection{More general nonlinearities}
As long as we deal with local estimates, we can apply the method to a larger class of operators and nonlinearities. (i)  First of all, namely we can treat local solutions of:
\begin{equation}\label{eq.gen}
-\nabla\cdot A(x, u, \nabla u)=\lambda\,u^p\,,
\end{equation}
where $A$ is a Carath\'eodory function such that
\[
\nu_1 |\xi|^2\le A(x,u,\xi)\cdot\xi\le \nu_2 |\xi|^2\qquad\mbox{and}\qquad |A(x,u,\xi)|^2\le \nu_2|\xi|^2
\]
for suitable constants $0<\nu_1<\nu_2$. The proofs of the inequalities are the same, and the results contain $\nu_1$ (resp. $\nu_2$) depending on whether you consider subsolutions (resp.  supersolutions).

\noindent(ii) Second we can consider supersolutions of the problem
\begin{equation}\label{eq.gen1}
-\nabla\cdot A(x, u, \nabla u)=f(x,u)\,,
\end{equation}
as long as $f(u)\ge a_0\,u^p$ with $a_0>0$, since they are supersolutions of $-\nabla\cdot A(x, u, \nabla u)= a_0\,u^p$.

\noindent (iii) We can consider subsolutions of \eqref{eq.gen1} with
$f(u)\le a_1 (u + b_1)^p\,,$ and $a_1,b_1\ge 0$. Then we can obtain an estimate for $v=u+b_1$.

The only thing that changes a bit are the energy estimates, and it is not so difficult to keep track of the new constants throughout the proof. We have decided here to consider the model case, to simplify the presentation and to focus on the main ideas.

\noindent (iv) Other  semilinear problems of this type are treated in the literature. Thus, Ambrosetti and Prodi's book \cite{AP} discusses right-hand sides of the form $f(x,u)= \lambda u + c(u) + h(x)$, with $a\in \RR$\,, $c(\cdot)\in C^0(\RR)\cap\LL^\infty(\RR)$ and $h\in C^{0,\alpha}(\overline{\Omega})$, for some $\alpha\in (0,1)$\,. Such nonlinearities can be treated with the methods presented here as well. We refrain from dealing with it in this work.


\section{Local Upper Bounds}\label{lub}
This section is devoted to the proof of the upper bounds and we will provide two kinds of estimates.
We prove local upper bounds for nonnegative subsolutions, then by Kato's inequality it is easy to extend such results to solutions with any sign.

\subsection{Local upper bounds I. The upper Moser iteration}
The local upper bounds follow from the local Sobolev imbedding theorem on balls $B_R\subset\RR^d$
\begin{equation}\label{Sobolev.BR}
\|f\|_{\LL^{2^*}(B_{R})}^2\le\mathcal{S}_2^2 \left(\|\nabla f\|_{\LL^2(B_{R})}^2 +\frac{1}{R^2}\|f\|_{\LL^2(B_{R})}^2\right)
\end{equation}
where $\mathcal{S}_2=\mathcal{S}_2(B_1)$ is the best constant and $2^*=2d/(d-2)$.  We are requiring hereafter without any further comment that $d\ge3$. The  Sobolev inequality combines with the energy inequalities of Lemma \ref{Lem.energy} which can be considered as local reverse Sobolev (or Poincar\'e) inequalities. The proof of the local upper bounds goes though the celebrated Moser iteration. We adopt the notation $\|f\|_{\LL^{q}(B_{R})}=\|f\|_{q,R}$, we recall that  $|B_R|=\omega_d R^d$ and that $\fint_Xf(x)\dx=\int_X f(x)\dx/|X|$. Throughout this section we are considering nonnegative subsolutions $u$ to $-\Delta u=\lambda u^p$, unless otherwise explicitly stated.

\begin{thm}[Local Upper Estimates]\label{thm.local.upper}
Let $\Omega\subseteq\RR^d$ and let $\lambda>0$. {\rm (i)} Let $u\ge 0$ be a local weak subsolution to $-\Delta u = \lambda u^p$ in $\Omega$, with $1< p<p_s=2^*-1=(d+2)/(d-2)$. Then, for any
    $q > \overline{q}:=d(p-1)_+/2$  and for any $B_{R_\infty}\subset B_{R_0}\subseteq\Omega$, the following bound holds true
\begin{equation}\label{upper.p+1}
\|u\|_{\LL^\infty(B_{R_\infty})}
    \le I_{\infty,q}\;
   \left(\fint_{B_{R_0}} u^{q}\dx\right)^{\frac{1+(p-1)\mu}{q}}
   \left({\fint_{B_{R_\infty}}u^{p-1}\dx}\right)^{-\mu }
\end{equation}
where $\mu= d/(2q-d(p-1))=d/2(q-\overline{q})$\,, and the constant $I_{\infty,q}>0$ depends on $d,p,q,R_0,R_\infty$, but not on $\lambda$.

\noindent {\rm (ii)} For $0\le p\le1$ the estimate simplifies into
\begin{equation}\label{upper.p<1}
\|u\|_{\LL^\infty(B_{R_\infty})}
    \le I_{\infty,q}\;
  \left( \fint_{B_{R_0}} u^{q}\dx\right)^{1/q}.
\end{equation}
valid for all $q>0$. $I_{\infty,q}>0$ has the same dependence as before, and it also depends on $\lambda $ when $p=1$, but not otherwise.
\end{thm}

\begin{figure}[ht]
\centering
\includegraphics[height=6cm, width=12.2cm]{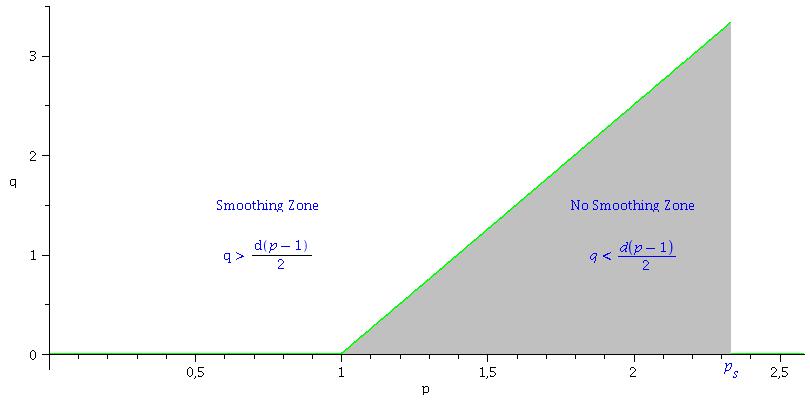}
 \noindent\textit{ Exponents of the local upper estimates.}
\end{figure}

\noindent {\bf Remarks on the result.} (i) Inequality \eqref{upper.p+1} is a kind of reverse H\"older inequality, indeed we can rewrite it as:
\begin{equation}\label{upper.p+1.2}
\|u\|_{\LL^{p-1}(B_{R_\infty})}^{\mu(p-1)}\|u\|_{\LL^\infty(B_{R_\infty})}
    \le C\;\|u\|_{\LL^{q}(B_{R_0})}^{1+\mu(p-1)}\,.
\end{equation}
Written in this form, it is clear from H\"older's inequality that a constant which makes \eqref{upper.p+1.2} true for a $q>\overline{q}$, make the same inequality true also for all $q'>q$\,. The same applies to \eqref{upper.p<1}\,.

\noindent(ii) The linear case $p=1$ is well known, cf. \cite{Evans,GT,Giusti}.

\noindent {\bf Remarks on the constant.} (i)  The proof below allows to find the following expression for the constant:
\begin{equation}\label{Const.Upper.q}\begin{split}
I_{\infty,q}=&\left[\frac{c_1\mathcal{S}_2^2
\omega_d^{\frac{2(p-1)_+}{d(p-1)}}}{(1-\rho)^2}\right]^{\frac{d}{2q-d(p-1)_+ }}\left\{\left(\frac{d}{d-2}\right)^d
\frac{2(d-2)}{\big(\sqrt{d}-\sqrt{d-2}\big)^2}\right.\times\\
&\times\left.\left[\Lambda_p+\frac{d-2}{q}+(1-\rho)^2
        \max\left\{\frac{d-2}{(dq)^2}|dq-(d-2)|,\,\frac14\right\}\right]\right\}^{\frac{d}{2q-d(p-1)_+ }}
        \\
    \end{split}
\end{equation}
where $\rho=R_\infty/R_0<1$ and we have used the convention $x_+/x=0$ when $x=0$ and, moreover, we have set $\Lambda_p=2$ if $p\ne 1$, $\Lambda_p= \lambda/4$ if $p=1$, with
\begin{equation}\label{c1.q}
c_1:=
\left\{\begin{array}{lll}
\frac{(d-2)q}{(d-2)q-d}\,&\mbox{if }q>\frac{d}{d-2}\,\\
\max\limits_{i=0,1}\frac{\left(\frac{d}{d-2}\right)^{k_0-1+i}\left[q-\frac{d(p-1)_+}{2}\right]+(p-1)_+\frac{d-2}{2}}
{\left|\left(\frac{d}{d-2}\right)^{k_0-1+i}\left[q-\frac{d(p-1)_+}{2}\right]+(p-1)_+\frac{d-2}{2}-1\right|}\,
&\mbox{if }0<q<\frac{d}{d-2}.
\end{array}
\right.
\end{equation}
(iii) When $q$ also satisfies $0<q<d/(d-2)$, we will require in the proof the additional condition
\begin{equation}\label{intero}
\frac{\log\frac{2^*-d(p-1)_+}{2q-d(p-1)_+}}{\log\frac{d}{d-2}}\qquad
\textrm{ is not an integer, and we let }\qquad k_0=i.p.\left[\frac{\log\frac{2^*-d(p-1)_+}{2q-d(p-1)_+}}{\log\frac{d}{d-2}}\right]\,,
\end{equation}
($i. p.$ is the integer part of a real number). Notice that taking $q=p+1>d(p-1)/2$ is possible if and only if $p<p_s=(d+2)/(d-2)$.

\noindent (iv) Of course, condition \eqref{intero} is not essential, in view of the remark after formula \eqref{upper.p+1.2}. In fact, let $q>\frac{d(p-1)_+}{2}$ be such that that $A(q):=\log\frac{2^*-d(p-1)_+}{2q-d(p-1)_+}/{\log\frac{d}{d-2}}$ is an integer. Take $\hat{q}\in(d(p-1)_+/2,q)$ such that $A(\hat{q})$ is not an integer. Then \eqref{upper.p+1} is valid with $\hat{q}$ instead of $q$.

\medskip

\noindent {\sl Proof.~}We are going to use the energy identity \eqref{local.energy.identity} for any $\alpha>-1$, $\alpha\not=0$, in the form \eqref{local.energy.identity.subsols} valid for subsolution, to prove $\LL^q-\LL^\infty$ local estimates via Moser iteration, keeping track of all the constants. We divide the proof in several steps.

\noindent$\bullet~$\textsc{Step 1. }Let $u$ as in Lemma \ref{Lem.energy} and $\varphi$ the test function of Lemma \ref{lem.test.funct}, which is supported in $B_{R_0}$ and such that $\varphi\equiv 1$ on $B_{R_1}$. The local Sobolev inequality \eqref{Sobolev.BR} on the ball $B_{R_1}$ applied to $f=u^{(\alpha+1)/2}$, together with the energy inequality \eqref{local.energy.identity.subsols} (we can take $\delta=0$ as in \eqref{local.energy.identity.subsols.2}), gives
\begin{equation}\label{upper.iter.0}\begin{split}
\left[\int_{B_{R_1}} u^{\frac{2^*}{2}(\alpha+1)}\dx\right]^{\frac{2}{2^*}}
     &\le \mathcal{S}_2^2 \left(\int_{B_{R_1}} \big|\nabla u^{\frac{\alpha+1}{2}}\big|^2\dx
         +\frac{1}{R_1^2}\int_{B_{R_1}} u^{\alpha+1}\dx\right)\\
   &\le \mathcal{S}_2^2 \left(\int_{B_{R_0}} \big|\nabla u^{\frac{\alpha+1}{2}}\big|^2\varphi\dx
       +\frac{1}{R_1^2}\int_{B_{R_1}} u^{\alpha+1}\dx\right)\\
    &=\mathcal{S}_2^2 \left(\frac{\lambda(\alpha+1)^2}{4|\alpha|}\int_{B_{R_0}} u^{p+\alpha}\varphi\dx
        +\frac{\alpha+1}{4|\alpha|}\int_{B_{R_0}} u^{\alpha+1}\Delta\varphi \dx
        +\frac{1}{R_1^2}\int_{B_{R_1}} u^{\alpha+1}\dx\right)\\
    &\le\mathcal{S}_2^2 \left(\frac{\lambda(\alpha+1)^2}{4|\alpha|}\int_{B_{R_0}} u^{p+\alpha}\dx
        +\left[\frac{(\alpha+1)\|\Delta \varphi\|_\infty}{4|\alpha|}+\frac{1}{R_1^2}\right]\int_{B_{R_0}} u^{\alpha+1} \dx\right)\\
    &\le\mathcal{S}_2^2 \left(\frac{\lambda(\alpha+1)^2}{4|\alpha|}\int_{B_{R_0}} u^{p+\alpha}\dx
        +\left[\frac{d(\alpha+1)}{|\alpha|(R_0-R_1)^2}+\frac{1}{R_1^2}\right]\int_{B_{R_0}} u^{\alpha+1} \dx\right)\\
\end{split}
\end{equation}
in the last step we have used the inequality $\|\Delta \varphi\|_\infty\le 4d/(R_0-R_1)^2$ of Lemma \ref{lem.test.funct}.

\noindent$\bullet~$\textsc{Step 2. }\textit{Caccioppoli estimates and the first iteration step. }Now we need to split two cases, namely $0\le p \le 1$ and $1<p<p_s$, and in both cases we will use the Caccioppoli estimate \eqref{Cacciopoli.Quantitativa} with $\delta=0$ which holds for any $p>0$ and reads
\begin{equation}\label{Caccioppoli.proof.2}
\lambda  \frac{ \|u\|_{p-1,R_{\infty}}^{p-1}}{|B_{R_{0}}|}\le \frac{8}{(R_0-R_\infty)^2}.
\end{equation}

\noindent \textit{Superlinear case: }$1<p<p_s$. We continue estimate \eqref{upper.iter.0} as follows:
\begin{equation}\label{upper.iter.1}\begin{split}
&\left[\int_{B_{R_1}} u^{\frac{2^*}{2}(\alpha+1)}\dx\right]^{\frac{2}{2^*}}\\ &\le\mathcal{S}_2^2\left(\frac{\lambda(\alpha+1)^2}{4|\alpha|}
        +\left[\frac{d(\alpha+1)}{|\alpha|(R_0-R_1)^2}+\frac{1}{R_1^2}\right]
            \frac{\int_{B_{R_0}} u^{\alpha+1} \dx}{\int_{B_{R_0}} u^{p+\alpha}\dx}\right)
            \int_{B_{R_0}} u^{p+\alpha}\dx\\
    &\le_{(a)}\mathcal{S}_2^2\left(\frac{\lambda(\alpha+1)^2}{4|\alpha|}
        +\left[\frac{d(\alpha+1)}{|\alpha|(R_0-R_1)^2}+\frac{1}{R_1^2}\right]\frac{|B_{R_0}|}{\int_{B_{R_0}} u^{p-1}\dx}\right)\int_{B_{R_0}} u^{p+\alpha}\dx\\
    &=\frac{\mathcal{S}_2^2|B_{R_0}|}{\|u\|_{p-1, R_0}^{p-1}}\left(\frac{\lambda(\alpha+1)^2}{4|\alpha|}
          \frac{\|u\|_{p-1, R_0}^{p-1}}{|B_{R_0}|}
        +\left[\frac{d(\alpha+1)}{|\alpha|(R_0-R_1)^2}+\frac{1}{R_1^2}\right]\right)
        \int_{B_{R_0}} u^{p+\alpha}\dx\\
   &\le_{(b)}\frac{\mathcal{S}_2^2|B_{R_0}|}{\|u\|_{p-1, R_0}^{p-1}}\left(\frac{2(\alpha+1)^2}{|\alpha|(R_0-R_1)^2}
        +\left[\frac{d(\alpha+1)}{|\alpha|(R_0-R_1)^2}+\frac{1}{R_1^2}\right]\right)
        \int_{B_{R_0}} u^{p+\alpha}\dx\\
   &=\frac{\mathcal{S}_2^2|B_{R_0}|}{(R_0-R_1)^2\|u\|_{p-1, R_0}^{p-1}}
        \left[\frac{1}{|\alpha|}\left(2(\alpha+1)^2+d(\alpha+1)\right)+\frac{(R_0-R_1)^2}{R_1^2}\right]
            \int_{B_{R_0}} u^{p+\alpha}\dx\,\\
\end{split}
\end{equation}
where in $(a)$ we have used the convexity in the variable $r>0$ of the function $N(r)=\log\|u\|_r^r$, the incremental quotient is increasing, hence choosing $\alpha+1 \ge\overline{\alpha}> 0$, we obtain
\[
\frac{N(p-1+\overline{\alpha})-N(\overline{\alpha})}{p-1}\le \frac{N(\alpha+p)-N(\alpha+1)}{p-1}\qquad\mbox{that is}\qquad
    \frac{\|u\|_{p-1+\overline{\alpha}}^{p-1+\overline{\alpha}}}{\|u\|_{\overline{\alpha}}^{\overline{\alpha}}} \le \frac{\|u\|_{\alpha+p}^{\alpha+p}}{\|u\|_{\alpha+1}^{\alpha+1}}
\]
Then we have
\[
\frac{\|u\|_{\alpha+p}^{\alpha+p}}{\|u\|_{\alpha+1}^{\alpha+1}}
\ge  \frac{\|u\|_{p-1+\overline{\alpha}}^{p-1+\overline{\alpha}}}{\|u\|_{\overline{\alpha}}^{\overline{\alpha}}}
=\frac{\|u\|_{p-1+\overline{\alpha}}^{\overline{\alpha}}}{\|u\|_{\overline{\alpha}}^{\overline{\alpha}}}\|u\|_{p-1+\overline{\alpha}}^{p-1}
\ge|B_{R_0}|^{\frac{-(p-1)}{\overline{\alpha}+p-1}}\, |B_{R_0}|^{\frac{p-1}{\overline{\alpha}+p-1}-1}\|u\|_{p-1}^{p-1}=\frac{\|u\|_{p-1}^{p-1}}{|B_{R_0}|}
\]
since by H\"older inequality:
\[
\frac{\|u\|_{p-1+\overline{\alpha}}}{\|u\|_{\overline{\alpha}}}\ge |B_R|^{\frac{-(p-1)}{\overline{\alpha}+p-1}}
\qquad\mbox{and}\qquad
\|u\|_{p-1+\overline{\alpha}}\ge |B_R|^{\frac{1}{\overline{\alpha}+p-1}-\frac{1}{p-1}}\|u\|_{p-1}.
\]
In $(b)$ we have used the Caccioppoli estimate \eqref{Caccioppoli.proof.2}.

\noindent \textit{Sublinear case: }$0\le p \le 1$. We first assume $0\le p<1$, we discuss the case $p=1$ separately. We continue estimate \eqref{upper.iter.0} as follows:
\begin{equation}\label{upper.iter.2}\begin{split}
\left[\int_{B_{R_1}} u^{\frac{2^*}{2}(\alpha+1)}\dx\right]^{\frac{2}{2^*}}    &\le\mathcal{S}_2^2\left(\frac{\lambda(\alpha+1)^2}{4|\alpha|} \frac{\int_{B_{R_0}} u^{p+\alpha}\dx}{\int_{B_{R_0}} u^{\alpha+1} \dx}
        +\frac{d(\alpha+1)}{|\alpha|(R_0-R_1)^2}+\frac{1}{R_1^2}\right)
            \int_{B_{R_0}} u^{\alpha+1}\dx\\
    &\le\mathcal{S}_2^2\left(\frac{2(\alpha+1)^2}{|\alpha|(R_0-R_\infty)^2}
        +\frac{d(\alpha+1)}{|\alpha|(R_0-R_1)^2}+\frac{1}{R_1^2}\right)\int_{B_{R_0}} u^{\alpha+1}\dx\\
    &=\frac{\mathcal{S}_2^2}{(R_0-R_\infty)^2}\left[\frac{1}{|\alpha|}\left(2(\alpha+1)^2+d(\alpha+1)\right)
        +\frac{(R_0-R_1)^2}{R_1^2}\right]\int_{B_{R_0}} u^{\alpha+1}\dx
\end{split}
\end{equation}
which follows by the convexity in the variable $r>0$ of the function $N(r)=\log\|u\|_r^r$, which implies that the incremental quotient is increasing, hence choosing $\alpha+1 \ge \overline{\alpha}:= \beta_0>0$, we obtain
\[
\frac{N(p-1+\overline{\alpha})-N(\overline{\alpha})}{p-1}\le \frac{N(\alpha+p)-N(\alpha+1)}{p-1}\qquad\mbox{that is}\qquad
    \frac{\|u\|_{p-1+\overline{\alpha}}^{p-1+\overline{\alpha}}}{\|u\|_{\overline{\alpha}}^{\overline{\alpha}}} \le \frac{\|u\|_{\alpha+p}^{\alpha+p}}{\|u\|_{\alpha+1}^{\alpha+1}}
\]
hence
\[\begin{split}
\frac{\int_{B_{R_0}} u^{p+\alpha}\dx}{\int_{B_{R_0}} u^{\alpha+1} \dx}
    &=\frac{\|u\|_{\alpha+p}^{\alpha+p}}{\|u\|_{\alpha+1}^{\alpha+1}}
    \le\frac{\|u\|_{\overline{\alpha}-(1-p)}^{\overline{\alpha}-(1-p)}}{\|u\|_{\overline{\alpha}}^{\overline{\alpha}}}
    \le \frac{|B_{R_0}|^{\frac{1-p}{\overline{\alpha}}}}{\|u\|_{\overline{\alpha}}^{1-p}}
    \le\frac{\|u\|^{p-1}_{p-1,R_0}}{|B_{R_0}|}\le\frac{8}{\lambda(R_0-R_\infty)^2}
\end{split}\]
again by H\"older inequalities, we just stress on the last step in which we have used that
\[
\frac{\|u\|_{p-1,R_0}}{|B_{R_0}|^{\frac{1}{p-1}}}\le \frac{\|u\|_{\overline{\alpha}}}{|B_{R_0}|^{\frac{1}{\overline{\alpha}}}}\qquad\mbox{hence}\qquad \frac{|B_{R_0}|^{\frac{1-p}{\overline{\alpha}}}}{\|u\|_{\overline{\alpha}}^{1-p}}\le
\frac{\|u\|^{p-1}_{p-1,R_0}}{|B_{R_0}|}\le\frac{8}{\lambda(R_0-R_\infty)^2}
\]
which is true since $p-1< 0<\overline{\alpha}$, and in the last step we have used the Caccioppoli estimate \eqref{Caccioppoli.proof.2}.

\noindent  Notice that when $p=1$, we obtain directly that
\begin{equation}\label{upper.iter.2.p=1}\begin{split}
&\left[\int_{B_{R_1}} u^{\frac{2^*}{2}(\alpha+1)}\dx\right]^{\frac{2}{2^*}}
\le\mathcal{S}_2^2\left(\frac{\lambda(\alpha+1)^2}{4|\alpha|}
        +\frac{d(\alpha+1)}{|\alpha|(R_0-R_1)^2}+\frac{1}{R_1^2}\right)
            \int_{B_{R_0}} u^{\alpha+1}\dx\\
    &=\frac{\mathcal{S}_2^2}{(R_0-R_\infty)^2}\left[\frac{1}{|\alpha|}\left(\frac{\lambda}{4}(\alpha+1)^2+d(\alpha+1)\right)
        +\frac{(R_0-R_1)^2}{R_1^2}\right]\int_{B_{R_0}} u^{\alpha+1}\dx
\end{split}
\end{equation}

\noindent\textit{The first iteration step. }We can write the first iteration step for all $p\ge 0$ in the following way: let $\beta=\alpha+1\ge \beta_0>0$ and recall that we are requiring $\beta\not=1$ as well, then inequalities \eqref{upper.iter.1} and \eqref{upper.iter.2} can be written as
\begin{equation}\label{iter.beta}\begin{split}
\left[\int_{B_{R_1}} u^{\frac{2^*}{2}\beta}\dx\right]^{\frac{2}{2^*}}
     &\le I(p,\beta,R_1,R_0)\int_{B_{R_0}} u^{\beta+(p-1)_+}\dx
\end{split}
\end{equation}
where
\begin{equation}\label{iter.beta.const}
I(p,\beta,R_1,R_0)=\frac{\mathcal{S}_2^2}{(R_0-R_1)^2}\frac{|B_{R_0}|}{\int_{B_{R_0}}u^{(p-1)_+}\dx}
        \left[\frac{\Lambda_p  \beta^2+d\beta}{|\beta-1|}+\frac{(R_0-R_1)^2}{R_1^2}\right]\,
\end{equation}
where $\Lambda_p=2$ if $p\ne 1$ and $\Lambda_p= \lambda/4$ if $p=1$.

\noindent$\bullet~$\textsc{Step 3. }\textit{The Moser iteration. }
Let us define the sequence of exponents $\beta_n>0$
so that
\[
\beta_n+(p-1)_+=\frac{2^*}{2}\beta_{n-1}\qquad\mbox{that is}\qquad \beta_n=\frac{2^*}{2}\beta_{n-1}-(p-1)_+
\]
it turns out that,  for any given $\beta_0$ and all $n\ge1$:
\begin{equation}\label{beta.n.def}\begin{split}
\beta_n
    &=\left[\frac{2^*}{2}\right]^n\,\left[\beta_0-(p-1)_+\sum_{k=0}^{n-1}\left(\frac{2^*}{2}\right)^{k-n}\right]
    =\left[\frac{2^*}{2}\right]^n\,\left[\beta_0-(p-1)_+\sum_{j=1}^{n}\left(\frac{2}{2^*}\right)^j\right]\\
    &=\left[\frac{2^*}{2}\right]^n\,\left[\beta_0-(p-1)_+\frac{d-2}{2}\left(1-\left(\frac{2}{2^*}\right)^n\right)\right]
    =\left[\frac{2^*}{2}\right]^n\,\left[\beta_0 - (p-1)_+\frac{d-2}{2}\right]+(p-1)_+\frac{d-2}{2}
\end{split}
\end{equation}
since $\sum_{j=1}^k s^j=(1-s^k)s/(1-s)$. Moreover we have that for all $p\ge 1$,
\[
\left(\frac{2^*}{2}\right)^{-n}\beta_n\xrightarrow[n\to\, \infty] \;\beta_0-\frac{d-2}{2}(p-1)_+.
\]
Requiring that $\beta_0>(p-1)_+(d-2)/2$, which will be assumed from now on, then implies that $\beta_n\to+\infty$  as $n\to+\infty$. We shall also require that $\beta_n\not=1$ for all $n$.

\noindent We will explicitly choose a decreasing sequence of radii $0<R_\infty<\ldots<R_n<R_{n-1}<\ldots<R_0$ in the next step,
in order to estimate explicitly the constants. The first iteration step then reads:
\begin{equation}\label{iter.beta.1}\begin{split}
\|u\|_{\frac{2^*}{2}\beta_n,R_n}&
    =\left[\int_{B_{R_n}} u^{\frac{2^*}{2}\beta_n}\dx\right]^{\frac{2}{2^*\beta_n}}
    \le I(p,\beta_n,R_n,R_{n-1})^{\frac{1}{\beta_n}} \left[\int_{B_{R_{n-1}}} u^{(p-1)_++\beta_n}\dx\right]^{\frac{1}{\beta_n}}\\
        &:=I_n^{\frac{1}{\beta_n}}\,\|u\|_{\beta_n+(p-1)_+,R_{n-1}}^{\frac{\beta_n+(p-1)_+}{\beta_n}}
         =I_n^{\frac{1}{\beta_n}}\,\|u\|_{\frac{2^*}{2}\beta_{n-1},R_{n-1}}^{\frac{2^*}{2}\frac{\beta_{n-1}}{\beta_n}}\,\\
\end{split}
\end{equation}
where the constants $I(p,\beta,R_1,R_0)$ are defined in \eqref{iter.beta.const}. Hence
\begin{equation}\label{Ik}
I_n=I(p,\beta_n,R_n,R_{n-1})=\frac{\mathcal{S}_2^2}{(R_{n-1}-R_n)^2}\frac{|B_{R_{n-1}}|}{\int_{B_{R_{n-1}}}u^{(p-1)_+}\dx}
        \left[\frac{2\beta_n^2+d\beta_n}{|\beta_n-1|}+\frac{(R_{n-1}-R_n)^2}{R_n^2}\right]
\end{equation}
Iterating the above inequality yields
\begin{equation}
\begin{split}
\|u_n\|_{\frac{2^*}{2}\beta_n,R_n}
    &\le I_n^{\frac{1}{\beta_n}}\,
        \|u_n\|_{\frac{2^*}{2}\beta_{n-1}\,\,R_{n-1}}^{\frac{2^*}{2}\frac{\beta_{n-1}}{\beta_n}}
     \le I_n^{\frac{1}{\beta_n}}\,I_{n-1}^{\frac{2^*}{2}\frac{1}{\beta_n}}
        \|u_n\|_{\frac{2^*}{2}\beta_{n-2}\,\,R_{n-2}}^{\left(\frac{2^*}{2}\right)^2\frac{\beta_{n-2}}{\beta_n}}\\
    &\le I_n^{\frac{1}{\beta_n}}\,I_{n-1}^{\frac{2^*}{2}\frac{1}{\beta_n}}
        \,\ldots\,I_1^{\left(\frac{2^*}{2}\right)^{n-1}\frac{1}{\beta_n}}
        \|u_n\|_{\frac{2^*}{2}\beta_0,R_{0}}^{\left(\frac{2^*}{2}\right)^n\frac{\beta_{0}}{\beta_n}}
    \le \prod_{j=1}^n I_j^{\left(\frac{2^*}{2}\right)^{n-j}\frac{1}{\beta_n}}
        \|u_n\|_{\frac{2^*}{2}\beta_0,R_{0}}^{\left(\frac{2^*}{2}\right)^n\frac{\beta_{0}}{\beta_n}}
\end{split}
\end{equation}
with
\[
\beta_0>\frac{d-2}{2}(p-1)_+\qquad\mbox{or}\qquad \overline{q}:=\frac{2^*}{2}\beta_0>\frac{d(p-1)_+}{2}.
\]
Taking the limit as $n\to\infty$ we obtain
\begin{equation}\label{local.upper.p}\begin{split}
\|u\|_{\infty,R_\infty}
    &=\lim_{n\to\infty}\|u\|_{\frac{2^*}{2}\beta_n,R_n}
    \le \lim_{n\to\infty}\prod_{k=1}^n I_k^{\left(\frac{2^*}{2}\right)^{n-k}\frac{1}{\beta_n}}
        \|u\|_{\frac{2^*}{2}\beta_0,R_{0}}^{\frac{\beta_{0}}{\beta_0-\frac{d-2}{2}(p-1)_+}}\\
    &\le \lim_{n\to\infty}\prod_{k=1}^n I_k^{\left(\frac{2^*}{2}\right)^{n-k}\frac{1}{\beta_n}}
        \|u\|_{\frac{2^*}{2}\beta_0,R_{0}}^{\frac{\beta_{0}}{\beta_0-\frac{d-2}{2}(p-1)_+}}
        =I_\infty \|u\|_{\overline{q},R_{0}}^{\frac{2\overline{q}}{2\overline{q}-d(p-1)_+}}
\end{split}
\end{equation}
notice that the penultimate passage follows because we shall see below that $\prod_{k=1}^n I_k^{\left(\frac{2^*}{2}\right)^{n-k}\frac{1}{\beta_n}}$ has a limit $I_\infty$ as $n\to+\infty$.

\noindent As a consequence of the above estimates $u\in\LL^\infty$, so that the above bounds holds for any $\overline{q}>d(p-1)_+/2$ as stated, provided we show that the constant $I_\infty$ is finite and can be estimated as in \eqref{Const.Upper.q}.

\noindent$\bullet~$\textsc{Step 4. }\textit{Estimating all the constants. }Now it remains to estimate $I_\infty$. We will prove later that
\begin{equation}\label{Ik.I0}
I_k\le I_0(p)\left[\frac{2^*}{2}\right]^{2k}
\end{equation}
where $I_0(p)$ will have the explicit form given in formula \eqref{I0}. Using such bound we show that
\[
\begin{split}
I_\infty
    &=\lim_{n\to\infty}\prod_{k=1}^n I_k^{\left(\frac{2^*}{2}\right)^{n-k}\frac{1}{\beta_n}}
        = \lim_{n\to\infty} \exp\left[\sum_{k=1}^{n} \log\left(I_k^{\left(\frac{2^*}{2}\right)^{-k}\left(\frac{2^*}{2}\right)^{n}\frac{1}{\beta_n}}\right)\right]\\
    &= \lim_{n\to\infty}\exp\left[\left(\frac{2^*}{2}\right)^{n}\frac{1}{\beta_n}
        \sum_{k=1}^n \left(\frac{2}{2^*}\right)^{k}\log\left(I_k\right)\right]\\
    &\le \lim_{n\to\infty}\exp\left[\left(\frac{2^*}{2}\right)^{n}\frac{1}{\beta_n}
        \sum_{k=1}^n \left(\frac{2}{2^*}\right)^{k}\log\left(I_0\left[\frac{2^*}{2}\right]^{2k}\right)\right]\\
    &=\lim_{n\to\infty}\exp\left[\left(\frac{2^*}{2}\right)^{n}\frac{1}{\beta_n}
        \left(\log(I_0)\sum_{k=1}^n \left(\frac{2}{2^*}\right)^{k}
        +2\log\left(\frac{2^*}{2}\right)\sum_{k=1}^n \left(\frac{2}{2^*}\right)^{k}k\,\right)\right]\\
\end{split}
\]
\[\begin{split}
    &= \exp\left[\frac{1}{\beta_0-\frac{d-2}{2}(p-1)_+ }
        \left(\log(I_0)\sum_{k=1}^{+\infty} \left(\frac{2}{2^*}\right)^{k}
        +2\log\left(\frac{2^*}{2}\right)\sum_{k=1}^{+\infty} \left(\frac{2}{2^*}\right)^{k}k\,\right)\right]\\
    &= \exp\left[\frac{2}{2\beta_0-(d-2)(p-1)_+ }
        \left(\log(I_0)\frac{d-2}{2}
        +2\log\left(\frac{2^*}{2}\right)\frac{d(d-2)}{4}\,\right)\right]\\
\end{split}
\]
\[
\begin{split}
    &= \exp\left[\frac{d-2}{2\beta_0-(d-2)(p-1)_+ }\log(I_0)
        +\frac{d(d-2)}{2\beta_0-(d-2)(p-1)_+ }
        \log\left(\frac{2^*}{2}\right)\,\right]\\
    &= I_0^{\frac{d-2}{2\beta_0-(d-2)(p-1)_+ }}\,\left(\frac{2^*}{2}\right)^{\frac{d(d-2)}{2\beta_0-(d-2)(p-1)_+ }}
    =\left[I_0\,\left(\frac{2^*}{2}\right)^d\right]^{\frac{d-2}{2\beta_0-(d-2)(p-1)_+ }}
\end{split}\]
We shall now obtain an explicit estimate for $I_0$ in order to finally obtain \eqref{Const.Upper.q}.

\noindent\textit{Estimating $I_k$. }We want to obtain estimates \eqref{Ik.I0}, and to this end we choose a decreasing sequence of radii $0<R_\infty<\ldots<R_k<R_{k-1}<\ldots<R_0$ such that
\[
(R_{k-1}-R_k)^2=\,(R_0-R_\infty)^2\,\frac{c_0^2}{\beta_k}
\qquad\mbox{with}\qquad
c_0=\left(\sum_{k=1}^\infty \sqrt{\frac{1}{\beta_k}}\right)^{-1} <+\infty\,
\]
so that
\[
\sum_{k=1}^\infty (R_{k-1}-R_k)=R_0-R_\infty
.
\]
We now estimate $I_k$:
\[\begin{split}
&I_k
    =\frac{\mathcal{S}_2^2}{(R_{k-1}-R_k)^2}\frac{|B_{R_{k-1}}|}{\int_{B_{R_{k-1}}}u^{(p-1)_+}\dx}
        \left[\frac{\Lambda_p\beta_k^2+d\beta_k}{|\beta_k-1|}+\frac{(R_{k-1}-R_k)^2}{R_k^2}\right]\\
    &=\frac{\mathcal{S}_2^2\beta_k^2}{|\beta_k-1|(R_{k-1}-R_k)^2}\frac{|B_{R_{k-1}}|}{\int_{B_{R_{k-1}}}u^{(p-1)_+}\dx}
        \left[\Lambda_p+\frac{d}{\beta_k}+\frac{(R_{k-1}-R_k)^2}{R_k^2}\frac{|\beta_k-1|}{\beta_k^2}\right]\\
    &\le_{(a)}\frac{\mathcal{S}_2^2\beta_k^3}{c_0^2|\beta_k-1|(R_0-R_\infty)^2}\frac{|B_{R_0}|}{\int_{B_{R_\infty}}u^{(p-1)_+}\dx}
        \left[\Lambda_p+\frac{d}{\beta_0}+\frac{(R_0-R_\infty)^2}{R_\infty^2}\max\left\{\frac{|\beta_0-1|}{\beta_0^2},\,\frac14\right\}\right]\\
    &\le_{(b)}\frac{c_1\mathcal{S}_2^2\beta_k^2}{c_0^2(R_0-R_\infty)^2}\frac{|B_{R_0}|}{\int_{B_{R_\infty}}u^{(p-1)_+}\dx}
        \left[\Lambda_p+\frac{d}{\beta_0}+\frac{(R_0-R_\infty)^2}{R_\infty^2}\max\left\{\frac{|\beta_0-1|}{\beta_0^2},\,\frac14\right\}\right]\\
    &\le_{(c)}\frac{2c_1\mathcal{S}_2^2\left[\beta_0 - (p-1)_+\frac{d-2}{2}\right]}{c_0^2(R_0-R_\infty)^2}\frac{|B_{R_0}|}{\int_{B_{R_\infty}}u^{(p-1)_+}\dx}
        \left[\Lambda_p+\frac{d}{\beta_0}+\frac{(R_0-R_\infty)^2}{R_\infty^2}\max\left\{\frac{|\beta_0-1|}{\beta_0^2},
        \,\frac14\right\}\right]\left[\frac{2^*}{2}\right]^{2n}\\
    &\le_{(d)}\frac{2(d-2)c_1\mathcal{S}_2^2|B_{R_0}|}{\big(\sqrt{d}-\sqrt{d-2}\big)^2(R_0-R_\infty)^2\int_{B_{R_\infty}}u^{(p-1)_+}\dx}
     \left[\Lambda_p+\frac{d}{\beta_0}+\frac{(R_0-R_\infty)^2}{R_\infty^2}\max\left\{\frac{|\beta_0-1|}{\beta_0^2},
        \,\frac14\right\}\right]\left[\frac{2^*}{2}\right]^{2n}\\
\end{split}\]
in $(a)$ we have used that
\begin{equation}\label{stima}
\frac{|\beta_k-1|}{\beta_k^2}\le\max\left\{\frac{|\beta_0-1|}{\beta_0^2},\,\frac14\right\}.
\end{equation}
In $(b)$ we have also used the inequality
\begin{equation}\label{beta.k.k-1}
\frac{\beta_k}{|\beta_k-1|}\le c_1:=
\left\{\begin{array}{lll}
\frac{\beta_0}{\beta_0-1}\, &\quad\mbox{if }\beta_0>1\,\\
\max\limits_{i=0,1}\frac{\beta_{k_0+i}}{|\beta_{k_0+i}-1|}
&\quad\mbox{if }0<\beta_0<1\,\mbox{ with }k_0=k_0=i.p.\left[\frac{\log\frac{1-(p-1)_+\frac{d-2}{2}}{\beta_0-(p-1)_+\frac{d-2}{2}}}{\log\frac{d}{d-2}}\right].
\end{array}
\right.
\end{equation}
which we state in the general case $p\not=1$ for later use and we now prove.
First notice that the numerical inequality
\[
\frac{s}{|s-1|}\le \max\left\{\frac{a}{1-a},\frac{b}{b-1}\right\}\qquad\mbox{for all }0< a<1<b<+\infty\, \mbox{ and all }s\in[0,a]\cup[b,\infty)\,
\]
holds true. When $\beta_0>1$  \eqref{beta.k.k-1} follows applying such numerical inequality to $s=\beta_k$ and noticing that
$\beta_k>\beta_0=b>1$ and that the function $x/|x-1|$ is decreasing when $x>1$. Suppose instead that $0<\beta_0<1$. Notice that, since we are
also requiring that $\beta_0>(p-1)_+(d-2)/2$, this is possible only when $0< p< p_c=d/(d-2)<p_s$. We define $k_0$ to
be the greatest integer for which $\beta_k<1$, so that $\beta_{k_0+1}>1$, so that
\[
\beta_{k_0}<1<\beta_{k_0+1}\qquad\mbox{with}\qquad k_0=i.p.\left[\frac{\log\frac{1-(p-1)_+\frac{d-2}{2}}{\beta_0-(p-1)_+\frac{d-2}{2}}}{\log\frac{d}{d-2}}\right]
\]
and we shall take $\beta_0\in(0,1)$ such that
\begin{equation}\label{non.integer}
\frac{\log\frac{1-(p-1)_+\frac{d-2}{2}}{\beta_0-(p-1)_+\frac{d-2}{2}}}{\log\frac{d}{d-2}}\  \rm is\ not\ an\ integer. \normalcolor
\end{equation}
The elementary properties of the function $x/|x-1|$ then show that, for all $k$:\normalcolor
\[\begin{split}
\frac{\beta_k}{|\beta_k-1|}\le \max\limits_{i=0,1}\frac{\beta_{k_0+i}}{|\beta_{k_0+i}-1|}
&=\max\limits_{i=0,1}\frac{\left(\frac{d}{d-2}\right)^{k_0+i}\left[\beta_0-(p-1)_+\frac{d-2}{2}\right]+(p-1)_+\frac{d-2}{2}}
{\left|\left(\frac{d}{d-2}\right)^{k_0+i}\left[\beta_0-(p-1)_+\frac{d-2}{2}\right]+(p-1)_+\frac{d-2}{2}-1\right|}\\
&=\max\limits_{i=0,1}\frac{\left(\frac{d}{d-2}\right)^{k_0-1+i}\left[\overline{q}-\frac{d(p-1)_+}{2}\right]+(p-1)_+\frac{d-2}{2}}
{\left|\left(\frac{d}{d-2}\right)^{k_0-1+i}\left[\overline{q}-\frac{d(p-1)_+}{2}\right]+(p-1)_+\frac{d-2}{2}-1\right|}
\end{split}\]
as claimed, where we have put $\beta_0=\frac{2}{2^*} \overline{q} =\frac{d-2}{d}\overline{q}$ and $\overline{q}$ has to be chosen such that \eqref{non.integer} holds.

\noindent In $(c)$ we have used that $\beta_k=\beta_0 (2^*/2)^k> \beta_0$
\begin{equation}\begin{split}
\beta_n
    =\left[\frac{2^*}{2}\right]^n\,\left[\beta_0 - (p-1)_+\frac{d-2}{2}\right]+(p-1)_+\frac{d-2}{2}\le 2\left[\frac{2^*}{2}\right]^n\,\left[\beta_0 - (p-1)_+\frac{d-2}{2}\right]
\end{split}
\end{equation}

\noindent Finally in $(d)$ we estimate $1/c_0^2$ as follows:
\[
\frac{1}{c_0^2}=\left(\sum_{k=1}^\infty \sqrt{\frac{1}{\beta_k}}\right)^2
\le \left(\sum_{k=1}^\infty \frac{1}{\left(\normalcolor\beta_0-(p-1)_+\frac{d-2}{2}\right)^{1/2}}
\left(\frac{2}{2^*}\right)^{\frac{k}{2}}\right)^2=\frac{1}{\big(\beta_0-(p-1)_+\frac{d-2}{2}\big)}\frac{d-2}{\big(\sqrt{d}-\sqrt{d-2}\big)^2}
\]
since  the explicit expression of $\beta_k$ shows that
\[
\beta_k\ge \left(\beta_0-(p-1)_+\frac{d-2}{2} \right)\left(\frac{2^*}{2}\right)^{k}\,
\]
and
\[
\sum_{k=1}^{+\infty}\left(\frac{2}{2^*}\right)^{k/2}=\sum_{k=1}^{+\infty}\left(\frac{d-2}{d}\right)^{k/2}
=\frac{\sqrt{d-2}}{\sqrt{d}-\sqrt{d-2}}.
\]
We conclude that we can take $I_0(p)$ as follows for any $p>0$:
\begin{equation}\label{I0}\begin{split}
I_0(p)&=\frac{2(d-2)}{\big(\sqrt{d}-\sqrt{d-2}\big)^2}\frac{c_1\mathcal{S}_2^2}{(R_0-R_\infty)^2}\frac{|B_{R_0}|}{\int_{B_{R_\infty}}u^{(p-1)_+}\dx}
        \left[\Lambda_p+\frac{d}{\beta_0}+\frac{(R_0-R_\infty)^2}{R_\infty^2}\max\left\{\frac{|\beta_0-1|}{\beta_0^2},
        \,\frac14\right\}\right]\\
    \end{split}
\end{equation}
and $c_1$ given by \eqref{beta.k.k-1} and we recall that  $\Lambda_p=2$ if $p\ne 1$ and $\Lambda_p= \lambda/4$ if $p=1$.  The proof is concluded once we let $\beta_0=2\overline{q}/2^*$ as in the previous step.\qed


\subsection{Local upper bounds II.  The linear case with unbounded coefficients}

The local upper bounds for nonnegative subsolutions to
\[
-\Delta u= b(x)\,u
\]
with $b\in \LL^r(B_R)$ eventually unbounded, follow from the local Sobolev imbedding theorem on balls $B_R\subset\RR^d$
\begin{equation}\label{Sobolev.BR.b}
\|f\|_{\LL^{2^*}(B_{R})}^2\le\mathcal{S}_2^2 \left(\|\nabla f\|_{\LL^2(B_{R})}^2 +\frac{1}{R^2}\|f\|_{\LL^2(B_{R})}^2\right)
\end{equation}
where $\mathcal{S}_2=\mathcal{S}_2(B_1)$ is the best constant and $2^*=2d/(d-2)$. In the case $f\in W^{1,2}_0(B_R)$, we have
\begin{equation}\label{Sobolev.BR.b.0}
\|f\|_{\LL^{2^*}(B_{R})}^2\le\mathcal{S}_2^2\,\|\nabla f\|_{\LL^2(B_{R})}^2.
\end{equation}
We are requiring hereafter without any further comment that $d\ge3$. We adopt the notation $\|f\|_{\LL^{q}(B_{R})}=\|f\|_{q,R}$ and $|B_R|=\omega_d R^d$.

\subsubsection{Energy Estimates and Reverse Poincar\'e inequalities}

\begin{lem}\label{lemma.Young}
Let $v\in L^{2^*}(B_R)$ and $b\in \LL^r(B_R)$ for some $r>d/2$. Then for any $\delta>0$ the following inequality holds
\begin{equation}\label{Sobolev.Improved}
\begin{split}
\int_{B_R} b(x)v^2(x)\dx\le \delta\left[\int_{B_R} v^{2^*}\dx\right]^{\frac{2}{2^*}}
    +\frac{K_{r,d}^{(1)}}{\delta^{\frac{d+r(d-2)}{2r-d}}} |B_R|^{\frac{2}{2^*}}
    \left[\int_{B_R} b^r(x)\dx\right]^{\frac{d}{2r-d}} \int_{B_R} v^2(x)\dx
\end{split}
\end{equation}
where
\begin{equation}\label{const.lemma.Young}
K_{r,d}^{(1)}:= \frac{2r-d}{rd}\left[\frac{rd}{d+r(d-2)}\right]^{\frac{d+r(d-2)}{2r-d}}
\end{equation}
\end{lem}
\noindent {\sl Proof.~}Let us estimate for any $0<\varepsilon<2$:
\[
\begin{split}
&\int_{B_R} bv^{(2-\varepsilon)+\varepsilon}\dx
\le_{(a)} \left[\int_{B_R} v^{(2-\varepsilon)\frac{2^*}{2}}\dx\right]^{\frac{2}{2^*}}
    \left[\int_{B_R} b^{\frac{d}{2}}v^{\varepsilon\frac{d}{2}}\dx\right]^{\frac{2}{d}}\\
&\le_{(b)} |B_R|^{\frac{\varepsilon}{2^*}}\left[\int_{B_R} v^{2^*}\dx\right]^{\frac{2-\varepsilon}{2^*}}
    \left[\int_{B_R} b^{\frac{d}{2}}v^{\varepsilon\frac{d}{2}}\dx\right]^{\frac{2}{d}}\\
&\le_{(c)}  \frac{\delta_0(2-\varepsilon)}{2} \left[\int_{B_R} v^{2^*}\dx\right]^{\frac{2}{2^*}}
    +\frac{\varepsilon}{2\delta_0^{\frac{2-\varepsilon}{\varepsilon}}} |B_R|^{\frac{2}{2^*}}
    \left[\int_{B_R} b^{\frac{d}{2}}v^{\varepsilon\frac{d}{2}}\dx\right]^{\frac{4}{d\varepsilon}}\\
&\le_{(d)}  \delta_0 \frac{d+r(d-2)}{rd}\left[\int_{B_R} v^{2^*}\dx\right]^{\frac{2}{2^*}}
    +\frac{2(2r-d)}{2rd\delta_0^{\frac{d+r(d-2)}{2r-d}}} |B_R|^{\frac{2}{2^*}}
    \left[\int_{B_R} b^r\dx\right]^{\frac{d}{2r-d}} \int_{B_R} v^2\dx\\
&\le_{(e)}  \delta\left[\int_{B_R} v^{2^*}\dx\right]^{\frac{2}{2^*}}
    +\frac{1}{\delta^{\frac{d+r(d-2)}{2r-d}}} \frac{2r-d}{rd}\left[\frac{rd}{d+r(d-2)}\right]^{\frac{d+r(d-2)}{2r-d}} |B_R|^{\frac{2}{2^*}}
    \left[\int_{B_R} b^r\dx\right]^{\frac{d}{2r-d}} \int_{B_R} v^2\dx
\end{split}
\]
where in the step $(a)$ we have used H\"older inequality with the conjugate exponents $s=2^*/2=d/(d-2)$ and $s'=s/(s-1)=d/2$. In $(b)$ we have used the inequality
\[
\left[\int_{B_R} v^{(2-\varepsilon)\frac{2^*}{2}}\dx\right]^{\frac{2}{2^*}}
\le |B_R|^{\frac{\varepsilon}{2^*}}\left[\int_{B_R} v^{2^*}\dx\right]^{\frac{2-\varepsilon}{2^*}}
\]
In $(c)$ we have applied the Young inequality, valid for every $\sigma>1$, $\delta_0>0$, $a,b\ge 0$:
$$
ab\leq \frac{\delta_0}{\sigma} a^\sigma+\frac{\sigma-1}\sigma\frac{b^{\frac\sigma{\sigma-1}}}{\delta_0^{\frac{1}{\sigma-1}}}\,
$$
with $\sigma=2/(2-\varepsilon)$, so that $\sigma/(\sigma-1)=2/\varepsilon$. In $(d)$ we have used the estimate
\[
\begin{split}
\left[\int_{B_R} b^{\frac{d}{2}}v^{\varepsilon\frac{d}{2}}\dx\right]^{\frac{4}{d\varepsilon}}
&\le \left[\int_{B_R} b^r\dx\right]^{\frac{d}{2r}\frac{4}{d\varepsilon}}
\left[\int_{B_R} v^{\varepsilon\frac{d}{2}\frac{2r}{2r-d}}\dx\right]^{\frac{2r-d}{2r}\frac{4}{d\varepsilon}}
= \left[\int_{B_R} b^r\dx\right]^{\frac{d}{2r-d}} \int_{B_R} v^2\dx
\end{split}
\]
where in the first step we have used H\"older inequality with the conjugate exponents $s=2r/d$ and $s'=s/(s-1)=2r/(2r-d)$ (notice that we are assuming $r>d/2$, hence $s>1$), while in the second step we have chosen $0<\varepsilon = 2(2r-d)/(rd)\le 2$. In $(e)$ we have put
\[
\delta= \delta_0 \frac{d+r(d-2)}{rd}\,
\]
notice that $\delta>0$ is in fact arbitrary since for every fixed $r$ we can choose appropriately $\delta_0$ to get any given value of $\delta$ by the above definition of $\delta$.\qed

\begin{thm}[Reverse Poincar\'e inequality for subsolutions]\label{Rev.Poinc}
Consider a weak subsolution $u$ to $-\Delta u= b(x)u$ on $B_R$\, with $b\in \LL^r({B_R})$ with $r>d/2$. Suppose that $u\in L^{\alpha+1}(B_R)$. Then for any positive test function $\varphi\in C^2_0(B_R)$ with $|\nabla \varphi|\equiv0$ on $\partial B_R$ we have that for any $R>0$ and $\alpha>0$:
\begin{equation}\label{Reverse.Poincare}
\int_{B_R}\big|\nabla u^{\frac{\alpha+1}{2}}\big|^2 \varphi^2 \dx \le K^{(2)}[b] \int_{B_R}u^{\alpha+1} \dx
\end{equation}
with
\begin{equation}\label{Reverse.Poincare.Lambda}\begin{split}
&K^{(2)}[b] =K^{(2)}(b,R,\alpha,\varphi,r,d) \\
&:=\frac{\alpha+1}{\alpha}
        \left[2\|\varphi\|_\infty\|\Delta\varphi\|_\infty+\|\nabla\varphi\|^2_\infty
    +\mathcal{S}_2^{\frac{2[d+r(d-2)]}{2r-d}}\left(\frac{(\alpha+1)^2}{2\alpha}\right)^{\frac{rd}{2r-d}}
        K_{r,d}\|\varphi\|_\infty^2|B_R|^{\frac{2}{2^*}}\|b\|_r^{\frac{dr}{2r-d}}\right].
\end{split}
\end{equation}
\end{thm}
{\bf Remark}. The requirement $u\in L^{1+\alpha}(B_R)$ will be dispensed with later, without further comment by using a Moser iteration technique.

\noindent {\sl Proof.~}It will divided into several steps.

\noindent$\bullet~$\textsc{Step 1. }\textit{Energy estimates. }Proceeding as in \eqref{local.energy.identity.subsols}, one shows that subsolutions to $-\Delta u\le b(x)u$, satisfy, even for any $\alpha\not=-1$:
\begin{equation}\label{Energy.b}
\frac{4\alpha}{(\alpha+1)^2}\int_{B_R}\big|\nabla u^{\frac{\alpha+1}{2}}\big|^2 \varphi^2 \dx\le \frac{1}{\alpha+1}\int_{B_R} u^{\alpha+1}\Delta\varphi^2\dx +\int_{B_R} b\,u^{\alpha+1}\varphi^2\dx.
\end{equation}
\noindent$\bullet~$\textsc{Step 2. }\textit{Sobolev inequality in $W_0^{1,2}(B_R)$. }We apply inequality \eqref{Sobolev.Improved} of Lemma \ref{lemma.Young} to $v=u^{(\alpha+1)/2}\varphi\in W^{1,2}_0(B_R)$ so that for any $\delta>0$:
\begin{equation}\label{Sobolev.Improved.2}
\begin{split}
\int_{B_R} bu^{\alpha+1}\varphi^2\dx\le \delta\left[\int_{B_R} \left(u^{\frac{\alpha+1}{2}}\varphi\right)^{2^*}\dx\right]^{\frac{2}{2^*}}
    +\frac{K_{r,d}}{\delta^{\frac{d+r(d-2)}{2r-d}}} |B_R|^{\frac{2}{2^*}}
    \left[\int_{B_R} b^r\dx\right]^{\frac{d}{2r-d}} \int_{B_R} u^{\alpha+1}\varphi^2\dx
\end{split}
\end{equation}
where $K_{r,d}$ is given in \eqref{const.lemma.Young}. We notice that $v=u^{(\alpha+1)/2}\varphi\in W^{1,2}_0(B_R)$, so that the Sobolev inequality \eqref{Sobolev.BR.b.0} reads
\[
\begin{split}
\left[\int_{B_R} \left(u^{\frac{\alpha+1}{2}}\varphi\right)^{2^*}\dx\right]^{\frac{2}{2^*}}
&\le \mathcal{S}_2^2 \int_{B_R}\big|\nabla u^{\frac{\alpha+1}{2}}\varphi\big|^2 \dx\\
&= \mathcal{S}_2^2 \left[\int_{B_R}\big|\nabla u^{\frac{\alpha+1}{2}}\big|^2 \varphi^2\dx
    + \int_{B_R}\big|\nabla \varphi \big|^2 u^{\alpha+1}\dx
    +\frac{1}{2}\int_{B_R}\nabla\varphi^2\cdot\nabla u^{\alpha+1}\dx \right]\\
&= \mathcal{S}_2^2 \left[\int_{B_R}\big|\nabla u^{\frac{\alpha+1}{2}}\big|^2 \varphi^2\dx
    - \int_{B_R}\varphi\big(\Delta\varphi\big) u^{\alpha+1}\dx\right]
\end{split}
\]
since $\Delta \varphi^2=2\varphi\Delta\varphi +2|\nabla\varphi|^2$. We combine the above Sobolev inequality with \eqref{Sobolev.Improved.2} to get
\begin{equation}\label{Sobolev.Improved.3}
\begin{split}
\int_{B_R} bu^{\alpha+1}\varphi^2\dx
&\le \delta\mathcal{S}_2^2 \left[\int_{B_R}\big|\nabla u^{\frac{\alpha+1}{2}}\big|^2 \varphi^2\dx
    - \int_{B_R}\varphi\big(\Delta\varphi\big) u^{\alpha+1}\dx\right]\\
&+\frac{K_{r,d}}{\delta^{\frac{d+r(d-2)}{2r-d}}} |B_R|^{\frac{2}{2^*}}
    \left[\int_{B_R} b^r\dx\right]^{\frac{d}{2r-d}} \int_{B_R} u^{\alpha+1}\varphi^2\dx
\end{split}
\end{equation}
where $K_{r,d}$ is given in \eqref{const.lemma.Young}.

\noindent$\bullet~$\textsc{Step 3. }Putting the pieces together, i.e. combining inequalities \eqref{Sobolev.Improved.3} and \eqref{Energy.b} we obtain
\[\begin{split}
\frac{4\alpha}{(\alpha+1)^2}\int_{B_R}\big|\nabla u^{\frac{\alpha+1}{2}}\big|^2 \varphi^2 \dx
&\le \frac{1}{\alpha+1}\int_{B_R} u^{\alpha+1}\Delta\varphi^2\dx +\int_{B_R} b\,u^{\alpha+1}\varphi^2\dx\\
&\le \frac{1}{\alpha+1}\int_{B_R} u^{\alpha+1}\Delta\varphi^2\dx\\
&+\delta\mathcal{S}_2^2 \left[\int_{B_R}\big|\nabla u^{\frac{\alpha+1}{2}}\big|^2 \varphi^2\dx
    - \int_{B_R}\varphi\big(\Delta\varphi\big) u^{\alpha+1}\dx\right]\\
&+\frac{K_{r,d}}{\delta^{\frac{d+r(d-2)}{2r-d}}} |B_R|^{\frac{2}{2^*}}
    \left[\int_{B_R} b^r\dx\right]^{\frac{d}{2r-d}} \int_{B_R} u^{\alpha+1}\varphi^2\dx
\end{split}\]
which thus implies
\[\begin{split}
&\left(\frac{4\alpha}{(\alpha+1)^2}-\delta\mathcal{S}_2^2\right)
 \int_{B_R}\big|\nabla u^{\frac{\alpha+1}{2}}\big|^2 \varphi^2 \dx\\
&\le \frac{1}{\alpha+1}\int_{B_R} u^{\alpha+1}\Delta\varphi^2\dx
        - \delta\mathcal{S}_2^2  \int_{B_R}\varphi\big(\Delta\varphi\big) u^{\alpha+1}\dx
    +\frac{K_{r,d}}{\delta^{\frac{d+r(d-2)}{2r-d}}} |B_R|^{\frac{2}{2^*}}
        \left[\int_{B_R} b^r\dx\right]^{\frac{d}{2r-d}} \int_{B_R} u^{\alpha+1}\varphi^2\dx\\
&\le \left[\left(\frac{2}{\alpha+1}+ \delta\mathcal{S}_2^2\right)\|\varphi\|_\infty\|\Delta\varphi\|_\infty+\frac{2}{\alpha+1}\|\nabla\varphi\|^2_\infty
    +\frac{K_{r,d}\|\varphi\|_\infty^2}{\delta^{\frac{rd}{2r-d}-1}} |B_R|^{\frac{2}{2^*}}
    \left(\int_{B_R} b^r\dx\right)^{\frac{d}{2r-d}}\right] \int_{B_R} u^{\alpha+1}\dx
\end{split}\]
Letting $\delta\mathcal{S}_2^2=\frac{2\alpha}{(\alpha+1)^2}$ gives the following reverse Poincar\'e inequality:
\[
\int_{B_R}\big|\nabla u^{\frac{\alpha+1}{2}}\big|^2 \varphi^2 \dx\le \Lambda_0\int_{B_R} u^{\alpha+1}\dx
\]
with the constant that we can estimate as follows
\[
\begin{split}
\Lambda_0
&=\frac{(\alpha+1)^2}{2\alpha}
        \left[\frac{2(2\alpha+1)}{(\alpha+1)^2}\|\varphi\|_\infty\|\Delta\varphi\|_\infty+\frac{2}{\alpha+1}\|\nabla\varphi\|^2_\infty\right.\\
    &+\left.\mathcal{S}_2^{\frac{2[d+r(d-2)]}{2r-d}}\frac{2\alpha}{(\alpha+1)^2}\left(\frac{(\alpha+1)^2}{2\alpha}\right)^{\frac{rd}{2r-d}}
        K_{r,d}\|\varphi\|_\infty^2|B_R|^{\frac{2}{2^*}}\|b\|_r^{\frac{dr}{2r-d}}\right]\\
&\le\frac{\alpha+1}{\alpha}
        \left[2\|\varphi\|_\infty\|\Delta\varphi\|_\infty+\|\nabla\varphi\|^2_\infty
    +\mathcal{S}_2^{\frac{2[d+r(d-2)]}{2r-d}}\left(\frac{(\alpha+1)^2}{2\alpha}\right)^{\frac{rd}{2r-d}}
        K_{r,d}\|\varphi\|_\infty^2|B_R|^{\frac{2}{2^*}}\|b\|_r^{\frac{dr}{2r-d}}\right]\\
        &\le\frac{\big(\alpha+1\big)^{1+\frac{rd}{2r-d}}}{\alpha}
         \left[2\|\varphi\|_\infty\|\Delta\varphi\|_\infty+\|\nabla\varphi\|^2_\infty
    +\mathcal{S}_2^{\frac{2[d+r(d-2)]}{2r-d}}\left(\frac{\alpha+1}{2\alpha}\right)^{\frac{rd}{2r-d}}
        K_{r,d}\|\varphi\|_\infty^2|B_R|^{\frac{2}{2^*}}\|b\|_r^{\frac{dr}{2r-d}}\right]=:K^{(2)}[b].\mbox{\qed}
\end{split}
\]
In fact, the last bound in the above formula for $K^{(2)}[b]$ could be avoided, but will make the following calculations somewhat easier.

\begin{lem}[Numerical Iteration]\label{lem.num.iter}
Let $Y_n\ge 0$ be a sequence of numbers such that
\begin{equation}\label{hyp.num}
Y_n\le I_{n-1}^{\sigma \theta^{n-1}}Y_{n-1}\qquad\mbox{with}\qquad I_{n-1}\le I_0 C^{\,n-1}\,
\end{equation}
for some $\sigma,I_0,C>0$, $\theta\in (0,1)$. Then $\{Y_n\}$ is a bounded sequence and one has
\begin{equation}\label{iteration.num}
Y_\infty:=\limsup_{n\to+\infty}Y_n\le I_0^{\frac{\sigma}{1-\theta}}\,C^{\frac{\sigma\,\theta}{(1-\theta)^2}}\,Y_0.
\end{equation}
\end{lem}
\noindent {\sl Proof.~}We iterate inequality \eqref{hyp.num} to get
\[\begin{split}
Y_n & \le I_{n-1}^{\sigma \theta^{n-1}}Y_{n-1}\le \left(I_0 C^{\,n-1}\right)^{\sigma \theta^{n-1}} Y_{n-1}
=I_0^{\sigma \theta^{n-1}} C^{\sigma (n-1)\theta^{n-1}} Y_{n-1}
\le \prod_{j=0}^{n-1}I_0^{\sigma \theta^j} C^{\sigma j\theta^j}\,Y_0\\
&= I_0^{\sigma \sum_{j=0}^{n-1}\theta^j} C^{\sigma \sum_{j=0}^{n-1}j\theta^j}\,Y_0\\
\end{split}
\]
We thus get, as $n\to+\infty$, $Y_\infty \le I_0^{\sigma \sum_{j=0}^{\infty}\theta^j} C^{\sigma \sum_{j=0}^{\infty}j\theta^j}\,Y_0
=I_0^{\frac{\sigma}{1-\theta}}\,C^{\frac{\sigma\,\theta}{(1-\theta)^2}}\,Y_0.\mbox{\qed}$

Now we are ready to perform the Moser iteration, by combining a local Sobolev inequality with the reverse Poincar\'e inequality of Theorem \ref{Rev.Poinc} and then using the above numerical Lemma.
\begin{thm}[Moser Iteration]\label{Moser.Upper.b.thm}
Let $u\ge 0$ be a weak subsolution  to $-\Delta u=b\,u$ on $B_R$\, with $b\in \LL^r({B_R})$ with $r>d/2$, and let $q>1$, $R_\infty<R_0<R$.
\begin{equation}\label{Moser.Upper.b}
\|u\|_{\infty, R_\infty}\le \frac{K^{(3)}_q[b]}{{(R_0-R_\infty)^{\frac{d}{q}}}}\,\|u\|_{q,R_0}
\end{equation}
with constant
\begin{equation}\label{Const.Upper.b}
\begin{split}
K^{(3)}_q[b]&=\left(\frac{qd^d}{2^d}\right)^{\frac{rd^2}{2(2r-d)q}}
         \left[8\frac{q(d+2)}{q-1}
        +\left(\frac{\mathcal{S}_2^2}{2}\right)^{\frac{rd}{2r-d}}\frac{2r-d}{rd}
        \left(\frac{q rd}{(q-1)[d+r(d-2)]}\right)^{1+\frac{rd}{2r-d}}\right.\\
    &\left.\times (R_0-R_\infty)^2|B_{R_0} |^{\frac{2}{2^*}}\|b\|_{\LL^r({B_{R_0} })}^{\frac{rd}{2r-d}}
         +\left(\frac{R_0-R_\infty}{R_\infty}\right)^2\right]^{\frac{d}{2q}}.
\end{split}
\end{equation}
\end{thm}
Notice that in the case of bounded coefficients $b(x)\in \LL^\infty(B_{R_0})$ we can pass to the limit as $r\to\infty$ in the above expression of $K^{(3)}_q[b]$ to get
\begin{equation}\label{Const.Upper.b.infty}
\begin{split}
K^{(3)}_q[b]&=\left(\frac{qd^d}{2^d}\right)^{\frac{d^2}{4q}}
         \left[8\frac{q(d+2)}{q-1}
        +\left(\frac{\mathcal{S}_2^2}{2}\right)^{\frac{d}{2}}\frac{2}{d}
        \left(\frac{qd}{(q-1)(d-2)}\right)^{1+\frac{d}{2}}\right.\\
    &\left.\times (R_0-R_\infty)^2|B_{R_0} |^{\frac{2}{2^*}}\|b\|_{\LL^\infty({B_{R_0} })}^{\frac{d}{2}}
         +\left(\frac{R_0-R_\infty}{R_\infty}\right)^2\right]^{\frac{d}{2q}}.
\end{split}
\end{equation}

\noindent {\sl Proof.~}The proof is divided in several steps.

\noindent$\bullet~$\textsc{Step 1. }\textit{Sobolev and Reverse Poincar\'e inequalities. }We start choosing radii $r_1,r_0$ with $R_\infty<r_1<r_0<R_0$ and use the test function of Lemma \ref{lem.test.funct} on the balls $B_{r_1}, B_{r_0}$. We use the Reverse Poincar\'e inequality \eqref{Reverse.Poincare} on the ball $B_{r_0}$ and the fact that $\varphi\equiv 1$ on $B_{r_1}$ to get
\begin{equation*}
\int_{B_{r_1}}\big|\nabla u^{\frac{\alpha+1}{2}}\big|^2 \dx\le \int_{B_{r_0}}\big|\nabla u^{\frac{\alpha+1}{2}}\big|^2 \varphi^2 \dx \le K^{(2)}[b] \int_{B_{r_0}}u^{\alpha+1} \dx
\end{equation*}
so that the local Sobolev inequality in $W^{1,2}(B_{r_1})$ applied to $f=v^{\frac{\alpha+1}{2}}$ for any $\alpha>0$ yields
\begin{equation}\label{iteration.1.b}
\begin{split}
\left(\int_{B_{r_1}}u^{\frac{2^*}{2}(\alpha+1)} \dx\right)^{\frac{2}{2^*}}
    &\le \mathcal{S}_2^2 \left[\int_{B_{r_1}}\big|\nabla u^{\frac{\alpha+1}{2}}\big|^2 \dx +\frac{1}{r_1^2}\int_{B_{r_1}}u^{\alpha+1} \dx\right]\\
    &\le\mathcal{S}_2^2\left(K^{(2)}[b]+\frac{1}{r_1^2}\right)\int_{B_{r_0}}u^{\alpha+1} \dx\\
\end{split}
\end{equation}
where the constant $K^{(2)}[b]$ is given by \eqref{Reverse.Poincare}, and we can estimate it as follows:
\begin{equation*}\begin{split}
K^{(2)}[b]&= \frac{\big(\alpha+1\big)^{1+\frac{rd}{2r-d}}}{\alpha}
         \Bigg[2\|\varphi\|_{\infty,r_0 }\|\Delta\varphi\|_{\infty,r_0 }+\|\nabla\varphi\|^2_{\infty,r_0 }\Bigg.\\
 & \left.+\mathcal{S}_2^{\frac{2[d+r(d-2)]}{2r-d}}\left(\frac{\alpha+1}{2\alpha}\right)^{\frac{rd}{2r-d}}
        \frac{2r-d}{rd}\left(\frac{rd}{d+r(d-2)}\right)^{1+\frac{rd}{2r-d}}\|\varphi\|_{\infty,r_0 }^2
        |B_{r_0} |^{\frac{2}{2^*}}\|b\|_{\LL^r({B_{r_0} })}^{\frac{rd}{2r-d}}\right]\\
 &\le_{(a)} \frac{\alpha+1}{\alpha}\big(\alpha+1\big)^{\frac{rd}{2r-d}}
         \left[\frac{8(d+2)}{(r_0-r_1)^2}
        +\mathcal{S}_2^{\frac{2[d+r(d-2)]}{2r-d}}\left(\frac{\alpha+1}{2\alpha}\right)^{\frac{rd}{2r-d}}\frac{2r-d}{rd}\right.\\
        &\left.\times\left(\frac{rd}{d+r(d-2)}\right)^{1+\frac{rd}{2r-d}}
         |B_{r_0} |^{\frac{2}{2^*}}\|b\|_{\LL^r({B_{r_0} })}^{\frac{rd}{2r-d}}\right]\\
 &\le_{(b)} \frac{\big(\alpha+1\big)^{\frac{rd}{2r-d}}}{{(r_0-r_1)^2}}
         \left[8(d+2)\frac{\alpha+1}{\alpha}
        +\frac{{\mathcal{S}_2}^{\frac{2[d+r(d-2)]}{2r-d}}}{2^\frac{rd}{2r-d}}\frac{2r-d}{rd}\right.\\
        &\left.\times\left(\frac{(\alpha+1)rd}{\alpha[d+r(d-2)]}\right)^{1+\frac{rd}{2r-d}}
         (R_0-R_\infty)^2|B_{R_0} |^{\frac{2}{2^*}}\|b\|_{\LL^r({B_{R_0} })}^{\frac{rd}{2r-d}}\right]\\
\end{split}
\end{equation*}
where in $(a)$ we have used the fact that the test function of Lemma \ref{lem.test.funct} satisfies $\|\varphi\|_{\infty,r_0}= 1$\, $\|\nabla\varphi\|_{\infty,r_0}\le 4/(r_0-r_1)$ and $\|\Delta\varphi\|_{\infty,r_0}\le 4d/(r_0-r_1)^2$, and in $(b)$ the fact that $0<R_\infty<r_1<r_0<R_0$. Finally we get:
\begin{equation}\begin{split}
\mathcal{S}_2^2 &\left(K^{(2)}[b]+\frac{1}{r_1^2}\right)
 \le \mathcal{S}_2^2\frac{\big(\alpha+1\big)^{\frac{rd}{2r-d}}}{{(r_0-r_1)^2}}
         \left[8(d+2)\frac{\alpha+1}{\alpha}
        +\frac{{\mathcal{S}_2}^{\frac{2[d+r(d-2)]}{2r-d}}}{2^\frac{rd}{2r-d}}\frac{2r-d}{rd}\right.\\
        &\left.\times\left(\frac{(\alpha+1)rd}{\alpha[d+r(d-2)]}\right)^{1+\frac{rd}{2r-d}}
         (R_0-R_\infty)^2|B_{R_0} |^{\frac{2}{2^*}}\|b\|_{\LL^r({B_{R_0} })}^{\frac{rd}{2r-d}}
         +\frac{1}{r_1^2}\frac{{(r_0-r_1)^2}}{\big(\alpha+1\big)^{\frac{rd}{2r-d}}}\right]\\
 &\le \mathcal{S}_2^2\frac{\big(\alpha+1\big)^{\frac{rd}{2r-d}}}{{(r_0-r_1)^2}}
         \left[8(d+2)\frac{\alpha+1}{\alpha}
        +\frac{{\mathcal{S}_2}^{\frac{2[d+r(d-2)]}{2r-d}}}{2^\frac{rd}{2r-d}}\frac{2r-d}{rd}\right.\\
        &\left.\times\left(\frac{(\alpha+1)rd}{\alpha[d+r(d-2)]}\right)^{1+\frac{rd}{2r-d}}
         (R_0-R_\infty)^2|B_{R_0} |^{\frac{2}{2^*}}\|b\|_{\LL^r({B_{R_0} })}^{\frac{rd}{2r-d}}
         +\left(\frac{R_0-R_\infty}{R_\infty}\right)^2\right]\\
\end{split}
\end{equation}
we have also used the fact that $\alpha>0$.

\noindent$\bullet~$\textsc{Step 2. }\textit{The Moser iteration. }We now fix $\beta_0=\alpha+1>1$, and we define the sequence
\[
\beta_n=\frac{2^*}{2}\beta_{n-1}=\left(\frac{2^*}{2}\right)^n\beta_0\,
\]
Next we pick a sequence of radii $R_\infty=r_\infty<\ldots<r_n<r_{n-1}<\ldots<r_0=R_0$, such that
\[
(r_{n-1}-r_n)^2=c_0^2\,(R_0-R_\infty)^2\,\left(\frac{2}{2^*}\right)^{\frac{rd\,n}{2r-d}}
\]
with
\begin{equation}\label{c000}
c_0=\left(\sum_{k=1}^\infty \left(\frac{2}{2^*}\right)^{\frac{rd}{2(2r-d)}k}\right)^{-1}
=\left(\frac{2^*}{2}\right)^{\frac{rd}{2(2r-d)}}-1\ge\left(\frac{2^*}{2}-1\right)^{\frac{rd}{2(2r-d)}}
=\left(\frac{2}{d-2}\right)^{\frac{rd}{2(2r-d)}}
\end{equation}
where the inequality in the above formula is easily shown to hold when $d\ge3$ and $r>d/2$ as assumed, so that
\[
\sum_{k=1}^\infty (r_{k-1}-r_k)=R_0-R_\infty,
\]
the above series being convergent. With these choices, inequality \eqref{iteration.1.b}  in which $\alpha+1$ is replaced by $\beta_{n-1}$, this being allowable since $\beta_n>1$ for all $n$, and $r_1,r_0$ replaced by $r_{n}, r_{n-1}$ reads, noticing in addition that $\beta_n/(\beta_n-1)\le\beta_0/(\beta_0-1)$ for all $n$,
\[
\begin{split}
&\left(\int_{B_{r_n}}u^{\frac{2^*}{2}\beta_{n-1}} \dx\right)^{\frac{2}{2^*}}
    \le \mathcal{S}_2^2\left(K^{(2)}[b]+\frac{1}{r_n^2}\right)\int_{B_{r_{n-1}}}u^{\beta_{n-1}} \dx\\
    &\le \frac{\mathcal{S}_2^2\beta_{n-1}^{\frac{rd}{2r-d}}}{{(r_{n-1}-r_n)^2}}
         \left[8(d+2)\frac{\beta_0}{\beta_0-1}
        +\frac{{\mathcal{S}_2}^{\frac{2[d+r(d-2)]}{2r-d}}}{2^\frac{rd}{2r-d}}\frac{2r-d}{rd}
        \left(\frac{\beta_0 rd}{(\beta_0-1)[d+r(d-2)]}\right)^{1+\frac{rd}{2r-d}}\right.\\
    &\left.\times (R_0-R_\infty)^2|B_{R_0} |^{\frac{2}{2^*}}\|b\|_{\LL^r({B_{R_0} })}^{\frac{rd}{2r-d}}
         +\left(\frac{R_0-R_\infty}{R_\infty}\right)^2\right]\int_{B_{r_{n-1}}}u^{\beta_{n-1}} \dx :=I_{n-1}\int_{B_{r_{n-1}}}u^{\beta_{n-1}}\dx
\end{split}
\]
Letting $Y_n:=\|u\|_{\beta_n, R_n}$, we have obtained
\[
Y_n=\|u\|_{\beta_n, R_n}\le I_{n-1}^{\frac{1}{\beta_{n-1}}}\|u\|_{\beta_{n-1}, R_{n-1}}=I_{n-1}^{\frac{1}{\beta_{n-1}}}Y_{n-1}
    =I_{n-1}^{\frac{1}{\beta_0}\left(\frac{2}{2^*}\right)^{n-1}}Y_{n-1}=I_{n-1}^{\sigma\theta^{n-1}}Y_{n-1}
\]
where we have set $\sigma=1/\beta_0$ and $\theta=2/2^*\in(0,1)$. We shall prove that $I_n\le I_0 C^n$. Indeed:
\begin{equation}\label{iteration.n.b}
\begin{split}
I_{n-1}
    &= \frac{\mathcal{S}_2^2\beta_{n-1}^{\frac{rd}{2r-d}}}{{(r_{n-1}-r_n)^2}}
         \left[8(d+2)\frac{\beta_0}{\beta_0-1}
        +\frac{{\mathcal{S}_2}^{\frac{2[d+r(d-2)]}{2r-d}}}{2^\frac{rd}{2r-d}}\frac{2r-d}{rd}
        \left(\frac{\beta_0 rd}{(\beta_0-1)[d+r(d-2)]}\right)^{1+\frac{rd}{2r-d}}\right.\\
    &\left.\times (R_0-R_\infty)^2|B_{R_0} |^{\frac{2}{2^*}}\|b\|_{\LL^r({B_{R_0} })}^{\frac{rd}{2r-d}}
         +\left(\frac{R_0-R_\infty}{R_\infty}\right)^2\right]\\
    &\le \frac{\beta_{0}^{\frac{rd}{2r-d}}}{{c_0^2(R_0-R_\infty)^2}}
         \left[8(d+2)\frac{\beta_0}{\beta_0-1}
        +\frac{{\mathcal{S}_2}^{\frac{2[d+r(d-2)]}{2r-d}}}{2^\frac{rd}{2r-d}}\frac{2r-d}{rd}
        \left(\frac{\beta_0 rd}{(\beta_0-1)[d+r(d-2)]}\right)^{1+\frac{rd}{2r-d}}\right.\\
    &\left.\times (R_0-R_\infty)^2|B_{R_0} |^{\frac{2}{2^*}}\|b\|_{\LL^r({B_{R_0} })}^{\frac{rd}{2r-d}}
         +\left(\frac{R_0-R_\infty}{R_\infty}\right)^2\right]\left(\frac{2^*}{2}\right)^{\frac{2rd\,n}{2r-d}}\\
    &\le \left(\frac{d-2}{2}\right)^{\frac{rd}{2r-d}}
        \frac{\beta_{0}^{\frac{rd}{2r-d}}}{{(R_0-R_\infty)^2}}
         \left[8(d+2)\frac{\beta_0}{\beta_0-1}
        +\frac{{\mathcal{S}_2}^{\frac{2[d+r(d-2)]}{2r-d}}}{2^\frac{rd}{2r-d}}
        \left(\frac{\beta_0 rd}{(\beta_0-1)[d+r(d-2)]}\right)^{1+\frac{rd}{2r-d}}\right.\\
    &\left.\times \frac{2r-d}{rd}(R_0-R_\infty)^2|B_{R_0} |^{\frac{2}{2^*}}\|b\|_{\LL^r({B_{R_0} })}^{\frac{rd}{2r-d}}
         +\left(\frac{R_0-R_\infty}{R_\infty}\right)^2\right]\left(\frac{2^*}{2}\right)^{\frac{2rd}{2r-d}n}:=I_0 C^{n-1}\\
\end{split}
\end{equation}
where in the last inequality we estimated $c_0$ as in \eqref{c000}.  Finally we use Lemma \ref{lem.num.iter} with the above choices of $\sigma$ and $\theta$, thus proving that
\begin{equation*}
Y_\infty\le I_0^{\frac{\sigma}{1-\theta}}\,C^{\frac{\sigma\,\theta}{(1-\theta)^2}}\,Y_0
\qquad\mbox{which is}\qquad
\|u\|_{\infty, R_\infty}\le I_0^{\frac{d}{2\beta_0}}\,C^{\frac{d(d-2)}{4\beta_0}}\,\|u\|_{\beta_0,R_0}=K^{(3)}_q[b]\,\|u\|_{\beta_0,R_0}
\end{equation*}
which is exactly \eqref{Moser.Upper.b} with
\[\begin{split}
K^{(3)}_q[b]&=\left(\frac{d-2}{2}\right)^{\frac{rd^2}{2(2r-d)\beta_0}}
        \frac{\beta_{0}^{\frac{rd^2}{2(2r-d)\beta_0}}}{(R_0-R_\infty)^{\frac{d}{\beta_0}}}
         \left[8(d+2)\frac{\beta_0}{\beta_0-1}
        +\frac{{\mathcal{S}_2}^{\frac{2[d+r(d-2)]}{2r-d}}}{2^\frac{rd}{2r-d}}
        \left(\frac{\beta_0 rd}{(\beta_0-1)[d+r(d-2)]}\right)^{1+\frac{rd}{2r-d}}\right.\\
    &\left.\times \frac{2r-d}{rd}(R_0-R_\infty)^2|B_{R_0} |^{\frac{2}{2^*}}\|b\|_{\LL^r({B_{R_0} })}^{\frac{rd}{2r-d}}
         +\left(\frac{R_0-R_\infty}{R_\infty}\right)^2\right]^{\frac{d}{2\beta_0}}
         \left(\frac{d}{d-2}\right)^{\frac{rd^2}{(2r-d)\beta_0}}\left(\frac{d}{d-2}\right)^{\frac{rd^2(d-2)}{2\beta_0(2r-d)}}\\
&\le \left(\frac{d}{2}\right)^{\frac{rd^3}{2(2r-d)\beta_0}}
        \frac{\beta_{0}^{\frac{rd^2}{2(2r-d)\beta_0}}}{(R_0-R_\infty)^{\frac{d}{\beta_0}}}
         \left[8(d+2)\frac{\beta_0}{\beta_0-1}
        +\frac{{\mathcal{S}_2}^{\frac{2[d+r(d-2)]}{2r-d}}}{2^\frac{rd}{2r-d}}
        \left(\frac{\beta_0 rd}{(\beta_0-1)[d+r(d-2)]}\right)^{1+\frac{rd}{2r-d}}\right.\\
    &\left.\times \frac{2r-d}{rd}(R_0-R_\infty)^2|B_{R_0} |^{\frac{2}{2^*}}\|b\|_{\LL^r({B_{R_0} })}^{\frac{rd}{2r-d}}
         +\left(\frac{R_0-R_\infty}{R_\infty}\right)^2\right]^{\frac{d}{2\beta_0}}\,,
\end{split}\]
as in \eqref{Const.Upper.b}. The proof is concluded once we let $\beta_0=q>1$.\qed

\subsubsection{Extending local upper bounds. A lemma by E. De Giorgi}

In this section extend the local upper bound of the previous section. More precisely we show that a bound of the type
\[
\|u\|_{\infty, r}\le \frac{A}{(R-r)^{\frac{d}{q}}}\,\|u\|_{q,R}
\]
which is valid for any $q>1$\, and any $a\le r<R\le b$\, indeed implies that
\[
\|u\|_{\infty, a}\le \frac{A}{(b-a)^{\frac{d}{q_0}}}\,\|u\|_{q_0,b}
\]
for all $q_0>0$ and and any $a\le r<R\le b$\, maybe with a different constant $A$. The proof relies on the following lemma, originally due to E. De Giorgi, whose proof is contained in several books and papers, see for example \cite{Giusti}, Lemma 6.1.
\begin{lem}[De Giorgi]\label{DG.Lem}
Let $Z(t)$ be a bounded non-negative function in the interval $[t_0,t_1]$. Assume that for $t_0\le t<s\le t_1$ we have
\begin{equation}\label{hyp.DG}
Z(t)\le \theta\,Z(s)+\frac{A}{(s-t)^{\alpha}}
\end{equation}
with $A\ge 0$, $\alpha>0$ and $0\le \theta<1$. Then
\begin{equation}\label{ts.DG}
Z(t_0)\le \frac{A\,c(\alpha,\lambda,\theta)}{(t_1-t_0)^{\alpha}}
\end{equation}
where
\[
c(\alpha,\lambda,\theta)=\frac{1}{(1-\lambda)^\alpha\left(1-\frac{\theta}{\lambda^\alpha}\right)}\,\qquad\mbox{for any}\qquad \lambda\in\left(\theta^{\frac{1}{\alpha}},1\right).
\]
\end{lem}
\noindent {\sl Proof.~}Consider the sequence $\{s_i\}$ defined by
\[
t_0=s_0,\ \ s_{i+1}=s_i+(1-\lambda)\lambda^i(t_1-t_0)
\]
so that $s_k=t_0+(1-\lambda)(t_1-t_0)\sum_{i=0}^{k-1}\lambda^i$ and in particular $s_k\uparrow t_1$ as $k\to+\infty$. From \eqref{hyp.DG}, by induction we get
\[
Z(t_0)\le \theta^k Z(s_k)
    +\frac{A}{(1-\lambda)^\alpha(t_1-t_0)^{\alpha}}\sum_{i=0}^{k-1}\left[\frac{\theta}{\lambda^\alpha}\right]^i
\]
Since $\theta/\lambda^\alpha<1$ by assumption, the series on the right-hand side converges and therefore taking the limit as $k\to\infty$ and using the boundedness of $Z$ we get \eqref{ts.DG}.\qed

The above Lemma has important consequences, indeed it allows to prove that if a reverse H\"older inequality holds for some $0<q<\overline{q}$, then it holds for any $0<q_0<\overline{q}$.

\begin{lem}[Extending Local Upper Bounds]\label{Lem.DG.Ext}Assume that the following bounds holds true:
\begin{equation}\label{hip.DG.Ext}
\|u\|_{\overline{q}, r}\le \frac{K}{(R-r)^\gamma}\,\|u\|_{\underline{q},R}
\end{equation}
for some $0<\underline{q}<\overline{q}$\,, $\gamma>0$ and for any $R_\infty\le r<R\le R_0$. Then we have that for all $0<q_0\le \underline{q}<\overline{q}$
\begin{equation}\label{DG.Ext.bdd}
\|u\|_{\overline{q}, R_\infty}\le 3\cdot 2^{\frac{\overline{q}(\underline{q}-q_0)}{q_0(\overline{q}-\underline{q})}}\,
        \left[\left(4\gamma\frac{\underline{q}(\overline{q}-q_0)}{q_0(\overline{q}-\underline{q})}
        \right)^\gamma
        \frac{K}{(R_0-R_\infty)^\gamma}\right]^{\frac{\underline{q}(\overline{q}-q_0)}{q_0(\overline{q}-\underline{q})}}
            \,\|u\|_{q_0,R_0}\,.
\end{equation}
\end{lem}
\noindent {\sl Proof.~}Define, for $t<R_0$, the bounded nonnegative function
\[
Z(t)=\|u\|_{\LL^{\overline{q}}(B_{t})}=\|u\|_{\overline{q},t}
\]
then \eqref{hip.DG.Ext} reads, for $s>t$,
\begin{equation}\label{DG33}
Z(t)=\|u\|_{\overline{q}, t}\le \frac{K}{(s-t)^\gamma}\,\|u\|_{\underline{q},s}\le \frac{K}{(s-t)^\gamma}\,\|u\|_{q_0,s}^{1-\sigma}
\|u\|_{\overline{q},s}^{\sigma}\,,
\end{equation}
where in the last step we have used that for all $0<q_0\le \underline{q}< \overline{q}\le +\infty$
\[
\|u\|_{\underline{q},s}
\le \|u\|_{q_0,s}^{1-\sigma}\|u\|_{\overline{q},s}^{\sigma}
=\|u\|_{q_0,s}^{\frac{q_0(\overline{q}-\underline{q})}{\underline{q}(\overline{q}-q_0)}}
\|u\|_{\overline{q},s}^{\frac{\overline{q}(\underline{q}-q_0)}{\underline{q}(\overline{q}-q_0)}}\,,\qquad\mbox{with }\sigma
=\frac{\overline{q}(\underline{q}-q_0)}{\underline{q}(\overline{q}-q_0)}\in [0,1)
\]
Inequality \eqref{DG33} gives then
\begin{equation}\label{DG44}\begin{aligned}
Z(t)&=\|u\|_{\overline{q}, t}
\le \frac{K}{(s-t)^\gamma}\,\|u\|_{q_0,s}^{1-\sigma} Z(s)^{\sigma}
\le \frac{1}{2}Z(s) + \frac{\left(2^\sigma K\right)^{\frac{1}{1-\sigma}}}{(s-t)^{\frac{\gamma}{1-\sigma}}}\,\|u\|_{q_0,s}\\
&\le \frac{1}{2}Z(s) + \frac{\left(2^\sigma K\right)^{\frac{1}{1-\sigma}}}{(s-t)^{\frac{\gamma}{1-\sigma}}}\,\|u\|_{q_0,R_0}
\end{aligned}
\end{equation}
where we have used Young's inequality valid for any $\nu>1$, $a,b\ge 0$, $\varepsilon>0$:
\[
ab\le \frac{\varepsilon}{\nu} a^{\nu}+\frac{\nu-1}{\nu}\frac{b^{\frac{\nu}{\nu-1}}}{\varepsilon^{\frac{1}{\nu-1}}}
\le \varepsilon a^{\nu}+\frac{b^{\frac{\nu}{\nu-1}}}{\varepsilon^{\frac{1}{\nu-1}}}
\]
with the choices
\[
\varepsilon=1/2\,\qquad a=Z(s)^{\sigma}, \qquad \nu= \frac{1}{\sigma}>1\qquad\mbox{and}\qquad b=\frac{K}{(s-t)^\gamma}\|u\|^{1-\sigma}_{q_0, R_0}.
\]
Inequality \eqref{DG44} is of the form appearing in Lemma \ref{DG.Lem} with $\alpha=\gamma/(1-\sigma)>0$, $\theta= 1/2$ and $A=\left(2^\sigma K\right)^{\frac{1}{1-\sigma}}\,\|u\|_{q_0,R_0}$. Thus we get
\[\begin{split}
\|u\|_{\overline{q},R_\infty}=Z(R_\infty) &\le
\frac{c(\alpha,\lambda,\theta)\left(2^\sigma K\right)^{\frac{1}{1-\sigma}}}{(R_0-R_\infty)^{\frac{\gamma}{1-\sigma}}}\,\|u\|_{q_0,R_0}
    \le 3\,\left(\frac{4\gamma}{1-\sigma}\right)^{\frac{\gamma}{1-\sigma}} \frac{\left(2^\sigma K
        \right)^{\frac{1}{1-\sigma}}}{(R_0-R_\infty)^{\frac{\gamma}{1-\sigma}}}\,\|u\|_{q_0,R_0}\\
    &= 3\cdot 2^{\frac{\overline{q}(\underline{q}-q_0)}{q_0(\overline{q}-\underline{q})}}\,
        \left[\left(4\gamma\frac{\underline{q}(\overline{q}-q_0)}{q_0(\overline{q}-\underline{q})}
        \right)^\gamma
        \frac{K}{(R_0-R_\infty)^\gamma}\right]^{\frac{\underline{q}(\overline{q}-q_0)}{q_0(\overline{q}-\underline{q})}}
            \,\|u\|_{q_0,R_0}
\end{split}
\]
noticing that
\[
\frac{\sigma}{1-\sigma}=\frac{\overline{q}(\underline{q}-q_0)}{q_0(\overline{q}-\underline{q})}\,,\qquad\mbox{and}\qquad
\alpha=\frac{\gamma}{1-\sigma}=\gamma\frac{\underline{q}(\overline{q}-q_0)}{q_0(\overline{q}-\underline{q})}
\]

which is the desired bound, once we notice that whenever $\theta<\lambda^\alpha<1$,
\[
c(\alpha,\lambda,\theta)=\frac{1}{(1-\lambda)^\alpha\left(1-\frac{\theta}{\lambda^\alpha}\right)}
=\frac{2(1+\theta)}{\left[2^{\frac{1}{\alpha}}-(1+\theta)^{\frac{1}{\alpha}}\right]^\alpha(1-\theta)}
=\frac{12}{\left[4^{\frac{1}{\alpha}}-3^{\frac{1}{\alpha}}\right]^\alpha}\le 12\frac{4^{\alpha}\,\alpha^\alpha}{4}=3\,(4\alpha)^\alpha
\]
since we can choose $1/2=\theta<\lambda^\alpha=(1+\theta)/2<1$, and since $\alpha=\gamma/(1-\sigma)>1$\,,
\[
\left(4^{1/\alpha}-3^{1/\alpha}\right)^\alpha \ge \frac{4}{4^{\alpha}\,\alpha^\alpha}
\]
since we know that $a^{1/\alpha}-b^{1/\alpha}\ge a^{1/\alpha}(a-b)/(\alpha\,a)$, for all $a\ge b\ge 0$ and $\alpha\ge 1$.\qed

The above lemma can be used to extend the local upper bounds \eqref{Moser.Upper.b.DG} of Theorem \ref{Moser.Upper.b.thm}.

\begin{thm}[Local Upper bounds, unbounded coefficient]\label{Thm.local.upper.b.DG}
Consider a weak subsolution $u$ to $-\Delta u= b\,u$ on $B_R$\, with $b\in \LL^r({B_R})$ with $r>d/2$ Let $0<R_\infty<R_0<R$. Then for any $q_0>0$\,, the following bound holds true
\begin{equation}\label{Moser.Upper.b.DG}
\|u\|_{\infty,R_\infty}\le
    \frac{A_{q_0}^{(1)}}{(R-R_\infty)^{\frac{d}{q_0}}}
    \left[A_{q_0}^{(2)}+A_{q_0}^{(3)}\|b\|_{\LL^r({B_{R_0} })}^{\frac{rd}{2r-d}}\right]^{\frac{d}{2q_0}}\;\|u\|_{q_0,R}
\end{equation}
with
\begin{equation}\label{A1.b}
A_{q_0}^{(1)}:=\left\{
\begin{array}{lll}
\left(\frac{q_0d^d}{2^d}\right)^{\frac{rd^2}{2(2r-d)q_0}}\,,\quad &\mbox{ if }q_0>1\,,\\
3\cdot 2^{\frac{2d+1}{q_0}}\,\left(\frac{d}{q_0}\right)^{\frac{d}{q_0}}
       \left(\frac{(q_0+1)d^d}{2^d}\right)^{\frac{rd^2}{2(2r-d)q_0}}\,,\quad &\mbox{ if }0< q_0 \le 1\,,\\
\end{array}
\right.
\end{equation}
\begin{equation}\label{A2.b}
A_{q_0}^{(2)}:=\left\{
\begin{array}{lll}
8\frac{q_0(d+2)}{q_0-1}+\left(\frac{R -R_\infty}{R_\infty}\right)^2\,,\quad &\mbox{ if }q_0>1\,,\\
8\frac{(q_0+1)(d+2)}{q_0}+\left(\frac{R -R_\infty}{R_\infty}\right)^2\,,\quad &\mbox{ if }0< q_0 \le 1\,,\\
\end{array}
\right.
\end{equation}
\begin{equation}\label{A3.b}
A_{q_0}^{(3)}:=\left\{
\begin{array}{lll}
\left(\frac{\mathcal{S}_2^2}{2}\right)^{\frac{rd}{2r-d}}\frac{2r-d}{rd}
        \left(\frac{q_0 rd}{(q_0-1)[d+r(d-2)]}\right)^{1+\frac{rd}{2r-d}}(R -R_\infty)^2|B_{R } |^{\frac{2}{2^*}}\,,\quad &\mbox{ if }q_0>1\,,\\
\left(\frac{\mathcal{S}_2^2}{2}\right)^{\frac{rd}{2r-d}}\frac{2r-d}{rd}
        \left(\frac{(q_0+1) rd}{q_0[d+r(d-2)]}\right)^{1+\frac{rd}{2r-d}}(R -R_\infty)^2|B_{R } |^{\frac{2}{2^*}}\,,\quad &\mbox{ if }0< q_0 \le 1\,.\\
\end{array}
\right.
\end{equation}
\end{thm}
\noindent {\sl Proof.~}The upper bounds \eqref{Moser.Upper.b.DG} of Theorem \ref{Moser.Upper.b.thm} can be rewritten as
\begin{equation}\label{DG1}
\|u\|_{\infty,r}\le \frac{K^{(3)}_q[b]}{(R-r)^{\frac{d}{q}}}
    \;\|u\|_{q,R}
\end{equation}
for any $q>1$ and $R_\infty\le r<R\le R_0$, where $K^{(3)}_q[b]$ is given by \eqref{Const.Upper.b}. It is clear that inequality \eqref{DG1} guarantees that we can use Lemma \ref{Lem.DG.Ext} with $0<\underline{q}=q<+\infty=\overline{q}$\,, $\gamma=d/q>1$\,, $K=K^{(3)}_q[b]$ and for any $R_\infty\le r<R\le R_0$. Then we have that for all $0<q_0\le \overline{q}=q$
\[
\|u\|_{\infty, R_\infty}
    \le3\cdot 2^{\frac{q-q_0}{q_0}}\,\left[\left(4\frac{d}{q}\frac{q}{q_0}\right)^{\frac{d}{q}}
        \frac{K^{(3)}_q[b]}{(R_0-R_\infty)^{\frac{d}{q}}}\right]^{\frac{q}{q_0}}\,\|u\|_{q_0,R_0}
    =3\cdot 2^{\frac{2d+1}{q_0}}\,\left(\frac{d}{q_0}\right)^{\frac{d}{q_0}}
        \frac{K^{(3)}_{q_0+1}[b]^{\frac{q_0+1}{q_0}}}{(R_0-R_\infty)^{\frac{d}{q_0}}}\,\|u\|_{q_0,R_0}
\]
since we can always choose $q=q_0+1>1$. Finally we notice that we can rewrite the upper bound for all $q_0>0$ in the following form:
\[
\|u\|_{\infty,r}\le
    \frac{A_{q_0}^{(1)}}{(R_0-R_\infty)^{\frac{d}{q_0}}}
    \left[A_{q_0}^{(2)}+A_{q_0}^{(3)}\|b\|_{\LL^r({B_{R_0} })}^{\frac{rd}{2r-d}}\right]^{\frac{d}{2q_0}}\;\|u\|_{q_0,R}
\]
where $A_q^{(j)}$ are as in \eqref{A1.b}, \eqref{A2.b} and \eqref{A3.b} respectively.\qed

\medskip

\noindent The above Theorem has the following important consequence, when applied to the equation $-\Delta u=\lambda u^p$.
\begin{thm}[Local Upper bounds, second form]\label{Thm.local.upper.b.p}
Consider a weak subsolution $u$ to $-\Delta u= \lambda u^p$ on $B_R$, with $\lambda>0$, $1< p<p_s=2^*-1=(d+2)/(d-2)$. Let $0<R_\infty<R_0<R$. If $u\in\LL^{\overline{r}}(B_{R_0})$ with $\overline{r}>d(p-1)/2:=\overline{q}$ then the following bound holds true for any $q_0>0$
\begin{equation}\label{Moser.Upper.b.p}
\|u\|_{\infty,r}\le
    \frac{A_{q_0}^{(1)}}{(R_0-R_\infty)^{\frac{d}{q_0}}}
    \left[A_{q_0}^{(2)}+A_{q_0}^{(3)}\lambda^{\frac{d(p-1)}{2\overline{r}-d(p-1)}}\|u\|_{\overline{r},R_0}^{\frac{d(p-1)\overline{r}}{2\overline{r}-d(p-1)}}\right]^{\frac{d}{2q_0}}\;\|u\|_{q_0,R}
\end{equation}
where $A_q^{(j)}$ are as in \eqref{A1.b}, \eqref{A2.b} and \eqref{A3.b} respectively.
\end{thm}
\noindent {\sl Proof.~}Since $u$ is a subsolution to $-\Delta u= \lambda u^p = bu$ with $b=\lambda u^{p-1}$, we need to assume that $u^{p-1}\in \LL^r$ with $r>d/2$, which amounts to require $u\in \LL^{\overline{r}}$ with $\overline{r}=r(p-1)>d(p-1)/2$,  so that
\[
\|b\|_{\LL^r({B_{R_0} })}^{\frac{rd}{2r-d}}=\left(\lambda \int_{B_{R_0}}u^{r(p-1)} \dx\right)^{\frac{d}{2r-d}}
=\lambda^{\frac{d(p-1)}{2\overline{r}-d(p-1)}}\|u\|_{\overline{r},R_0}^{\frac{d(p-1)\overline{r}}{2\overline{r}-d(p-1)}}
\]
Finally, we can apply the bounds of Theorem \ref{Thm.local.upper.b.DG} to get the bounds \eqref{Moser.Upper.b.p} with the constants written above.\qed

\normalcolor


\section{Lower bounds}\label{lb}

The lower bounds can be obtained in two steps: first we perform a Moser iteration, then we need reverse H\"older inequalities, which are a consequence of the celebrated John-Nirenberg Lemma.

\subsection{A short reminder about the spaces $M^p(\Omega)$. }We recall here some basic definitions and properties of suitable functional spaces, that will be used in the sequel. We omit the proofs, but we give appropriate references.

We say that a measurable function on $\Omega\subseteq \RR^d$ belong to the space $M^p(\Omega)$ if and only if there exists a constant $K\ge 0$ such that
\[
\int_{\Omega\cap B_R(x_0)}|f|\dx\le K R^{\frac{d(p-1)}{p}}\,\qquad\mbox{for all }B_R(x_0),
\]
and we define the norm on $M^p(\Omega)$ as follows
\[
\|f\|_{M^p(\Omega)}=\inf\left\{K>0\;:\;\int_{\Omega\cap B_R(x_0)}|f|\dx\le K R^{\frac{d(p-1)}{p}}\,\quad\mbox{for all }B_R(x_0) \right\}.
\]
One can easily check the strict inclusion $\LL^p(\Omega)\subset M^p(\Omega)$ for all $1<p<\infty$, and when $\Omega$ is bounded, the equalities $\LL^1(\Omega)=M^1(\Omega)$ and $\LL^\infty(\Omega)=M^\infty(\Omega)$. Moreover it is easy to check that when $\Omega$ is bounded one has:
\begin{equation}\label{Mp.L1}
\|f\|_{\LL^1(\Omega)}\le \diam(\Omega)^{\frac{d(p-1)}{p}}\|f\|_{M^p(\Omega)}.
\end{equation}
We now proceed with a series of results that relate the Marcinckievitz norm with the Riesz potential
\begin{equation}\label{Riesz.Potential}
\mathcal{V}_\mu [f](x):=\int_\Omega\frac{f(y)}{|x-y|^{d(1-\mu)}}\dy\,\qquad\mbox{with }\mu\in(0,1].
\end{equation}
We collect hereafter some well known results, whose proof can be found for instance in \cite{GT}.
\begin{lem}\label{GT.Lemmata}
Let $\mathcal{V}_\mu$ be defined as above. Then the following holds.

\noindent(i) The operator $\mathcal{V}_\mu$ maps continuously $\LL^s(\Omega)$ into $\LL^r(\Omega)$ for any $1\le r\le \infty$ satisfying
\[
0\le \frac{1}{s}-\frac{1}{r}<\mu.
\]
Moreover, for any $f\in\LL^p(\Omega)$,
\[
\|\mathcal{V}_\mu f\|_r\le\left(\frac{s(r+1)-r}{s(\mu\,r+1)-r}\right)^{\frac{s(r+1)-r}{sr}}\omega_d^{1-\mu}|\Omega|^{\frac{s(\mu\,r+1)-r}{sr}}\|f\|_s.
\]
\noindent(ii) Let $f\in M^p(\Omega)$, with $p>1/\mu\ge 1$. Then
\[
\left|\mathcal{V}_\mu[f](x)\right|\le \frac{p-1}{p\mu-1}\diam(\Omega)^{\frac{d}{p}(p\mu-1)}\|f\|_{M^p(\Omega)}.
\]
\noindent(iii) {\em A ``potential'' version of the Morrey inequality. }Let $\Omega$ be a convex bounded subset of $\RR^d$. Then for all $f\in W^{1,1}(\Omega)$ the following inequality holds
\begin{equation}\label{Morrey.potential}
|f(x)-f_{\Omega'}|\le \frac{\diam(\Omega)^d}{d\,|\Omega'|}\,\left|\mathcal{V}_{\frac{1}{d}}[\,|\nabla f|](x)\right|
\end{equation}
for any measurable $\Omega'\subseteq\Omega$ with
\[
f_{\Omega'}=\int_{\Omega'}f\frac{\dx}{|\Omega'|}
\]
\end{lem}
\noindent {\sl Proof.~}Part (i) is exactly Lemma 7.12 of \cite{GT}, part (ii) is exactly Lemma 7.18 of \cite{GT} and  part (iii) is exactly Lemma 7.16 of \cite{GT}.\qed

\subsection{The John-Nirenberg Lemma and reverse H\"older inequalities. }The Caccioppoli estimates proved in Corollary \ref{Lem.energy2} show that the gradient of the logarithm of the solution belongs to the Marcinckievitz space $M^d(\Omega)$, see Proposition \ref{JN.cor} below. Such  $M^d-$regularity then guarantees the validity of the celebrated John-Nirenberg lemma which as a consequence give a reverse H\"older inequality of the form
\[
\frac{\|u\|_{q,R_{0}}}{\|u\|_{-q,R_0}}\le \kappa_1^{2/q}
\]
for some $0<q<1$ and some constant $\kappa_1$.

We need a lemma concerning estimates on the Riesz potential $\mathcal{V}_\mu$ defined in \eqref{Riesz.Potential}. It is a quantified version of Lemma 7.20 of \cite{GT}.

\begin{lem}[A ``potential'' version of the Moser-Trudinger imbedding.]\label{Lem.Pot.Mos-Tru}
Let $f\in M^p(\Omega)$ with $p>1$ and suppose $\|f\|_{M^p(\Omega)}\le K$. Then there exist two constants $\kappa_2$ and $\kappa_3$ such that
\begin{equation}\label{Moser.Trudinger.potential}
\int_\Omega\exp\left[\frac{\left|\mathcal{V}_{\frac{1}{p}}[f](x)\right|}{\kappa_2\, K}\right]\,\dx\le \kappa_3.
\end{equation}
One can take
\[
\kappa_2>(p-1)\,\ee \qquad\mbox{and}\qquad
\kappa_3=|\Omega|+\frac{\diam(\Omega)^d}{\sqrt{2\pi}}\frac{p\,\ee\,\omega_d}{\kappa_2-(p-1)\,\ee}.
\]
\end{lem}
\noindent {\sl Proof.~} Let $q\ge 1$, $\mu=1/p$ and $g=\mathcal{V}_{\mu}[f]$. Then
\[
|x-y|^{d(\mu-1)}=|x-y|^{\frac{d}{q}\left(\frac{\mu}{q}-1\right)}|x-y|^{d\left(1-\frac{1}{q}\right)\left(\frac{\mu}{q}+\mu-1\right)}
\]
and by H\"older inequality we obtain
\begin{equation}\label{GT.00}
|\mathcal{V}_{\mu}[f]|\le \left|\mathcal{V}_{\frac{\mu}{q}}[f]\right|^{\frac{1}{q}}
    \left|\mathcal{V}_{\mu+\frac{\mu}{q}}[f]\right|^{1-\frac{1}{q}}.
\end{equation}
Applying now estimates (i) of Lemma \ref{GT.Lemmata} with $s=r=1$, to $\mathcal{V}_{\frac{\mu}{q}}[f]$, we obtain,
\begin{equation}\label{GT.01}\begin{split}
\|\mathcal{V}_{\frac{\mu}{q}} f\|_1
&\le\frac{q\,\omega_d^{1-\frac{\mu}{q}}}{\mu}|\Omega|^{\frac{\mu}{q}}\|f\|_1
 \le p\,q\,\omega_d^{1-\frac{1}{pq}}|\Omega|^{\frac{1}{pq}}\diam(\Omega)^{\frac{d(p-1)}{p}}\|f\|_{M^p(\Omega)}\\
&\le p\,q\,\omega_d\diam(\Omega)^{d\left(1-\frac{1}{p}+\frac{1}{pq}\right)}\|f\|_{M^p(\Omega)}
\le p\,q\,\omega_d\diam(\Omega)^{d\left(1-\frac{1}{p}+\frac{1}{pq}\right)}K
\end{split}
\end{equation}
where we have used inequality \eqref{Mp.L1} together with the fact that $|\Omega|\le \omega_d \diam(\Omega)^d$. Next we apply estimates (ii) of Lemma \ref{GT.Lemmata} to $\mathcal{V}_{\mu+\frac{\mu}{q}}[f]$ (the operator $\mathcal{V}_{\nu}$ is well-defined on L$^1$, if $\Omega$ is bounded, for $\nu>1$ as well) and we obtain
\begin{equation}\label{GT.02}
\left|\mathcal{V}_{\mu+\frac{\mu}{q}}[f](x)\right|\le \frac{p-1}{p\left(\mu+\frac{\mu}{q}\right)-1}\diam(\Omega)^{\frac{d}{p}\left(p\left(\mu+\frac{\mu}{q}\right)-1\right)}\|f\|_{M^p(\Omega)}
\le q(p-1)\,\diam(\Omega)^{\frac{d}{pq}}K
\end{equation}
for all $x\in\Omega$, hence the same bound is valid for the $\LL^\infty(\Omega)$-norm, provided $p(\mu+\mu/q)>1$ which indeed holds true since $\mu=1/p$. Joining now inequalities \eqref{GT.00}, \eqref{GT.01} and \eqref{GT.02}, we obtain
\[\begin{split}
\|\mathcal{V}_{\mu}[f]\|_q^q
     &\le  \left\|\mathcal{V}_{\mu+\frac{\mu}{q}}[f]\right\|_{\LL^\infty}^{q-1}
            \left\|\mathcal{V}_{\frac{\mu}{q}}[f](x)\right\|_{\LL^1(\Omega)}
      \le \frac{p\,\omega_d}{p-1}\left[(p-1)\,K \, q\right]^q\diam(\Omega)^d
\end{split}
\]
Letting now $1\le q=k\in\mathbb{N}$ we get, for $k_2$ as in the statement,
\[\begin{split}
\int_\Omega \sum_{k=1}^{\infty}\frac{|g|^k}{k!(\kappa_2 K)^k}\dx
    &\le \frac{p\,\omega_d}{p-1}\diam(\Omega)^d
        \sum_{k=1}^{\infty} \frac{\left[(p-1)\,K \, k\right]^k}{k!(\kappa_2 K)^k}
     \le \frac{p\,\omega_d}{p-1}\diam(\Omega)^d
        \sum_{k=1}^{\infty} \left[\frac{p-1}{\kappa_2}\right]^k\frac{k^k}{k!}\\
    &\le \frac{p\,\omega_d}{p-1}\diam(\Omega)^d
        \sum_{k=1}^{\infty} \left[\frac{(p-1)\,\ee}{\kappa_2}\right]^k\frac{1}{\sqrt{2\pi k}}\\
    &\le \frac{p\,\omega_d}{p-1}\frac{\diam(\Omega)^d}{\sqrt{2\pi}}\frac{(p-1)\,\ee}{\kappa_2-(p-1)\,\ee}
        =\frac{\diam(\Omega)^d}{\sqrt{2\pi}}\frac{p\,\ee\,\omega_d}{\kappa_2-(p-1)\,\ee}\\
\end{split}
\]
we have used Stirling's formula:
\begin{equation}
n!=\sqrt{2\pi\, n}\,\left[\frac{n}{e}\right]^n\,\ee^{\alpha_n}\qquad\mbox{with}\qquad \frac{1}{12n+1}\le \alpha_n \le\frac{1}{12n}\,\mbox{.\qed}
\end{equation}

We prove hereafter a simplified but quantitative version of the celebrated John-Nirenberg Lemma, which holds in convex domains. Indeed we will use it only on balls and in such case the constants simplify a bit.
\begin{lem}[John-Nirenberg]\label{John.Nirenberg}
Let $f\in W^{1,1}(\Omega)$ where $\Omega$ is convex, and suppose there exists a constant $K$ such that
\[
\int_{B_R\cap\Omega} \big|\nabla f\big|\dx \le K\, R^{d-1}\qquad\mbox{for all balls }B_R
\]
Then the following inequality holds true
\begin{equation}\label{JN.ineq}
\int_\Omega\exp\left[ \frac{|f-f_\Omega|}{\kappa_0 K}\right] \dx\le \kappa_1
\end{equation}
where for any $\kappa_2>(d-1)\,\ee$
\[
\kappa_0=\frac{d\,|\Omega|}{\diam(\Omega)^d}\kappa_2\, \qquad
\kappa_1=\frac{\omega_d\,\diam(\Omega)^d\big(\kappa_2+\ee\big)}{\kappa_2-(d-1)\,\ee}\qquad\mbox{and}\qquad
f_\Omega=\int_\Omega f\frac{\dx}{|\Omega|}.
\]
\end{lem}
\noindent {\sl Proof.~}The proof relies on the previous Lemma \ref{Lem.Pot.Mos-Tru} in the special case $p=d$. Indeed inequality \eqref{Moser.Trudinger.potential} in that case takes the form
\begin{equation}\label{Moser.Trudinger.potential.d}
\int_\Omega\exp\left[\frac{\left|\mathcal{V}_{\frac{1}{d}}[|\nabla f|](x)\right|}{\kappa_2 K}\right]\,\dx\le \frac{\diam(\Omega)^d}{\sqrt{2\pi}}\frac{d\,\ee\,\omega_d}{\kappa_2-(d-1)\,\ee}+|\Omega|\le \kappa_3
\end{equation}
where
\[
\kappa_2>(d-1)\,\ee \qquad\mbox{and}\qquad
\kappa_3=\omega_d\,\diam(\Omega)^d \left[\frac{d\,\ee}{\kappa_2-(d-1)\,\ee}+1\right]
= \frac{\omega_d\,\diam(\Omega)^d\big(\kappa_2+\ee\big)}{\kappa_2-(d-1)\,\ee}.
\]
We combine this latter inequality with inequality \eqref{Morrey.potential} (which requires convexity of the domain) with $\Omega'=\Omega$ and $|\nabla f|\in M^d(\Omega)$.\qed

\medskip

\noindent The John-Nirenberg Lemma has an important consequence when applied to $f=\log(u+\delta)$:
\begin{prop}[Reverse H\"older inequalities]
Let $\delta\ge 0$ and $u$ be a positive measurable function such that $\log (u+\delta)\in W^{1,1}(\Omega)$, where $\Omega$ is convex, and suppose there exists a constant $K$ such that
\begin{equation}\label{log}
\int_{B_R\cap\Omega} \big|\nabla \log (u+\delta)\big|\dx \le K\, R^{d-1}\qquad\mbox{for all balls }B_R.
\end{equation}
Then the following inequality
\begin{equation}\label{rev.holder.sols}
\frac{\|u+\delta\|_{q,\Omega}}{\|u+\delta\|_{-q,\Omega}}\le \kappa_1^{2/q}\qquad\mbox{holds true for any
}\qquad 0<q\le \frac{1}{\kappa_0\,K}
\end{equation}
where the constants $\kappa_i$ are given in Lemma \ref{John.Nirenberg}.
\end{prop}
\noindent {\sl Proof.~} Let $\delta >0$. The validity of \eqref{log} for $u$ entails the validity of the same inequality for $u+\delta$. Notice now that
\[
\frac{\|u+\delta\|_{q,\Omega}}{\|u+\delta\|_{-q,\Omega}}\le \kappa\qquad\iff\qquad
\left(\int_{\Omega} (u+\delta)^{q}\dx\right)\left(\int_{\Omega} (u+\delta)^{-q}\dx\right)\le \kappa^{q}
\]
Then, letting $f=\log(u+\delta)$:
\[\begin{split}
&\left(\int_{\Omega}(u+\delta)^{q}\dx\right)\left(\int_{\Omega} (u+\delta)^{-q}\dx\right)
     =\left(\int_{\Omega}\ee^{\left[q\log(u+\delta)\right]}\dx\right)
        \left(\int_{\Omega}\ee^{\left[-q\log(u+\delta)\right]}\dx\right)\\
    &=\left(\int_{\Omega}\ee^{qf}\dx\right)\left(\int_{\Omega}\ee^{-qf}\dx\right)
    =\left(\int_{\Omega}\ee^{q\left(f-f_\Omega\right)}\dx\right)
         \left(\int_{\Omega}\ee^{-q\left(f-f_\Omega\right)}\dx\right)\\
    &\le \left(\int_{\Omega}\ee^{q\left|f-f_\Omega\right|}\dx\right)^2 \le \kappa_1^2
\end{split}\]
where we used \eqref{JN.ineq} for $f=\log(u+\delta)$, and have assumed $q\le 1/(\kappa_0 K)$ in order to ensure its validity. The case $\delta=0$ is also true, just by taking the limit $\delta\to 0$.\qed

We conclude this section by showing that reverse H\"older inequalities holds for local supersolutions to our problem, as a consequence of Caccioppoli estimates.

\begin{prop}[Reverse H\"older inequalities for supersolutions]\label{JN.cor}
Let $\Omega\subset\RR^d$ and let $\lambda>0$. Let $u$ be a local weak supersolution to $-\Delta u = \lambda u^p$, with $1\le p<p_s=2^*-1=(d+2)/(d-2)$. Then for any $\varepsilon>0$ the following inequality holds true for any $\delta\ge 0$
\[
\left[\frac{\varepsilon}{2^d \,(\ee\,d+\varepsilon)\, }\right]^{2/q}
\frac{\|u+\delta\|_{q,R_{0}}}{|B_{R_0}|^{\frac{1}{q}}}\le
\frac{\|u+\delta\|_{-q,R_0}}{|B_{R_0}|^{-\frac{1}{q}}}\qquad\mbox{for all}\qquad
0<q\le\frac{{2^{{\frac{d-3}{2} }} }}{d\omega_d^2[\ee(d-1)+\varepsilon]}.
\]
\end{prop}
\noindent {\sl Proof.~} The Caccioppoli estimates \eqref{Cacciopoli} with $R_0$ replaced by $2r$ and $R$ replaced by $r$ imply the hypothesis of the above Lemma, in fact:
\begin{equation}\label{RevHold.hyp}\begin{split}
\int_{B_{r}\cap B_{R_0}} \big|\nabla \log (u+\delta)\big|\dx
&\le \int_{B_{r}} \big|\nabla \log (u+\delta)\big|\dx\\
&\le |B_{r}|^{\frac{1}{2}}\left[\int_{B_{r}} \big|\nabla \log (u+\delta)\big|^2\dx\right]^\frac{1}{2}
\le 2^{{\frac{d+3}{2} }}\omega_d {r}^{d-1}:=K\,{r}^{d-1}.
\end{split}\end{equation}
Therefore putting $K=2^{{\frac{d+3}{2} }}\omega_d$, taking an $\varepsilon>0$ and choosing $\kappa_2=\ee(d-1)+\varepsilon$, we obtain that
\[
\frac{1}{\kappa_0 K}=\frac{2^{{\frac{d-3}{2} }}}{d\omega_d^2\,[\ee(d-1)+\varepsilon]}\qquad\mbox{and}\qquad
\kappa_1={2^d}\omega_dR_0^d\frac{\varepsilon+\ee\,d}{\varepsilon}=|B_{R_0}|2^d\frac{\varepsilon+\ee\,d}{\varepsilon}\,\mbox{.\qed}
\]

\subsection{Lower Moser iteration}

Now we are ready to run the Moser iteration to obtain quantitative local lower bounds in the form:

\begin{thm}[Local Lower Estimates]\label{thm.local.lower}
Let $\Omega\subseteq\RR^d$ and let $\lambda>0$. Let $u$ be a nonnegative local weak supersolution in $B_{R_0}\subseteq\Omega$ to $-\Delta u = \lambda u^p$, with $0\le p<p_s=2^*-1=(d+2)/(d-2)$. Then for any $\varepsilon>0$ and for any
\begin{equation}\label{q0-figure}
0<\underline{q}\le\frac{2^{{\frac{d-3}{2}}}}{d\omega_d^2[\ee(d-1)+\varepsilon]}=q_0
\end{equation}
the following bound holds true
\begin{equation}\label{lower.local.thm}
\inf_{x\in B_{R_\infty}}u(x) =\|u\|_{-\infty,R_\infty}
    \ge I_{-\infty,\underline{q}}\frac{\|u\|_{{\underline{{q}}},R_{0}}}{|B_{R_0}|^{\frac{1}{\underline{q}}}}.
\end{equation}
where
\begin{equation}\label{Const.lower.local.thm}
I_{-\infty,\underline{q}}
=\left[2^d\mathcal{S}_2^2\left(\frac{dR_0^2}{(R_0-R_\infty)^2}+\frac{R_0^2}{R_\infty^2}\right)
    \right]^{-\frac{d}{2\underline{q}}}
    \left[\frac{\varepsilon}{2^d \,(\ee\,d+\varepsilon)\,\sqrt{\omega_d} }\right]^{\frac{2}{\underline{q}}}.
\end{equation}
\end{thm}
\noindent\textbf{Remark. } One can see that when the dimension $d$ is sufficiently low one has $q_0<1$ whereas $q_0>1$ in higher dimensions. Notice also that the equality $\inf\limits_{x\in B_{R_\infty}}u(x) =\|u\|_{-\infty,R_\infty}$ holds since $u$ is nonnegative.

\begin{figure}[ht]
\centering
\includegraphics[height=5.10cm, width=7cm]{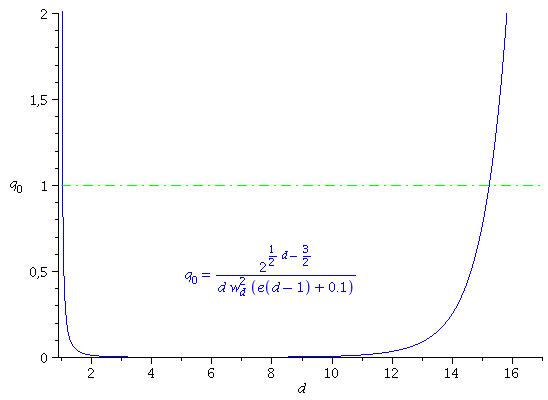}
\includegraphics[height=4.81cm, width=7cm]{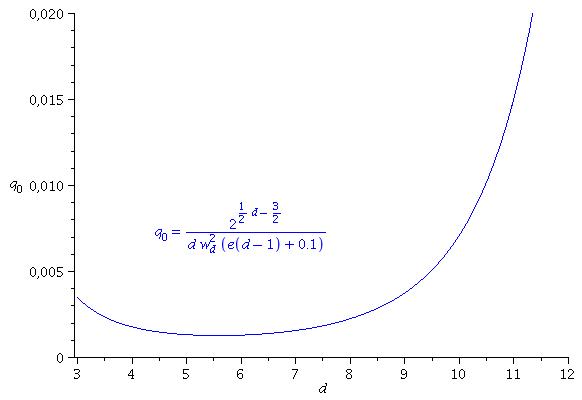}
\flushleft \noindent\textit{Left figure: Plot of $q_0(d)$ defined in \eqref{q0-figure}, for $1\le d\le 16$, with $\varepsilon=0.1$.\\
Right figure: Zoom for the plot of the same $q_0(x)$ near its minimum that lies in $(5\,,\,6)$.}
\end{figure}

\noindent {\sl Proof.~} The proof is divided in two steps. We always consider a local supersolution $u$ of $-\Delta u \ge \lambda u^p$.

\noindent$\bullet~$\textsc{Step 1.} In this step we consider $\alpha<0$, and we want to prove $\LL^{-q}-\LL^{-\infty}$ local estimates via Moser iteration. The the energy inequality \eqref{local.energy.identity.supersols} for $\alpha< -1$ and $\delta>0$ gives the estimate
\begin{equation}\label{local.energy.ineq.negative}\begin{split}
\int_\Omega \big|\nabla \left((u+\delta)^{\frac{\alpha+1}{2}}\right)\big|^2\varphi\dx
    &\le\frac{\lambda(\alpha+1)^2}{4\alpha}\int_\Omega u^p(u+\delta)^{\alpha}\varphi\dx
        +\frac{\alpha+1}{4\alpha}\int_\Omega (u+\delta)^{\alpha+1}\Delta\varphi \dx\\
    &\le \frac{\alpha+1}{4\alpha}\int_\Omega (u+\delta)^{\alpha+1}\big|\Delta\varphi\big| \dx
    \end{split}
\end{equation}
Applying now the Sobolev inequality \eqref{Sobolev.BR} on the ball $B_{R_1}$ and the properties of the test function $\varphi$ defined in Lemma \ref{lem.test.funct}, one gets
\begin{equation}\label{alpha.neg}\begin{split}
\left[\int_{B_{R_1}} (u+\delta)^{\frac{2^*}{2}(\alpha+1)}\dx\right]^{\frac{2}{2^*}}
    &\le \mathcal{S}_2^2 \left(\int_{B_{R_1}} \big|\nabla (u+\delta)^{\frac{\alpha+1}{2}}\big|^2\dx
         +\frac{1}{R_1^2}\int_{B_{R_1}} (u+\delta)^{\alpha+1}\dx\right)\\
    &\le \mathcal{S}_2^2 \left(\int_{\Omega} \big|\nabla (u+\delta)^{\frac{\alpha+1}{2}}\big|^2\varphi\dx
         +\frac{1}{R_1^2}\int_{B_{R_1}} (u+\delta)^{\alpha+1}\dx\right)\\
    &\le \mathcal{S}_2^2 \left(\frac{\alpha+1}{4\alpha}\int_\Omega (u+\delta)^{\alpha+1}\big|\Delta\varphi\big| \dx
         +\frac{1}{R_1^2}\int_{B_{R_1}} (u+\delta)^{\alpha+1}\dx\right)\\
    &\le \mathcal{S}_2^2 \left(\frac{\alpha+1}{4\alpha}\big\|\Delta\varphi\big\|_\infty
         +\frac{1}{R_1^2}\right)\int_{B_{R_0}} (u+\delta)^{\alpha+1}\dx\\
    &\le \mathcal{S}_2^2 \left(
         \frac{d}{ (R_0-R_1)^2}
         +\frac{1}{R_1^2}\right)\int_{B_{R_0}} (u+\delta)^{\alpha+1}\dx\\
\end{split}
\end{equation}
Let, for a given $\gamma_0<0$, $\gamma_n:=\left[\frac{2^*}{2}\right]^n\gamma_0$ so that $\gamma_n=\frac{2^*}{2}\gamma_{n-1}$. Notice that $\gamma_n\to-\infty$ monotonically. Consider the above inequality for $\alpha=\alpha_n$ and let  $\alpha_n+1=\gamma_{n-1}$ so that
\begin{equation}\begin{split}
\|u+\delta\|_{\gamma_n,R_n}
    &=\|u+\delta\|_{\frac{2^*}{2}\gamma_{n-1},R_n}
     =\left[\int_{B_{R_n}} (u+\delta)^{\frac{2^*}{2}\gamma_{n-1}}\dx\right]^{\frac{2}{2^*\gamma_{n-1}}}\\
    &\ge\left[\mathcal{S}_2^2 \left(
         \frac{d}{ (R_{n-1}-R_n)^2}
         +\frac{1}{R_n^2}\right)\right]^{\frac{1}{\gamma_{n-1}}}\left[\int_{B_{R_{n{-1}}}} (u+\delta)^{\gamma_{n-1}}\dx\right]^{\frac{1}{\gamma_{n-1}}}\\
    &\ge\left[
         \mathcal{S}_2^2 \left(
         \frac{d}{ (R_{n-1}-R_n)^2}
         +\frac{1}{R_n^2}\right)
         \right]^{\frac{1}{\gamma_{n-1}}}\|u+\delta\|_{\gamma_{n-1},R_{n-1}}
     :=I_n^{\frac{1}{\gamma_{n-1}}}\|u+\delta\|_{\gamma_{n-1},R_{n-1}}\\
\end{split}
\end{equation}
Hence, iterating the above inequality:
\begin{equation}\label{iter.neg}
\|u+\delta\|_{\gamma_n,R_n}\ge I_n^{\frac{1}{\gamma_{n-1}}}I_{n-1}^{\frac{1}{\gamma_{n-2}}}\ldots I_1^{\frac{1}{\gamma_0}}\|u+\delta\|_{\gamma_0,R_0}=\prod_{k=1}^{n}I_k^{\frac{1}{\gamma_{k-1}}}\|u+\delta\|_{\gamma_0,R_0}
\end{equation}
where have chosen $0<R_\infty<\ldots<R_{n+1}<R_n<\ldots<R_0$ such that
\[
\sum_{k=1}^\infty ( R_{k-1}-R_k ) =R_0-R_\infty\qquad\mbox{and}\qquad R_{k-1}-R_k=\frac{R_0-R_\infty}{2^k}
\]
so that
\[
I_k=\mathcal{S}_2^2 \left(
         \frac{d}{ (R_{n-1}-R_n)^2}
         +\frac{1}{R_n^2}\right)\le
         \mathcal{S}_2^2 \left(
         \frac{d}{ (R_{0}-R_\infty)^2}
         +\frac{1}{R_{\infty}^2}\right)4^k:=I_0 4^k
\]
and
\[\begin{split}
\prod_{k=1}^{n}I_k^{\frac{1}{\gamma_{k-1}}}
&=\exp\left[\sum_{k=1}^{n} \frac{1}{\gamma_{k-1}}\log I_k\right]
 =\exp\left[\frac{2^*}{2\gamma_0}\sum_{k=1}^{n}\left[\frac{2}{2^*}\right]^k\log I_k\right]\\
 &=\exp\left[\frac{2^*}{2\gamma_0}\sum_{k=1}^{n}\left[\frac{2}{2^*}\right]^k\log I_0 +\frac{2^*\log 4}{2\gamma_0}\sum_{k=1}^{n}k\left[\frac{2}{2^*}\right]^k\right]\\
 &\ge I_0^{\frac{2^*}{2\gamma_0}\sum_{k=1}^{n}\left[\frac{2}{2^*}\right]^k}
 4^{\frac{d^2}{4\gamma_0}}.
\end{split}
\]
Taking limits we obtain
\[
\prod_{k=1}^{\infty}I_k^{\frac{1}{\gamma_{k-1}}}\ge I_0^{\frac{2^*}{2\gamma_0}\frac{d-2}{2}}
 4^{\frac{d^2}{4\gamma_0}}=\left(2^d\,I_0\right)^{\frac{d}{2\gamma_0}}.
\]
We can now take the limit in \eqref{iter.neg}  to get for any $\gamma_0<0$:
\begin{equation}\label{lower.bounds}\begin{split}
    \|u+\delta\|_{-\infty,R_\infty}
    &\ge \prod_{k=1}^{{\infty}}I_k^{\frac{1}{\gamma_{k-1}}}\|u+\delta\|_{\gamma_0,R_0}
    \ge \left(2^d\,I_0\right)^{\frac{d}{2\gamma_0}}\|u+\delta\|_{\gamma_0,R_0}\\
    &=\left[2^d\mathcal{S}_2^2 \left(
         \frac{d}{(R_{0}-R_\infty)^2}
         +1\right)\right]^{\frac{d}{2\gamma_0}}\|u+\delta\|_{\gamma_0,R_0}.
\end{split}\end{equation}
Now we need some Reverse H\"older inequalities, which is the subject of the next step.

\noindent$\bullet~$\textsc{Step 2. }\textit{Reverse H\"older inequalities.} The John-Nirenberg lemma implies reverse H\"older inequalities for super-solutions, in the form of Proposition \ref{JN.cor}: for any $\varepsilon>0$ the following inequality holds true
\begin{equation}\label{rev.holder.JN}
\left[\frac{\varepsilon}{2^d \,(\ee\,d+\varepsilon)}\right]^{\frac{2}{\underline{q}}}
\frac{\|u+\delta\|_{\underline{q},R_{0}}}{|B_{R_0}|^{\frac{2}{\underline{q}}}}\le
\|u+\delta\|_{-\underline{q},R_0}
\qquad\mbox{for all}\qquad
0<\underline{q}\le\frac{2^{\frac{d-3}{2}}}{d\omega_d^2[\ee(d-1)+\varepsilon]}.
\end{equation}
Joining inequality \eqref{lower.bounds} and \eqref{rev.holder.JN} and letting $\gamma_0=-\underline{q}$ with $\underline{q}$ as in \eqref{rev.holder.JN} we obtain
\begin{equation}\label{lower.bounds.1}\begin{split}
\|u+\delta\|_{-\infty,R_\infty}
    &\ge \left[2^d\mathcal{S}_2^2 \left(
         \frac{d}{(R_{0}-R_\infty)^2}
         +\frac{1}{R_\infty^2}\right)\right]^{-\frac{d}{2\underline{q}}}\|u+\delta\|_{-\underline{q},R_0}\\
    &\ge\left[2^d\mathcal{S}_2^2 \left(
         \frac{d}{(R_{0}-R_\infty)^2}
         +\frac{1}{R_\infty^2}\right)\right]^{-\frac{d}{2\underline{q}}}
    \left[\frac{\varepsilon}{2^d \,(\ee\,d+\varepsilon)}\right]^{\frac{2}{\underline{q}}}
        \frac{\|u+\delta\|_{\underline{q},R_{0}}}{|B_{R_0}|^{\frac{2}{\underline{q}}}}\\
    &= \left[2^d\mathcal{S}_2^2\left(\frac{dR_0^2}{(R_0-R_\infty)^2}+\frac{R_0^2}{R_\infty^2}\right)
    \right]^{-\frac{d}{2\underline{q}}}
    \left[\frac{\varepsilon}{2^d \,(\ee\,d+\varepsilon)\,\sqrt{\omega_d} }\right]^{\frac{2}{\underline{q}}}
        \frac{\|u+\delta\|_{\underline{q},R_{0}}}{|B_{R_0}|^{\frac{1}{\underline{q}}}}\\
    &:= I_{-\infty,\underline{q}}\frac{\|u+\delta\|_{\underline{q},R_{0}}}{|B_{R_0}|^{\frac{1}{\underline{q}}}}.
\end{split}
\end{equation}
Finally we observe that we can let $\delta\to 0^+$, and obtain the desired result.\qed

\subsection{Reverse H\"older inequalities and additional local lower bounds when $1<p<p_c$}\label{Sect.Finite.Moser}
In this section we will first prove more quantitative \textit{reverse H\"older inequalities}, when $p>1$. We have obtained a reverse smoothing effect from $\LL^q$ to $\LL^{-\infty}$, for a suitable explicit $q$ which may be close to zero, if we seek for a bound valid for any dimension. In order to be able to join local upper and lower estimates to get a clean form of Harnack inequality, we need to reach those values of $q$ which are above $d(p-1)/2$, and this is possible only when $1<p<p_c=d/(d-2)$.

\begin{prop}[Reverse H\"older inequalities for $1<p<p_c$]\label{rev.holder.pc}
Let $\Omega\subseteq\RR^d$ and let $\lambda>0$. Let $u$ be a  nonnegative local weak supersolution in $\Omega$ to $-\Delta u = \lambda u^p$, with $1< p<p_c=d/(d-2)$. Let $B_{\overline{R}}\subset B_{R_0}\subset\Omega$.  Then we have that
\begin{equation}\label{iter.neg.2}
\frac{\|u\|_{\overline{q},\overline{R}}}{|B_{\overline{R}}|^{\frac{1}{\overline{q}}}}\le I_{\overline{q},q_0}\frac{\,\|u\|_{q_0,R_0}}{|B_{R_0}|^{\frac{1}{q_0}}}\qquad\mbox{for any }\quad q_0\in(0,\overline{q}]\normalcolor \mbox{ and } d(p-1)/2<\overline{q}<d/(d-2)
\end{equation}
where
\begin{equation}\label{I.pc}
I_{\overline{q},q_0}:=\left\{\begin{array}{lll}
\left[\frac{2d\,\overline{q}\,\mathcal{S}_2^2}{(2^*-2\overline{q})}
    +\,\mathcal{S}_2^2\frac{(R_0-\overline{R})^2}{\overline{R}^2}\right]^\frac{2^*}{2\overline{q}}
    \left[\frac{\omega_d^{1/d}R_0}{R_0-\overline{R}}\right]^{\frac{2^*}{\overline{q}}}
    \left[\frac{R_0}{\overline{R}}\right]^{\frac{d}{\overline{q}}}
    & \mbox{ if }\frac{d-2}{d}\overline{q}\le q_0\le \overline{q} \,, \\
3\cdot 2^{\frac{(d-2)\overline{q}}{2q_0}-\frac{d}{2}}\,
        \,\left[\frac{2d\,\overline{q}\,\mathcal{S}_2^2}{(2^*-2\overline{q})}\frac{\overline{R}^2}{(R_0-\overline{R})^2}+\,\mathcal{S}_2^2
            \right]^{{\frac{\overline{q}-q_0}{\overline{q}\,q_0}\frac{d}{2}}}
        \,\left(4\omega_d^{\frac{1}{d}}\frac{\overline{q}-q_0}{q_0\overline{q}}\right)^{\frac{d}{q_0}-\frac{d}{\overline{q}}}
            \left[\frac{\overline{R}}{R_0}\right]^{\frac{d}{q_0}}\,,
    & \mbox{ if }0<q_0<\frac{d-2}{d}\overline{q}.\\
\end{array}\right.
\end{equation}
\end{prop}
\noindent {\sl Proof.~} Consider the energy identity for supersolutions with $-1<\alpha<0$ (we can take $\delta=0$ in such a range of $\alpha$), which gives the following estimate for any positive test function $\varphi\in C^2_0(\Omega)$ with $\nabla \varphi\equiv0$ on $\partial\Omega$:
\begin{equation}\begin{split}
\frac{4|\alpha|}{(\alpha+1)^2}\int_\Omega \big|\nabla u^{\frac{\alpha+1}{2}}\big|^2\varphi\dx+\lambda\int_\Omega u^{p+\alpha}\varphi\dx
    &\le \frac{1}{|\alpha+1|}\int_\Omega u^{\alpha+1}|\Delta\varphi|\dx\\
    \end{split}
\end{equation}
that implies, using the test function $\varphi$ of Lemma \ref{lem.test.funct} with $R_\infty<R_0$
\begin{equation}\label{parte1}
\int_{B_{R_\infty}} \big|\nabla u^{\frac{\alpha+1}{2}}\big|^2\dx\leq \frac{d|\alpha +1|}{|\alpha|\,(R_0-R_\infty)^2}\int_{B_{R_0}} u^{\alpha+1}\dx
\end{equation}
Applying now the Sobolev inequality \eqref{Sobolev.BR.b} on the ball $B_{R_\infty}$ we arrive at
$$
\left[\int_{B_{R_\infty}} u^{\frac{2^*}{2}(\alpha+1)}\dx\right]^{\frac{2}{2^*}}
\le \mathcal{S}_2^2\left[ \frac{d|\alpha +1|}{|\alpha|\,(R_0-R_\infty)^2}+\frac{1}{R_\infty^2}\right]\int_{B_{R_0}} u^{\alpha+1}\dx
$$
Letting now $0<\alpha+1=\beta<1$ we get
\begin{equation}\label{parte2}
\left[\int_{B_{R_\infty}} u^{\frac{2^*}{2}\beta}\dx\right]^{\frac{2}{2^*\beta}}
\le\left[\frac{\mathcal{S}_2}{R_0-R_\infty}\right]^\frac{2}{\beta}
\left[\frac{d\,|\beta|}{(1-\beta)}+\frac{(R_0-R_\infty)^2}{R_\infty^2}\right]^\frac{1}{\beta}\left[\int_{B_{R_0}} u^{\beta}\dx\right]^\frac{1}{\beta}\,.
\end{equation}
Choosing $\beta>(d-2)(p-1)/2$ is compatible with $\beta<1$, if and only if $p<d/(d-2)=p_c$ and this is the point where the well known Serrin's exponent $p_c$ enters. We now let $d(p-1)/2<\overline{q}=2^*\beta/2<2^*/2$ and we see that \eqref{parte2} implies
\begin{equation}\label{DG55}
\|u\|_{\overline{q},r}\le
\left[\frac{2d\,\overline{q}\,\mathcal{S}_2^2}{(2^*-2\overline{q})}+\mathcal{S}_2^2\frac{(R-r)^2}{r^2}\right]^\frac{2^*}{2\overline{q}}
\frac{\|u\|_{\frac{2}{2^*}\overline{q},R}}{(R-r)^\frac{2^*}{\overline{q}}}
\le\left[\frac{2d\,\overline{q}\,\mathcal{S}_2^2}{(2^*-2\overline{q})}+\mathcal{S}_2^2\frac{(R_0-R_\infty)^2}{R_\infty^2}\right]^\frac{2^*}{2\overline{q}}
\frac{\|u\|_{\frac{2}{2^*}\overline{q},R}}{(R-r)^\frac{2^*}{\overline{q}}}
\end{equation}
for any $R_\infty\le r< R\le R_0$. Let $\underline{q}=2\overline{q}/2^*<\overline{q}$. We consider separately the case $\underline{q}\le q_0\le \overline{q}$ and the case $0<q_0<\underline{q}<\overline{q}$. In the first case we can use H\"older inequality in \eqref{DG55}:
\[
\begin{split}
\|u\|_{\overline{q},r}
&\le\left[\frac{2d\,\overline{q}\,\mathcal{S}_2^2}{(2^*-2\overline{q})}+\,\mathcal{S}_2^2\frac{(R_0-R_\infty)^2}{R_\infty^2}\right]^\frac{2^*}{2\overline{q}}
    \frac{\|u\|_{\frac{2}{2^*}\overline{q},R}}{(R-r)^\frac{2^*}{\overline{q}}}\\
&\le\normalcolor\left[\frac{2d\,\overline{q}\,\mathcal{S}_2^2}{(2^*-2\overline{q})}
    +\,\mathcal{S}_2^2\frac{(R_0-R_\infty)^2}{R_\infty^2}\right]^\frac{2^*}{2\overline{q}}
    \left[\frac{\omega_d^{1/d}R}{R-r}\right]^{\frac{2^*}{\overline{q}}}
    \frac{|B_R|^{\frac{1}{\overline{q}}}}{|B_R|^{\frac{1}{q_0}}}\|u\|_{q_0,R}
\end{split}
\]
which is \eqref{iter.neg.2} when $\underline{q}\le q_0\le \overline{q}$, once we let $R=R_0$ and $r=\overline{R}$. On the other hand, when $0<q_0<\underline{q}<\overline{q}$\,, we can use inequality \eqref{DG55} rewritten as
\begin{equation}\label{DG66}
\|u\|_{\overline{q},r}
\le\left[\frac{2d\,\overline{q}\,\mathcal{S}_2^2}{(2^*-2\overline{q})}+\,\mathcal{S}_2^2\frac{(R_0-R_\infty)^2}{R_\infty^2}\right]^\frac{2^*}{2\overline{q}}
\frac{\|u\|_{\frac{2}{2^*}\overline{q},R}}{(R-r)^\frac{2^*}{\overline{q}}}
:=\frac{K}{(R-r)^\frac{2^*}{\overline{q}}}\|u\|_{\frac{2}{2^*}\overline{q},R}
\end{equation}
so that Lemma \ref{Lem.DG.Ext} with $\gamma=2^*/\overline{q}$ gives that for all $0<q_0\le \underline{q}<\overline{q}$ (recall that $\underline{q}=2\overline{q}/2^*$)
\begin{equation}\label{DG.77}\begin{split}
\|u\|_{\overline{q}, R_\infty}
&\le 3\cdot 2^{\frac{\overline{q}(\underline{q}-q_0)}{q_0(\overline{q}-\underline{q})}}\,
        \left[\left(4\gamma\frac{\underline{q}(\overline{q}-q_0)}{q_0(\overline{q}-\underline{q})}
        \right)^\gamma \frac{K}{(R_0-R_\infty)^\gamma}\right]^{\frac{\underline{q}(\overline{q}-q_0)}{q_0(\overline{q}-\underline{q})}}
            \,\|u\|_{q_0,R_0}\\
&= 3\cdot 2^{\frac{(d-2)\overline{q}}{2q_0}-\frac{d}{2}}\,
        \left[\left(4\frac{2d}{(d-2)\overline{q}}\frac{d-2}{d}\overline{q}\frac{(\overline{q}-q_0)}{q_0\frac{2}{d}\overline{q}}
        \right)^{\gamma} \frac{K}{(R_0-R_\infty)^{\gamma}}\right]^{\frac{\underline{q}(\overline{q}-q_0)}{q_0(\overline{q}-\underline{q})}}
            \,\|u\|_{q_0,R_0}\\
&= 3\cdot 2^{\frac{(d-2)\overline{q}}{2q_0}-\frac{d}{2}}\,
        \,K^{\frac{\underline{q}(\overline{q}-q_0)}{q_0(\overline{q}-\underline{q})}}
        \,\left(4d\normalcolor\frac{\overline{q}-q_0}{q_0\overline{q}}\frac{1}{R_0-R_\infty}\right)^{\gamma\frac{\underline{q}(\overline{q}-q_0)}{q_0(\overline{q}-\underline{q})}}
            \,\|u\|_{q_0,R_0}\\
&= 3\cdot 2^{\frac{(d-2)\overline{q}}{2q_0}-\frac{d}{2}}\,
        \,K^{\frac{\overline{q}-q_0}{q_0}\frac{d-2}{2}}
        \,\left(4d\normalcolor\frac{\overline{q}-q_0}{q_0\overline{q}}\frac{\omega_d^{1/d}R_0}{R_0-R_\infty}\right)^{\frac{d}{q_0}-
        \frac{d}{\overline{q}}}
            \frac{|B_{R_0}|^{\frac{1}{\overline{q}}}}{|B_{R_0}|^{\frac{1}{q_0}}}\,\|u\|_{q_0,R_0}\\
&=3\cdot 2^{\frac{(d-2)\overline{q}}{2q_0}-\frac{d}{2}}\,
        \,\left[\frac{2d\,\overline{q}\,\mathcal{S}_2^2}{(2^*-2\overline{q})}+\,\mathcal{S}_2^2\frac{(R_0-R_\infty)^2}{R_\infty^2}
            \right]^{{\frac{\overline{q}-q_0}{\overline{q}\,q_0}\frac{d}{2}}}
        \,\left(4d\normalcolor\frac{\overline{q}-q_0}{q_0\overline{q}}\frac{\omega_d^{1/d}R_0}{R_0-R_\infty}
        \right)^{\frac{d}{q_0}-\frac{d}{\overline{q}}}\\
        &\times \left[\frac{R_\infty}{R_0}\right]^{\frac{d}{q_0}}\normalcolor
            \frac{|B_{R_0}|^{\frac{1}{\overline{q}}}}{|B_{R_\infty}|^{\frac{1}{q_0}}}\,\|u\|_{q_0,R_0}\\
&=3\cdot 2^{\frac{(d-2)\overline{q}}{2q_0}-\frac{d}{2}}\,
        \,\left[\frac{2d\,\overline{q}\,\mathcal{S}_2^2}{(2^*-2\overline{q})}\frac{R_\infty^2}{(R_0-R_\infty)^2}+\,\mathcal{S}_2^2
            \right]^{{\frac{\overline{q}-q_0}{\overline{q}\,q_0}\frac{d}{2}}}
        \,\left(4d\normalcolor\omega_d^{\frac{1}{d}}\frac{\overline{q}-q_0}{q_0\overline{q}}\right)^{\frac{d}{q_0}-\frac{d}
        {\overline{q}}}\\
        &\times\left[\frac{R_\infty}{R_0}\right]^{\frac{d}{q_0}}\normalcolor
            \frac{|B_{R_0}|^{\frac{1}{\overline{q}}}}{|B_{R_\infty}|^{\frac{1}{q_0}}}\,\|u\|_{q_0,R_0}\\
\end{split}
\end{equation}
whence the statement follows upon relabeling $R_\infty$ as ${\overline R}$.
\qed

As a first consequence of the above inequalities, we can improve the local lower bounds of Theorem \ref{thm.local.lower} in this good supercritical range.

\begin{thm}[Local Lower Estimates when $1<p<p_c$]\label{thm.local.lower.pc}
Let $\Omega\subseteq\RR^d$ and let $\lambda>0$. Let $u$ be a nonnegative local weak supersolution in $B_{R_0}\subseteq\Omega$ to $-\Delta u = \lambda u^p$, with $1< p<p_c=d/(d-2)$.
\begin{equation}\label{lower.est.pc}
\inf_{x\in B_{R_\infty}}u(x) =\|u\|_{-\infty,R_\infty}
    \ge\frac{I_{-\infty,{\underline q}}}{I_{\overline{q},{\underline q}}}\frac{\|u\|_{\overline{q},\overline{R}}}{|B_{\overline{R}}|^{\frac{1}{\overline{q}}}} \qquad\mbox{with}\qquad d(p-1)/2<\overline{q}<d/(d-2)
\end{equation}
for any $0<R_\infty<\overline{R}<R_0$, where ${\underline q}\in(0,q_0\wedge\overline{q}]$, $q_0$ and
$I_{-\infty,{\underline q}}$ are given in \eqref{Const.lower.local.pc.q0} and $I_{\overline{q},{\underline q}}$ is given by \eqref{I.pc.lower}.
\end{thm}
\noindent {\sl Proof.~}We use the local lower bounds of Theorem \ref{thm.local.lower} for $\underline{q}\in(0,q_0]$, $\varepsilon=\ee$, with the definition of $q_0$ to be recalled below, so that
\begin{equation}\label{lower.local.pc.q0}
\inf_{x\in B_{R_\infty}}u(x) =\|u\|_{-\infty,R_\infty}
    \ge I_{-\infty,{\underline q}}\frac{\|u\|_{{\underline q},R_{0}}}{|B_{R_0}|^{\frac{1}{{\underline q}}}}.
\end{equation}
where
\begin{equation}\label{Const.lower.local.pc.q0}
{\underline q}\le q_0:=\frac{2^{{\frac{d-3}{2}}}}{d^2\omega_d^2\ee}\quad\mbox{ and }\quad I_{-\infty,{\underline q}}
=\left[2^d\mathcal{S}_2^2\left(\frac{dR_0^2}{(R_0-R_\infty)^2}+\frac{R_0^2}{R_\infty^2}\right)
    \right]^{-\frac{d}{2{\underline q}}}
    \left[\frac{\ee}{2^d \,\ee\,(d+1)\,\sqrt{\omega_d} }\right]^{\frac{2}{{\underline q}}}.
\end{equation}
Recall the reverse H\"older inequalities of Proposition \ref{rev.holder.pc}
\begin{equation}\label{lower.q.pc}
\frac{\,\|u\|_{{\underline q},R_0}}{|B_{R_0}|^{\frac{1}{{\underline q}}}}\ge \frac{\|u\|_{\overline{q},\overline{R}}}{I_{\overline{q},{\underline q}}|B_{\overline{R}}|^{\frac{1}{\overline{q}}}} \qquad\mbox{for any $0<\overline{R}<R_0$ } {\underline q}\in(0,{\overline q}] \mbox{ and } d(p-1)/2<\overline{q}<d/(d-2)
\end{equation}
where
\begin{equation}\label{I.pc.lower}
I_{\overline{q},{\underline q}}:=\left\{\begin{array}{lll}
\left[\frac{2d\,\overline{q}\,\mathcal{S}_2^2}{(2^*-2\overline{q})}
    +\,\mathcal{S}_2^2\frac{(R_0-\overline{R})^2}{\overline{R}^2}\right]^\frac{2^*}{2\overline{q}}
    \left[\frac{\omega_d^{1/d}R_0}{R_0-\overline{R}}\right]^{\frac{2^*}{\overline{q}}}
    \left[\frac{R_0}{\overline{R}}\right]^{\frac{d}{\overline{q}}}
    & \mbox{ if }\frac{d-2}{d}\overline{q}\le {\underline q}\le \overline{q} \,, \\
3\cdot 2^{\frac{(d-2)\overline{q}}{2{\underline q}}-\frac{d}{2}}\,
        \,\left[\frac{2d\,\overline{q}\,\mathcal{S}_2^2}{(2^*-2\overline{q})}\frac{\overline{R}^2}{(R_0-\overline{R})^2}+\,\mathcal{S}_2^2
            \right]^{{\frac{\overline{q}-{\underline q}}{\overline{q}\,{\underline q}}\frac{d}{2}}}
        \,\left(4\omega_d^{\frac{1}{d}}\frac{\overline{q}-{\underline q}}{{\underline q}\overline{q}}\right)^{\frac{d}{{\underline q}}-\frac{d}{\overline{q}}}
            \left[\frac{R_0}{\overline{R}}\right]^{\frac{d}{\overline{q}}}\,,
    & \mbox{ if }0<{\underline q}<\frac{d-2}{d}\overline{q}\\
\end{array}\right.
\end{equation}
with $q_0$ as in \eqref{Const.lower.local.pc.q0}. Combining inequalities \eqref{lower.local.pc.q0} and \eqref{lower.q.pc} we obtain \eqref{lower.est.pc}.\qed

\noindent\textbf{Remark. }The above lower bounds turn our to be important when applied to solutions, since they will imply directly a clean form of Harnack inequality when $1<p<p_c$ and then local absolute bounds, which is a novelty and a typical feature of the ``good'' superlinear case $1<p<p_c$. We stress the fact that in the upper range $p_c\le p<p_s$ such absolute bounds can not be true, as explicit counter-examples show. We will give more details on these counterexamples in the next section.


\section{Harnack inequalities}\label{HI}
In this section we will show in a quantitative way how upper and lower bounds can be joined to form Harnack inequalities for solutions, and to obtain as a consequence absolute local upper ($1<p<p_c$) and absolute local lower bounds ($0<p<1$), which are new, as far as we know. We first join local bounds of Theorems \ref{thm.local.upper}, \ref{Thm.local.upper.b.p} (upper) and \eqref{thm.local.lower} (lower), to obtain a general form for Harnack inequalities, which at a first sight appear to be weaker than what expected, because its constant depends on local $\LL^q$-norms of the solution itself. This is the only form of Harnack inequality that can hold for all $0\le p<p_s=(d+2)/(d-2)$. To eliminate this quotient and to obtain Harnack inequalities in a more classical form one has to assume that $0<p<p_c=d/(d-2)$.

This fact might seem puzzling, but there are very weak (distributional) solutions in the range $p_c \le p<p_s$ that are not bounded, cf. \cite{MP1, P1,P2,P3,P4}, even when one prescribes zero Dirichlet boundary conditions. According to Mazzeo and Pacard \cite{MP1}, in this range there are solutions with a singularity of the type $|x-x_0|^{-2/(p-1)}$ at a point $x_0\in\Omega$. Such solutions are not locally in $\LL^q$ with $q>d(p-1)_+/2$ if  $p>p_c$, hence the local upper estimate fails for them when applied to a ball that contains the singularity. In this range there appears in a clear form the difference between weak and very weak solutions, which helps  understanding these critical exponents. Regarding boundary behaviour, the range to  consider is $p_1\le p<p_s$, where $p_1=(d+1)/(d-1)$ is the exponent introduced by Brezis and Turner \cite{BT}. In this range there exist very weak solutions which are not weak (energy) solutions and can have a singularity at some points of the boundary and satisfy elsewhere on the boundary the prescribed condition in a suitable trace sense, not necessarily in a continuous fashion, cf. del Pino et al. \cite{dPMP}.

\begin{thm}[Harnack inequality for $0\le p <p_s$]\label{Thm.Harnack.ps}
Let $\Omega\subseteq\RR^d$ and let $\lambda>0$. Let $u$ be a nonnegative local weak solution in $B_{R_0}\subseteq\Omega$ to $-\Delta u = \lambda u^p$, with $0\le  p < p_s=(d+2)/(d-2)$. Given $R_\infty<R_0$ and $\varepsilon>0$ we assume
\begin{equation}\label{q0-figure.h}
0<\underline{q}\le q_0:=\frac{2^{{\frac{d-3}{2}}}}{d\omega_d^2[\ee(d-1)+\varepsilon]},\quad \overline{q} > \frac{d(p-1)_+}{2}.
\end{equation}
If $0<\overline{q}<d/(d-2)$ we also assume
\[\left[\frac{\log\frac{2^*-d(p-1)_+}{2\overline{q}-d(p-1)_+}}{\log\frac{d}{d-2}}\right]
\textrm{not integer}.
\]
Then the following bound holds true
\begin{equation}\label{Harnack.ps}
\sup_{x\in B_{R_\infty}}u(x)\le \mathcal{H}_p[u]\inf_{x\in B_{R_\infty}}u(x)
\end{equation}
where $\mathcal{H}_p[u]$ depends on $u$ through some local norms as follows
\begin{equation}\label{Harnack.ps.const}
\mathcal{H}_p[u]=\mathcal{H}_p[u](d,\overline{q},\underline{q},\varepsilon, R_0, R_\infty)
= \frac{ I_{\infty,\overline{q}}}{I_{-\infty,\underline{q}}} \;
    \left(\frac{\left(\fint_{B_{R_0}} u^{q}\dx\right)^{\frac{(p-1)_+}{q}}}{   {\fint_{B_{R_\infty}}u^{(p-1)_+}\dx}}\right)^{\frac{d}{{2\overline{q}-d(p-1)_+}}}
    \frac{\left(\fint_{B_{R_0}} u^{\overline{q}}\dx\right)^{\frac{1}{\overline{q}}} }{\left(\fint_{B_{R_0}} u^{\underline{q}}\dx\right)^{\frac{1}{\underline{q}}}}\,.
\end{equation}
with $I_{\infty,\overline{q}}$ given by \eqref{Const.Upper.q},  $I_{-\infty,\underline{q}}$ is given by \eqref{Const.lower.local.thm}.
\end{thm}
\noindent {\sl Proof.~} We recall the local upper bounds of Theorem \ref{thm.local.upper}: for any $B_{R_\infty}\subset B_{R_0}\subseteq\Omega$
\begin{equation}\label{upper.harnack}
\|u\|_{\infty,R_\infty}
    \le I_{\infty,\overline{q}}\;
    \left(\frac{\left(\fint_{B_{R_0}} u^{q}\dx\right)^{\frac{(p-1)_+}{q}}}{   {\fint_{B_{R_\infty}}u^{(p-1)_+}\dx}}\right)^{\frac{d}{{2\overline{q}-d(p-1)_+}}}\frac{\|u\|_{\overline{q},R_{0}}}{|B_{R_0}|^{\frac{1}{\overline{q}}}}
    \qquad\mbox{for any}\qquad
    \overline{q} > \frac{d(p-1)_+}{2}
\end{equation}
where $I_{\infty,\overline{q}}$ is given by \eqref{Const.Upper.q} and when $0<\overline{q}<d/(d-2)$ we require the additional condition \eqref{intero} on $\overline{q}$.
We also recall the lower bounds of Theorem \ref{thm.local.lower}: for any $\varepsilon>0$ and for any $\underline{q}$ as in \eqref{q0-figure.h}, the following bound holds true
\begin{equation}\label{lower.harnack}
\frac{\inf_{x\in B_{R_\infty}}u(x)}{I_{-\infty,\underline{q}}}  \frac{|B_{R_0}|^{\frac{1}{\underline{q}}}}{\|u\|_{{\underline{{q}}},R_{0}}}
    \ge 1.
\end{equation}
where $I_{-\infty,\underline{q}}$ is given by \eqref{Const.lower.local.thm}. Joining \eqref{upper.harnack} and \eqref{lower.harnack} gives
\eqref{Harnack.ps}\,.\qed

\begin{thm}[Harnack inequality, $0\le p \le 1$]\label{Harnack.p<1}
Let $\Omega\subseteq\RR^d$ and let $\lambda>0$. Let $u$ be a nonnegative local weak solution in $B_{R_0}\subseteq\Omega$ to $-\Delta u = \lambda u^p$, with $0\le p\le 1$. For all $R_\infty<R_0$ the following bound holds true
\[
\sup_{x\in B_{R_{\infty}}}u(x)\le \mathcal{H}_p\inf_{x\in B_{R_{\infty}}}u(x)
\]
where $\mathcal{H}_p$ does not depend on $u$\,, and is given by
\begin{equation}\label{Harnack.sublin.const}
\begin{split}\mathcal{H}_p
&=\left[\frac{2^d\mathcal{S}_2^4 R_0^2}{(R_0-R_\infty)^2}\left(\frac{dR_0^2}{(R_0-R_\infty)^2}+\frac{R_0^2}{R_\infty^2}\right)
    \right]^{\frac{d}{2q_0}}
    \left[\frac{2^d \,\left(\left(\frac{d}{d-2}\right)^{n_0-\frac{1}{2}}\frac{2^{{\frac{d-3}{2}}}}{d\omega_d^2}+\ee\right)\,\sqrt{\omega_d} }{\left(\frac{d}{d-2}\right)^{n_0-\frac{1}{2}}\frac{2^{{\frac{d-3}{2}}}}{d\omega_d^2}-e(d-1)}\right]^{\frac{2}{q_0}}\times\\
    &\times\left\{\left(\frac{d}{d-2}\right)^d \frac{2(d-2)\sqrt{d}}{\big(\sqrt{d}-\sqrt{d-2}\big)^3}\left[\Lambda_p+\frac{d-2}{q_0}+\frac{(R_0-R_\infty)^2}{R_\infty^2}
        \max\left\{\frac{d-2}{(dq_0)^2}|dq_0-(d-2)|,\,\frac14\right\}\right]\right\}^{\frac{d}{2q_0}}\\
\end{split}
\end{equation}
with
\begin{equation}\label{q0-figure.h.1}
q_0=\left(\frac{d-2}{d}\right)^{n_0-\frac{1}{2}}\qquad\mbox{and}\qquad n_0=i.p.\left[\frac{\log\left(\ee(d-1)\frac{d\omega_d^2}{2^{{\frac{d-3}{2}}}}\right)}{\log\frac{d}{d-2}}\,+\frac{3}{2}\right]
\end{equation}
\end{thm}
\noindent {\sl Proof.~}The goal of the proof is to simplify the quotient of $\LL^q$-norms in the expression of the constant $\mathcal{H}_p[u]$ of the Harnack inequality \eqref{Harnack.ps}. Since we are dealing with the range $0\le p\le 1$, we can choose
\[
0<\overline{q}=\underline{q}=q_0=q_0(\varepsilon):=\frac{2^{{\frac{d-3}{2}}}}{d\omega_d^2[\ee(d-1)+\varepsilon]}\quad\mbox{and}\quad \overline{q} > 0
\quad\mbox{with}\quad
\left[\frac{\log\frac{2^*}{2\overline{q}}}{\log\frac{d}{d-2}}\right]
\textrm{not integer}.
\]
In fact, we shall arrive, with a suitable choice of the parameter $\varepsilon$, to a value of $q_0$ smaller than $d/(d-2)$, so that the requirement $\left[\log\frac{2^*}{2\overline{q}}/{\log\frac{d}{d-2}}\right]$ not being integer is necessary. The latter condition means $q_0(\varepsilon)\neq [(d-2)/d]^n$ for all $n\in \mathbb{N}$, and this is possible since we can always choose $\varepsilon$
\[
0<\varepsilon=\left(\frac{d}{d-2}\right)^{n_0-\frac{1}{2}}\frac{2^{{\frac{d-3}{2}}}}{d\omega_d^2}-e(d-1)\qquad\mbox{so that}\qquad
q_0=\left(\frac{d-2}{d}\right)^{n_0-\frac{1}{2}}
\]
where $n_0$ is the first integer $n$ such that $\varepsilon(n)>0$, which is
\[
n_0=i.p.\left[\frac{\log\left(\ee(d-1)\frac{d\omega_d^2}{2^{{\frac{d-3}{2}}}}\right)}{\log\frac{d}{d-2}}\,+\frac{1}{2}\right]+1.
\]
The constants become in this case
\begin{equation}\label{Const.Upper.q.harnack}\begin{split}
I_{\infty,\overline{q}}
&=\left[\frac{c_1\mathcal{S}_2^2 R_0^2}{(R_0-R_\infty)^2}\right]^{\frac{d}{2q_0}}
\left\{\left(\frac{d}{d-2}\right)^d \frac{2(d-2)}{\big(\sqrt{d}-\sqrt{d-2}\big)^2}\right.\times\\
&\times\left.\left[\Lambda_p+\frac{d-2}{q_0}+\frac{(R_0-R_\infty)^2}{R_\infty^2}
        \max\left\{\frac{d-2}{(dq_0)^2}|dq_0-(d-2)|,\,\frac14\right\}\right]\right\}^{\frac{d}{2q_0}}\\
&=\left[\frac{\mathcal{S}_2^2 R_0^2}{(R_0-R_\infty)^2}\right]^{\frac{d}{2q_0}}
\left\{\left(\frac{d}{d-2}\right)^d \frac{2(d-2)\sqrt{d}}{\big(\sqrt{d}-\sqrt{d-2}\big)^3}\right.\times\\
&\times\left.\left[\Lambda_p+\frac{d-2}{q_0}+\frac{(R_0-R_\infty)^2}{R_\infty^2}
        \max\left\{\frac{d-2}{(dq_0)^2}|dq_0-(d-2)|,\,\frac14\right\}\right]\right\}^{\frac{d}{2q_0}},\\
\end{split}
\end{equation}
where $\Lambda_p=2$ if $p\ne 1$, $\Lambda_p= \lambda/4$ if $p=1$ and, since $q_0<d/(d-2)$,
\begin{equation}\label{c1.q.h}
c_1:=\max\limits_{i=0,1}\frac{q_0\left(\frac{d}{d-2}\right)^{k_0-1+i}}
{\left|q_0\left(\frac{d}{d-2}\right)^{k_0-1+i}-1\right|}
=\max\limits_{i=0,1}\frac{\left(\frac{d}{d-2}\right)^{i+\frac{1}{2}}}
{\left(\frac{d}{d-2}\right)^{i+\frac{1}{2}}-1}
=\frac{\sqrt{d}}
{\sqrt{d}-\sqrt{d-2}}
\end{equation}
since $k_0$ is given by:
\[
k_0=i.p.\left[\frac{\log\frac{2^*}{2q_0}}{\log\frac{d}{d-2}}\right]
=i.p.\left[1+\frac{\log\frac{1}{q_0}}{\log\frac{d}{d-2}}\right]
=i.p.\left[1+n_0-\frac{1}{2}\right]=n_0+1
\]
and the last step in \eqref{c1.q.h} follows by an explicit calculation. Moreover $I_{-\infty,\underline{q}}$ given by formula \eqref{Const.lower.local.thm} takes the form
\[\begin{split}
I_{-\infty,q_0}
    &=\left[2^d\mathcal{S}_2^2\left(\frac{dR_0^2}{(R_0-R_\infty)^2}+\frac{R_0^2}{R_\infty^2}\right)
    \right]^{-\frac{d}{2q_0}}
    \left[\frac{\varepsilon}{2^d \,(\ee\,d+\varepsilon)\,\sqrt{\omega_d} }\right]^{\frac{2}{q_0}}\\
    &=\left[2^d\mathcal{S}_2^2\left(\frac{dR_0^2}{(R_0-R_\infty)^2}+\frac{R_0^2}{R_\infty^2}\right)
    \right]^{-\frac{d}{2q_0}}
    \left[\frac{\left(\frac{d}{d-2}\right)^{n_0-\frac{1}{2}}\frac{2^{{\frac{d-3}{2}}}}{d\omega_d^2}-e(d-1)}{2^d \,\left(\left(\frac{d}{d-2}\right)^{n_0-\frac{1}{2}}\frac{2^{{\frac{d-3}{2}}}}{d\omega_d^2}+\ee\right)\,\sqrt{\omega_d} }\right]^{\frac{2}{q_0}}.
\end{split}
\]
Hence we get the expression of $\mathcal{H}_p=I_{\infty,q_0}/I_{-\infty,q_0}$ given in \eqref{Harnack.sublin.const}\,.\qed

When $p>1$ we can not join the upper and the lower bound so easily, we need the improved lower bounds of Theorem \ref{thm.local.lower.pc}, valid only when $p<p_c$.

\begin{thm}[Harnack Inequalities when $1<p<p_c$]\label{Harnack.p<pc}
Let $\Omega\subseteq\RR^d$ and let $\lambda>0$. Let $u$ be a nonnegative local weak solution to $-\Delta u = \lambda u^p$ in $B_{R_0}\subseteq\Omega$, with $1< p<p_c=d/(d-2)$. Then for any $0<R_\infty<\overline{R}<R_0$ there exists an explicit constant $\mathcal{H}_p>0$ such that
\begin{equation}\label{Harnack.est.p<pc}
\sup_{x\in B_{R_\infty}}u(x)\le \mathcal{H}_p\inf_{x\in B_{R_\infty}}u(x)
\end{equation}
where $\mathcal{H}_p$ does not depend on $u$\,, and is given by
\begin{equation}\label{Harnack.const.p<pc}
\mathcal{H}_p=I_{\infty,\overline{q}}\;
        \left( \frac{I_{\overline{q},{\underline q}}}{I_{-\infty,{\underline q}}}   \right)^{\frac{2\overline{q}}{{2\overline{q}-d(p-1)}}}\,,
        \qquad\mbox{with}\qquad \frac{d(p-1)}{2}<\overline{q}<\frac{d}{d-2}
\end{equation}
where the constants ${\underline q}\in(0,q_0\wedge\overline{q}]$, $q_0$ and
$I_{-\infty,{\underline q}}$ are given in \eqref{Const.lower.local.pc.q0}, $I_{\overline{q},{\underline q}}$ is given by \eqref{I.pc.lower}, $I_{\infty,\overline{q}}$ is given by \eqref{Const.Upper.q}; moreover, since $\overline{q}<d/(d-2)$ we require the additional condition \eqref{intero}.
\end{thm}
\noindent {\sl Proof.~}We first consider the lower bounds of Theorem \ref{thm.local.lower.pc}.
Let $\Omega\subseteq\RR^d$ and let $\lambda>0$. Let $u$ be a nonnegative local weak supersolution in $B_{R_0}\subseteq\Omega$ to $-\Delta u = \lambda u^p$, with $1< p<p_c=d/(d-2)$. Then
\begin{equation}\label{lower.est.pc.h}
\frac{\|u\|_{\overline{q},\overline{R}}}{|B_{\overline{R}}|^{\frac{1}{\overline{q}}}}\le \frac{I_{\overline{q},{\underline q}}}{I_{-\infty,{\underline q}}}\inf_{x\in B_{R_\infty}}u(x)
\end{equation}
for any $0<R_\infty<\overline{R}<R_0$, where $d(p-1)/2<\overline{q}<d/(d-2)$, ${\underline q}\in(0,q_0\wedge\overline{q}]$, $q_0$ and
$I_{-\infty,{\underline q}}$ are given in \eqref{Const.lower.local.pc.q0} and $I_{\overline{q},{\underline q}}$ is given by \eqref{I.pc.lower}.
Then we recall the upper bounds of Theorem \ref{thm.local.upper} which we rewrite as
\begin{equation}\label{upper.p+1.h}\begin{split}
\|u\|_{\infty,R_\infty}
 &\le I_{\infty,\overline{q}}\;
    \left(\frac{\|u\|_{\overline{q},\overline{R}}^{(p-1)_+}}{|B_{\overline{R}}|^{\frac{(p-1)_+}{\overline{q}}}}\frac{|B_{R_\infty}|}{\int_{B_{R_\infty}}u^{(p-1)_+}\dx}
    \right)^{\frac{d}{{2\overline{q}-d(p-1)_+}}}\frac{\|u\|_{\overline{q},\overline{R}}}{|B_{\overline{R}}|^{\frac{1}{\overline{q}}}}\\
 &\le I_{\infty,\overline{q}}\;
    \left(\frac{\|u\|_{\overline{q},\overline{R}}}{|B_{\overline{R}}|^{\frac{1}{\overline{q}}}}
    \frac{1}{\inf_{x\in B_{R_\infty}}u(x)  }
    \right)^{\frac{d(p-1)}{{2\overline{q}-d(p-1)}}}\frac{\|u\|_{\overline{q},\overline{R}}}{|B_{\overline{R}}|^{\frac{1}{\overline{q}}}}\\
 &\le I_{\infty,\overline{q}}\;
        \left( \frac{I_{\overline{q},{\underline q}}}{I_{-\infty,{\underline q}}}   \right)^{\frac{d(p-1)}{{2\overline{q}-d(p-1)}}+1}
        \inf_{x\in B_{R_\infty}}u(x)
\end{split}\end{equation}
for any $\overline{q} > \frac{d(p-1)_+}{2}$, where $I_{\infty,\overline{q}}$ is given by \eqref{Const.Upper.q} and since $0<\overline{q}<d/(d-2)$ we require the additional condition \eqref{intero}. In the third step we have used the lower bound \eqref{lower.est.pc.h}. \qed

\noindent\textbf{Remark. } Notice that the constant $\mathcal{H}_p$ does not depend on $u$ in the range $0\le p<p_c$\,, and it does not depend on $\lambda>0$ when moreover $p\ne 1$.


\section{Local Absolute bounds }\label{sect.absolute}

In this section we will prove local absolute lower bounds when $0<p<1$ and local absolute upper bounds when $1<p<p_c$ as a consequence of the Harnack inequalities of the previous section together with the Caccioppoli estimates \eqref{Abs.Upper.p-1}.
\begin{thm}[Local Absolute bounds]\label{Thm.Absolute.0.pc}
Let $\Omega\subseteq\RR^d$ and let $\lambda>0$. Let $u$ be a local nonnegative weak solution to $-\Delta u = \lambda u^p$ in $B_{R_0}\subseteq\Omega$, with $0< p<p_c=d/(d-2)$. Then for any $0<R_\infty<\overline{R}<R_0$ there exists a constant $\mathcal{H}_p>0$ that does not depend on $u$, such that
\begin{equation}\label{Absolute.upper.p<pc}
\sup_{x\in B_R(x_0)}u(x)\le \mathcal{H}_p \left(\frac{8R_0^d}{\lambda(R_0-R)^2 R^d}\right)^{\frac{1}{p-1}}\qquad\mbox{when }1<p<p_c=\frac{d}{d-2}\,,
\end{equation}
and, if $u\not\equiv 0$ on $B_{R_0}$
\begin{equation}\label{Absolute.lower.p<1}
\inf_{x\in B_R(x_0)}u(x)\ge \mathcal{H}_p^{-1} \left(\frac{\lambda(R_0-R)^2 R^d}{8R_0^d}\right)^{\frac{1}{1-p}}\qquad\mbox{when }0<p<1\,.
\end{equation}
The constant $\mathcal{H}_p$ is given by \eqref{Harnack.sublin.const} when $0<p<1$ and by \eqref{Harnack.const.p<pc}  when $1<p<p_c$.
\end{thm}
\noindent {\sl Proof.~}We combine the quantitative Harnack inequalities of Theorems \ref{Harnack.p<1} and \ref{Harnack.p<pc} together with the quantitative Caccioppoli estimates \eqref{Abs.Upper.p-1}
\begin{equation*}
\lambda\int_{B_R} u^{p-1}\dx\le\frac{8\omega_d R_0^d}{(R_0-R)^2}\,
\end{equation*}
which implies, when $p>1$,
\begin{equation}\label{Cacc.h.p>1}
\inf_{x\in B_R}u(x)\le \left(\frac{1}{|B_R|}\int_{B_R} u^{p-1}\dx\right)^{\frac{1}{p-1}}
\le\left(\frac{8R_0^d}{\lambda(R_0-R)^2 R^d}\right)^{\frac{1}{p-1}}
\end{equation}
and when $0<p<1$ as
\begin{equation}\label{Cacc.h.p<1}
\left(\frac{\lambda(R_0-R)^2R^d}{8 R_0^d}\right)^{\frac{1}{1-p}}
\le \left(\frac{|B_R|}{\int_{B_R} u^{p-1}\dx}\right)^{\frac{1}{1-p}}
\le \left(\frac{1}{\sup\limits_{x\in B_R}u(x)^{p-1}}\right)^{\frac{1}{1-p}}=\sup_{x\in B_R}u(x)
\end{equation}
The above inequalities can be now combined with the corresponding Harnack inequalities of Theorems \ref{Harnack.p<1} and \ref{Harnack.p<pc}, which have the form
\[
\sup_{x\in B_R}u(x)\le \mathcal{H}_p \inf_{x\in B_R}u(x)
\]
to obtain the desired bounds in both cases. The constant $\mathcal{H}_p$ is given by \eqref{Harnack.sublin.const} when $0<p<1$ and by \eqref{Harnack.const.p<pc} $1<p<p_c$.\qed

\noindent\textbf{Remark. } These bounds are new as far as we know. Notice that they depend explicitly on $\lambda$.


\section{Regularity. Local bounds for the gradients}

In this section we will prove $\LL^\infty$ bounds for the gradients, to conclude that solutions to $-\Delta u=\lambda u^p$ are indeed local Lipschitz functions. The strategy to prove such results is to show that the incremental quotients $u_{h,i}$ satisfy the equation $-\Delta u_{h,i}\le b(x)u_{h,i}$ for a suitable $b(x)$, so that we can apply the local $\LL^\infty$ bounds of Theorem \ref{Thm.local.upper.b.DG}\,.

\noindent \textbf{Short reminder about incremental quotients in $W^{1,q}$.} Here we follow Giusti \cite{Giusti}. \normalcolor
It is well known that if $u\in W^{1,q}(\Omega)$ then its incremental quotients is defined as
\[
u_{h,i}:=\frac{u(x+h e_i)-u(x)}{h}
\]
where $e_i$ denotes the unit vector in the direction $x_i$\,, cf. \cite{Evans,Giusti}\,. Let us recall some properties of the incremental quotients:\\
\noindent(i) If $u\in W^{1,q}(\Omega)$\,, then its incremental quotient $u_{h,i}$ is defined in the set
\[
\Omega_{|h|}:=\left\{x\in\Omega\;\big|\; \dist(x,\partial\Omega)>|h|\right\}\,,\qquad\mbox{moreover}\qquad
u_{h,i}\in W^{1,q}(\Omega_{|h|})\,.
\]

\noindent(ii) If $u\in W^{1,q}(\Omega)$ for $1\le q\le\infty$ and $\Sigma\subset\subset\Omega$\,, then for any $|h|<\dist(\Sigma,\Omega)/(10 \sqrt{d})$ we have
\begin{equation}\label{incr.1}
\|u_{h,i}\|_{\LL^q(\Sigma)}\le 5^\frac{d}{q} \, \|\partial_i u\|_{\LL^q(\Omega)}\,.
\end{equation}
for a proof of the latter fact we refer to Lemma 8.1 of \cite{Giusti}\,.

\noindent(iii) Let $u\in\LL^q(\Omega)$, $1<q<\infty$\,, and assume that there is a constant $K$ such that for every $h$ small enough we have
$\|u_{h,i}\|_{\LL^q(\Omega_{|h|})}\le K$. Then $\partial_i u\in\LL^q(\Omega)$ and $\|\partial_i u\|_{\LL^q(\Omega)}\le K$\,. Moreover $u_{h,i}\to \partial_i u$ in $\LL_{\rm loc}^q(\Omega)$ as $h\to 0$\,. For a proof of this fact we refer to Lemma 8.2 of \cite{Giusti}\,.

\medskip

We can now state and prove the following theorem.
\begin{thm}[Local upper bounds for the gradient]\label{thm.local.upper.grad}
Let $\Omega\subseteq\RR^d$ and let $\lambda>0$. Let $u$ be a local nonnegative weak solution to $-\Delta u = \lambda u^p$ in $B_{R_0}\subseteq\Omega$, with $0< p<p_c=d/(d-2)$. Then for any $0<R_\infty<R_0$ we have
\begin{equation}\label{Upper.grad.u.thm}
\|\nabla u\|_{\infty, R_\infty}
\le K[u]\,\|u\|_{2, R_0}
\end{equation}
where
\begin{equation}\label{Const.Upper.uhi.2}
\begin{split}
K[u]&= \left(\frac{15}{R_0-R_\infty}\right)^{\frac{d}{2}}
\left[\lambda\,b_{p,R_0}[u] +\frac{18d}{(R_0-R_\infty)^2} \right]^{\frac{1}{2}}\left(\frac{2d^d}{2^d}\right)^{\frac{d^2}{8}}\\
&\times
         \left[16(d+2)+\frac{(R_0-R_\infty)^2}{9R_\infty^2}
        +\left(\frac{d\,\mathcal{S}_2^2 (p\vee 1)}{d-2}\right)^{\frac{d}{2}}
        \frac{4(R_0-R_\infty)^2}{9(d-2)}\,|B_{R_0} |^{\frac{d-2}{d}}\,\left(\lambda b_{p,R_0}[u]\right)^{\frac{d}{2}}
         \right]^{\frac{d}{4}}\\
\end{split}
\end{equation}
with
\begin{equation}\label{bp0}
b_{p,R_0}[u]\le \left\{
\begin{array}{lll}
1\,,& \mbox{if } p=1\,,\\
\dfrac{8R_0^d\,\mathcal{H}_p^{|p-1|}}{\lambda(R_0-R_\infty)^2 R_\infty^d}\,, & \mbox{if } 0\le p< p_c\mbox{ and }p\neq 1\,,\\[4mm]
\|u\|_{\infty, R_0}^{p-1}\,, & \mbox{if } p_c\le p< p_s\,,\\
\end{array}
\right.
\end{equation}
where the constant $\mathcal{H}_p$ is given by \eqref{Harnack.sublin.const} when $0<p<1$ and by \eqref{Harnack.const.p<pc} when $1<p<p_c$.
\end{thm}
\noindent {\sl Proof.~}The proof is divided into several steps.  We start fixing $h_0>0$ small enough.

\noindent$\bullet~$\textsc{Step 1. }\textit{The equation satisfied by the incremental quotients. }First we deduce formally the equation for the positive and negative part, then we justify it rigorously at the end of this step, using Kato's inequality. If $u$ is a solution to $-\Delta u=\lambda u^p$, then the equation satisfied by $u_{h,i}^+$ is
\begin{equation}\label{eq.uhi+}
-\Delta u_{h,i}^+= b^+(x,h)\,u_{h,i}^+\le \lambda\, (p\vee1)\, b_{p}\, u_{h,i}^+\,,\quad\mbox{for all } |h|\le h_0\,,
\end{equation}
where
\begin{equation}\label{bp}
b_p=b_{p,R_0}[u]:=\left\{
\begin{array}{lll}
\sup\limits_{B_{R_0}}u^{p-1} & \mbox{if } 1\le p< p_s\\[5mm]
\left[\inf\limits_{B_{R_0}}u^{1-p}\right]^{-1} & \mbox{if } 0\le p< 1\\
\end{array}
\right.
\end{equation}
and we observe that $b_{p,R}[u]\le b_{p,R_0}[u]$ for any $0<R<R_0$. Indeed, when $u_{h,i}\ge 0$:
\[\begin{split}
-\Delta u_{h,i}^+= \lambda\frac{u^p(x+he_i)-u^p(x)}{h}
&=\lambda\frac{u^p(x+he_i)-u^p(x)}{u(x+he_i)-u(x)}\frac{u(x+he_i)-u(x)}{h} := b^+(x,h)\,u_{h,i}^+\\
&\le \lambda(p\vee1)\,\max\big\{u^{p-1}(x+he_i), u^{p-1}(x)\big\}u_{h,i}^+\\
\end{split}
\]
by using the numerical inequality \eqref{num.ineq.abp}, namely $(a-c)(a^p-c^p)\le (p\vee1)\,\max\big\{a^{p-1}, c^{p-1}\big\}(a-c)^2$ valid for any $p> 0$ and all $a,c\ge 0$ to estimate
\[
b^+(x,h)=\lambda\frac{u^p(x+he_i)-u^p(x)}{u(x+he_i)-u(x)}\le \lambda(p\vee1)\,\max\big\{u^{p-1}(x+he_i), u^{p-1}(x)\big\}\\
\]
we have used the fact that $u^p(x+he_i)-u^p(x)$ and $u(x+he_i)-u(x)$ have the same sign.\\
When $p\ge 1$ we have
\[
-\Delta u^+_{h,i}=b^+(x,h)\,u_{h,i}^+ \le  \lambda(p\vee1)\,  \sup\limits_{B_{R+h_0}}\big(u^{p-1}\big)\,u^+_{h,i}\,,
\]
while when $0\le p<1$ we have
\[
-\Delta u^+_{h,i}=b^+(x,h)\,u_{h,i}^+ \le  \frac{\lambda(p\vee1)}{\inf\limits_{B_{R+h_0}}u^{1-p}}\,u^+_{h,i}
\]
On the other hand, if $u$ is a solution to $-\Delta u=\lambda u^p$, then the equation satisfied by $u_{h,i}^-$ is
\begin{equation}\label{eq.uhi-}
-\Delta u_{h,i}^-= b^-(x,h)\,u_{h,i}^-\le \lambda\, (p\vee1)\, b_{p}\, u_{h,i}^-\,,\quad\mbox{for all } |h|\le h_0\,,
\end{equation}
where $b_p$ is given by \eqref{bp}. Indeed when $u_{h,i}\le 0$ we have that
\[\begin{split}
-\Delta u_{h,i}^-= -\lambda\frac{u^p(x+he_i)-u^p(x)}{h}
&=-\lambda\frac{u^p(x+he_i)-u^p(x)}{u(x+he_i)-u(x)}\frac{u(x+he_i)-u(x)}{h} := b^-(x,h)\,u_{h,i}^-\\
&\le \lambda(p\vee1)\,\max\big\{u^{p-1}(x+he_i), u^{p-1}(x)\big\}u_{h,i}^-\\
\end{split}
\]
for the same arguments as above. Now it remains to justify the formal calculations made above. First we recall Kato's inequality: if $j:\RR\to\RR$ is a convex function such that $j(0)=0$, $j'(v)>0$ if $v>0$, then $\Delta j(v)\ge j'(v)\Delta v$, in the weak sense, whenever $\Delta v\in \LL^1_{\rm loc}(\Omega)$\,. Consider a sequence of convex function $j_\varepsilon$ that approximate $j(u_{h,i})=u_{h,i}^+$ and such that $j_\varepsilon(0)=0$, $j'_\varepsilon(u_{h,i})>0$ if $u_{h,i}>0$. Then by Kato's inequality, we have that indeed $u^+$ satisfy the weak formulation
\[\begin{split}
\int_K \nabla\varphi \cdot\nabla j_\varepsilon(u_{h,i})\dx
&= -\int_K \varphi \Delta j_\varepsilon(u_{h,i})\dx
\le -\int_K \varphi\, j'_\varepsilon(u_{h,i})\Delta u_{h,i} \dx
= \int_K \varphi\, j'_\varepsilon(u_{h,i})b^+(x)u_{h,i} \dx\\
&\le\int_K \varphi\, \big(j_\varepsilon(u_{h,i})+\varepsilon\big)b^+(x) \dx
\end{split}
\]
for any subdomain  with compact closure $K\subset \Omega$, and all bounded $0\le\varphi \in C^1_0(K)$. Passing to the limit as $\varepsilon\to 0$ proves that $u_{h,i}^+$ is a weak subsolution to  $-\Delta u_{h,i}^+\le  b^+(x)\,u_{h,i}^+$. A similar procedure can be applied to $u_{h,i}^-$, therefore all the formal calculations made above are justified.

\noindent$\bullet~$\textsc{Step 2. }\textit{$\LL^\infty$-bounds for the gradients. }Since $|u_{h,i}|=u_{h,i}^++u_{h,i}^-$ is a weak nonnegative subsolution to $-\Delta\big|u_{h,i}\big|\le \lambda(p\vee1)\,b_p\,|u_{h,i}|:=b(x)\,|u_{h,i}|$\,,  we can apply the upper bounds of Theorem \ref{Moser.Upper.b.thm} that read
\begin{equation}\label{Moser.Upper.uhi}
\|u_{h,i}\|_{\infty, R}\le \frac{K^{(3)}_2[b]}{{h_0^{\frac{d}{2}}}}\,\|u_{h,i}\|_{2,R+h_0}
\end{equation}
with $q=2$ and the expression of the constant obtained by letting $r\to\infty$, since $b(x)\in \LL^\infty(B_{R+h_0})$:
\begin{equation}\label{Const.Upper.uhi}
\begin{split}
K^{(3)}_2[b]
&=\left(\frac{2d^d}{2^d}\right)^{\frac{d^2}{8}}
         \left[16(d+2)+\frac{h_0^2}{R^2}
        +\left(\frac{d\,\mathcal{S}_2^2}{d-2}\right)^{\frac{d}{2}}
        \frac{4\, h_0^2}{d-2}\,|B_{R+h_0} |^{\frac{d-2}{d}}\|b\|_{\infty, R+h_0}^{\frac{d}{2}}
         \right]^{\frac{d}{4}}\\
&\le\left(\frac{2d^d}{2^d}\right)^{\frac{d^2}{8}}
         \left[16(d+2)+\frac{h_0^2}{R^2}
        +\left(\frac{d\,\mathcal{S}_2^2 (p\vee1)}{d-2}\right)^{\frac{d}{2}}
        \frac{4\, h_0^2}{d-2}\,|B_{R+h_0} |^{\frac{d-2}{d}}\,\big(\lambda\,b_p\big)^{\frac{d}{2}}
         \right]^{\frac{d}{4}}\\
\end{split}
\end{equation}
since
\[
\|b\|_{\infty, R+h_0}^{\frac{d}{2}}=\big(\lambda(p\vee1)\,\,b_p\big)^{\frac{d}{2}}\,.
\]
Next we observe that by inequality \eqref{incr.1} it follows that for any $\delta>0$ and any $|h|<\delta/(10\sqrt{d})$ we have
\[
\|u_{h,i}\|_{2,R+h_0}\le 5^{\frac{d}{2}}\|\partial_i u\|_{2,R+h_0+\delta}.
\]
Finally, since
\[
\frac{\|u_{h,i}\|_{s, R}}{|B_R|^{\frac{1}{s}}}\le \|u_{h,i}\|_{\infty, R}\le K
\]
holds for any $|h|\le h_0$ with $K$ that do not depend on $s$, then by remark (iii) above we have that
\[
\frac{\|\partial_i u\|_{s, R}}{|B_R|^{\frac{1}{s}}}\le K\,.
\]
Letting now $s\to \infty$ in the above expression gives $\|\partial_i u\|_{\infty, R}\le K$.
Therefore we have proven that
\begin{equation}\label{Upper.pi}
\|\partial_i u\|_{\infty, R}\le \frac{K^{(3)}_2[b]}{{h_0^{\frac{d}{2}}}}\,5^{\frac{d}{2}}\|\partial_i u\|_{2,R+h_0+\delta}\,,
\end{equation}
with $K^{(3)}_2[b]$ as in \eqref{Const.Upper.uhi} which implies
\begin{equation}\label{Upper.grad}
\|\nabla u\|_{\infty, R}\le \frac{K^{(3)}_2[b]}{{h_0^{\frac{d}{2}}}}\,5^{\frac{d}{2}}\|\nabla u\|_{2,R+h_0+\delta}\,,
\end{equation}

\noindent$\bullet~$\textsc{Step 3. }\textit{Energy inequalities. }We now need the energy inequalities \eqref{local.energy.identity.subsols} to estimate the $\LL^2$ norm of the gradient of $u$ in terms of $u$ itself.  We choose $\alpha=1$ there so that the choice $\delta=0$ is admissible.
\begin{equation}\label{local.energy.identity.subsols.grad}\begin{split}
\int_{B_{R+h_0+\delta}} \big|\nabla u\big|^2\varphi\dx &\le \int_\Omega \big|\nabla u\big|^2\varphi\dx
\le \lambda\int_\Omega u^{p+1}\varphi\dx
    +\frac{1}{2}\int_\Omega u^2\Delta\varphi \dx\\
&\le \lambda\int_{B_{R+h_0+2\delta}} u^{p+1}\dx
    +\frac{2d}{\delta^2}\int_{B_{R+h_0+2\delta}} u^2\dx\\
&\le \left(\lambda\,b_p +\frac{2d}{\delta^2} \right)\int_{B_{R+h_0+2\delta}} u^2\dx\\
&\le \left(\lambda\,b_p +\frac{2d}{\delta^2} \right)\|u\|_{2, R+h_0+2\delta}^2
\end{split}
\end{equation}
since we have used the fact that $u^{p-1}\le b_p$ for any $0\le p<p_s$ and the test function $\varphi$ of Lemma \ref{lem.test.funct} with the choice of balls $B_{R+h_0+\delta}\subset B_{R+h_0+2\delta}$.

\noindent$\bullet~$\textsc{Step 4. }Putting all the pieces together, we have obtained
\begin{equation}
\|\nabla u\|_{\infty, R}\le \frac{K^{(3)}_2[b]}{{h_0^{\frac{d}{2}}}}\,5^{\frac{d}{2}}\|\nabla u\|_{2,R+h_0+\delta}
\le\frac{K^{(3)}_2[b]}{{h_0^{\frac{d}{2}}}}\,5^{\frac{d}{2}}
\left(\lambda\,b_p +\frac{2d}{\delta^2} \right)^{\frac{1}{2}}\|u\|_{2, R+h_0+2\delta}.
\end{equation}
We finally choose $h_0=\delta>0$ and we let $R_\infty=R$\,, $R_0=R+h_0+2\delta=R+3\delta$, so that $\delta=(R_0-R_\infty)/3$ and we have obtained
\begin{equation}\label{Upper.grad.u}
\|\nabla u\|_{\infty, R_\infty}
\le K^{(3)}_2[b]\,\left(\frac{15}{R_0-R_\infty}\right)^{\frac{d}{2}}
\left(\lambda\,b_p +\frac{18d}{(R_0-R_\infty)^2} \right)^{\frac{1}{2}}\|u\|_{2, R_0}:=K[u]\|u\|_{2, R_0}
\end{equation}
where we recall that, with the above choices of $h_0, \delta$ we have (see \eqref{Const.Upper.uhi})
\begin{equation}\label{Const.Upper.uhi.3}
\begin{split}
K^{(3)}_2[b]
&\le\left(\frac{2d^d}{2^d}\right)^{\frac{d^2}{8}}
         \left[16(d+2)+\frac{(R_0-R_\infty)^2}{9R_\infty^2}
        +\left(\frac{d\,\mathcal{S}_2^2 (p\vee1)}{d-2}\right)^{\frac{d}{2}}
        \frac{4(R_0-R_\infty)^2}{9(d-2)}\,|B_{R_0} |^{\frac{d-2}{d}}\,\left(\lambda\,\,b_p\right)^{\frac{d}{2}}
         \right]^{\frac{d}{4}}.\\
\end{split}
\end{equation}
Finally we observe that $b_p$ can be bounded depending on the values of $p$ as follows:\\
\noindent(i) If $0\le p<1$ we can use the absolute bounds \eqref{Absolute.lower.p<1} to get
\begin{equation}
b_p=\frac{1}{\inf\limits_{B_{R_0}}u^{1-p}}\le \mathcal{H}_p^{1-p}\frac{8R_0^d}{\lambda(R_0-R_\infty)^2 R_\infty^d},
\end{equation}
the constant $\mathcal{H}_p$ being given in this case by \eqref{Harnack.sublin.const}.

\noindent(ii) If $p=1$ then $b_p=1$\,.

\noindent(iii) If $1< p<p_c$ we can use the absolute bounds \eqref{Absolute.upper.p<pc}
\begin{equation}
b_p=\sup_{x\in B_R(x_0)}u^{p-1}(x)\le \mathcal{H}^{p-1}_p \frac{8R_0^d}{\lambda(R_0-R_\infty)^2 R_\infty^d},
\end{equation}
the constant $\mathcal{H}_p$ being given in this case by \eqref{Harnack.const.p<pc}.

\noindent(iv) If $p_c\le p<p_s$, we just leave $b_p=\|u\|_{\infty, R_0}^{p-1}$\,.\qed
\medskip

\noindent When $1<p<p_c$ we have local absolute bounds for the gradients, which seem to be new.
\begin{thm}[Local absolute bounds for the gradient when $1<p<p_c$]\label{thm.local.abs.grad}
Let $\Omega\subseteq\RR^d$ and let $\lambda>0$. Let $u$ be a local nonnegative weak solution to $-\Delta u = \lambda u^p$ in $B_{R_0}\subseteq\Omega$, with $1< p<p_c=d/(d-2)$. Then for any $0<R_\infty<R_0$ we have
\begin{equation}\label{Upper.grad.u.thm.abs}
\|\nabla u\|_{\infty, R_\infty}
\le K
\end{equation}
where
\begin{equation}\label{Const.Upper.uhi.4}
\begin{split}
K&= \left(\frac{d^d}{2^{d-1}}\right)^{\frac{d^2}{8}}\frac{\left(15\right)^{\frac{d}{2}}\mathcal{H}_p\omega_d^{\frac12}R_\infty^{\frac d2}}{(R_0-R_\infty)^{1+\frac d2+\frac2{p-1}}}
\left[\dfrac{8R_0^d\,\mathcal{H}_p^{p-1}}{R_\infty^d} +18d \right]^{\frac{1}{2}}  \left(\frac{8R_0^d}{\lambda R_\infty^d}\right)^{\frac{1}{p-1}}\\
&\times \left[16(d+2)+\frac{(R_0-R_\infty)^2}{9R_\infty^2}
        +\left(\frac{d\,\mathcal{S}_2^2 p}{d-2}\right)^{\frac{d}{2}}
        \frac{2^{2+\frac32 d}\,\omega_d^{\frac{(d-2)}{d}}\,R_0^{\frac{d^2}2+(d-2)}}{9(d-2)(R_0-R_\infty)^{2(d-1)}R_\infty^{\frac{d^2}2}}\,\mathcal{H}_p^{\frac{d(p-1)}2}
         \right]^{\frac{d}{4}}\\
\end{split}
\end{equation}
where the constant $\mathcal{H}_p$ is given by \eqref{Harnack.const.p<pc} and depends on $R_0,R_\infty$ as well.
\end{thm}

\section{Table of results}
Let us resume the main results of this paper: recall that $d\ge 3$ and
\[
p_c=\frac{d}{d-2}\,,\quad p_s=\frac{d+2}{d-2}\,,\quad \overline{q}=\frac{d(p-1)_+}{2}\,,\quad q_0=\frac{2^{{\frac{d-3}{2}}}}{d\omega_d^2[\ee(d-1)+\varepsilon]}, \quad \forall \varepsilon>0\,.
\]
\noindent \small
\begin{tabular}{ | c | c | c | c | c | c | c |}
 \hline
    &  \bf Upper I & \bf Upper II & \bf Lower & \bf Harnack & \bf Absolute & \bf Gradient \\
 \hline
   $0\le p<1$ & $0<q \to\infty$ & $q_0>0\,,\; r>0$& $0<q<q_0$ & $\mathcal{H}_p$ & lower & upper\\
              & Thm. \ref{thm.local.upper} & Thm. \ref{Thm.local.upper.b.p}& Thm. \ref{thm.local.lower} & Thm. \ref{Harnack.p<1}& Thm. \ref{Thm.Absolute.0.pc} & Thm. \ref{thm.local.upper.grad}\\
 \hline
   $p=1$ & $0<q \to\infty$  & $q_0>0$\,,\; $b\in\LL^r$\,,\;$r>\frac{d}{2}$ & $0<q<q_0$ & $\mathcal{H}_1$ & No  &upper\\
         & Thm. \ref{thm.local.upper}  &Thm. \ref{Thm.local.upper.b.DG} & Thm. \ref{thm.local.lower} &Thm. \ref{Harnack.p<1} & & Thm. \ref{thm.local.upper.grad} \\
 \hline
    $1< p<p_c$ & $\overline{q}<q\to \infty$ & $q_0>0$\,,\;  $b=\lambda u^{p-1}\in\LL^r$\,,\; $r>\overline{q}$& $\overline{q}<q<p_c$ & $\mathcal{H}_p$ & upper & absolute \\
               & Thm. \ref{thm.local.upper}  & Thm. \ref{Thm.local.upper.b.p}& Thm. \ref{thm.local.lower.pc} & Thm. \ref{Harnack.p<pc} & Thm. \ref{Thm.Absolute.0.pc} & Thm. \ref{thm.local.abs.grad}\\
 \hline
    $p_c< p<p_s$ & $\overline{q}<q\to \infty$ & $q_0>0$\,,\;  $b=\lambda u^{p-1}\in\LL^r$\,,\; $r>\overline{q}$& $0<q<q_0$ & $\mathcal{H}_p[u]$& No & upper\\
                 & Thm. \ref{thm.local.upper}  &Thm. \ref{Thm.local.upper.b.p} &Thm. \ref{thm.local.lower}  & Thm. \ref{Thm.Harnack.ps} & & Thm. \ref{thm.local.upper.grad} \\
 \hline
 \end{tabular}

\

\normalsize \noindent Recall the bounds:\\[3mm]
$
\mbox{\bf Upper I}\qquad\|u\|_{\LL^\infty(B_{R_\infty})} \, \|u\|_{\LL^{p-1}(B_{R_\infty})}^{\mu(p-1)_+}
    \le I_{\infty,q} \;\frac{\|u\|_{\LL^{q}(B_{R_0})}^{1+\mu(p-1)_+}}{|B_{R_0}|^{\frac{1}{q}}}\qquad\mbox{with }\quad\mu=\frac{d}{2q-d(p-1)_+}
$
\\[3mm]
$
\mbox{\bf Upper II}\qquad
\|u\|_{\infty,R_\infty}\le
    \frac{A_{q_0}^{(1)}}{(R-R_\infty)^{\frac{d}{q_0}}}
    \left[A_{q_0}^{(2)}+A_{q_0}^{(3)}\|b\|_{\LL^r({B_{R_0} })}^{\frac{rd}{2r-d}}\right]^{\frac{d}{2q_0}}\;\|u\|_{q_0,R}
$
\\[3mm]
$
\mbox{\bf Lower}\qquad\inf_{x\in B_{R_\infty}}u(x) =\|u\|_{\LL^{-\infty}(B_{R_\infty})}
    \ge I_{-\infty,q} \frac{\|u\|_{\LL^{q}(B_{R_0})}}{|B_{R_0}|^{\frac{1}{q}}}.
$
\\[3mm]
$
\mbox{\bf Harnack}\qquad\sup_{x\in B_{R_\infty}}u(x)\le \mathcal{H}_p[u]\inf_{x\in B_{R_\infty}}u(x)
$ \\[3mm]
where $\mathcal{H}_p[u]$ depends on $u$ only when $p_c\le p<p_s$ through some local norms as follows\\[3mm]
$
\mathcal{H}_p[u]=\mathcal{H}_p[u](d,\overline{q},\underline{q},\varepsilon, R_0, R_\infty)
= \frac{ I_{\infty,\overline{q}}}{I_{-\infty,\underline{q}}} \;
    \left(\frac{\left(\fint_{B_{R_0}} u^{q}\dx\right)^{\frac{(p-1)_+}{q}}}{   {\fint_{B_{R_\infty}}u^{(p-1)_+}\dx}}\right)^{\frac{d}{{2\overline{q}-d(p-1)_+}}}
    \frac{\left(\fint_{B_{R_0}} u^{\overline{q}}\dx\right)^{\frac{1}{\overline{q}}} }{\left(\fint_{B_{R_0}} u^{\underline{q}}\dx\right)^{\frac{1}{\underline{q}}}}\,.
$\\
whereas $\mathcal{H}_p[u]$ can be taken to be independent of $u$ if $p\in[0,p_c)$, see \eqref{Harnack.sublin.const}, \eqref{Harnack.const.p<pc}.
\\[3mm]
$
\mbox{\bf Gradient}\qquad\|\nabla u\|_{\infty, R_\infty}
\le K[u]\,\|u\|_{2, R_0}\,.
$

\section{Appendix. Numerical Identities and Inequalities}

\noindent\textbf{Sum of some series.}
\[\begin{split}
&\sum\limits_{j=1}^\infty
        \left(\frac{2}{2^*}\right)^{j}=\sum\limits_{j=1}^\infty
        \left(\frac{d-2}{d}\right)^{j}=\frac{d-2}{2}\,\qquad
        \sum\limits_{j=1}^\infty j\left(\frac{2}{2^*}\right)^j=\frac{d(d-2)}{4}\\
        &\sum\limits_{j=k+1}^\infty \left(\frac{2}{2^*}\right)^j=\frac{d}{2}\left(\frac{2}{2^*}\right)^{k+1}\,\qquad
        \sum\limits_{j=1}^k \left(\frac{2}{2^*}\right)^j=\frac{d-2}{2}\left(\frac{2}{2^*}\right)^k
\end{split}
\]
since for any $0<s<1$ we have, for any $0\le N\in\mathbb{N}$
\[
\sum_{j=0}^\infty j^N\,s^j =\underbrace{s\frac{\rd}{\rd s}\left[s\frac{\rd}{\rd s}\left(\dots s\frac{\rd}{\rd s}\left(\frac{1}{1-s}\right)\dots\right)\right]}_{\mbox{N-times}}=\left[s\frac{\rd}{\rd s}\right]^{(N)}\left(\frac{1}{1-s}\right)\mbox{.\,\qed}
\]

\noindent\textbf{Stirling's formula:}
\begin{equation}\label{Stirling}
n!=\sqrt{2\pi\, n}\,\left[\frac{n}{e}\right]^n\,\ee^{\alpha_n}\qquad\mbox{with}\qquad \frac{1}{12n+1}\le \alpha_n \le\frac{1}{12n}.
\end{equation}
We recall that
\[
\omega_d=\frac{\pi^{d/2}}{\Gamma\left(1+\frac{d}{2}\right)}
\sim\frac{\big(2\,\ee\sqrt{\pi}\big)^d}{d^d\,\ee^{\alpha_d}\,\sqrt{d\,\pi}}
\qquad\mbox{with}\qquad \frac{1}{6d+1}\le \alpha_d \le\frac{1}{6d}.
\]\normalcolor

\begin{lem}
The following inequality holds for any $a,b\ge 0$
\begin{equation}\label{num.ineq.abp}
(a-b)(a^p-b^p)\le (p\vee 1)\,\max\big\{a^{p-1}, b^{p-1}\big\}(a-b)^2\,,\qquad\mbox{and for any $p\ge0$\,.}
\end{equation}
Moreover the following inequality holds for any $a,b\ge 0$ and $p\ge 1$:
\begin{equation}\label{num.ineq.abp2}
a^p-b^p\ge p\, b^{p-1}(a-b).
\end{equation}
\end{lem}
\noindent {\sl Proof.~} If $a\ge b$ the validity of \eqref{num.ineq.abp} is equivalent, setting $x=\frac ba$, to the validity of $(1-x)(1-x^p)\le p(1-x)^2$ for all $x\in [0,1]$, that is to $1-x^p\le p(1-x)$ for all $x\in[0,1]$, which does in fact hold if $p\ge1$ by the concavity of $g(x):=1-x^p$, since the line $h(x):=p(1-x)$ is the tangent to $g$ at $x=1$. The case $a<b$ follows as well by interchanging the role of $a$ and $b$. The case $0<p<1$ can be proven analogously: if fact the stated inequality is equivalent to $1-x^p\le1-x$ for any $x\in[0,1]$, which holds true by the convexity of $h(x)=1-x^p$ for any $p\in(0,1)$.

The second inequality \eqref{num.ineq.abp2} follows by the inequality $x^p-1\ge p(x-1)$ for all $x\ge 0$
which is valid since $x^p-1$ is convex so that its graph lies above its tangent at $x=1$.\qed

\normalcolor

\

\noindent {\large \sc Acknowledgment}

\noindent  The first and third author have been partially funded by Projects MTM2008-06326-C02-01 and MTM2011-24696 (Spain).

\vskip1cm


\small
\bibliographystyle{amsplain} 

\end{document}